\documentclass[12pt]{article}
\pdfoutput=1    
\usepackage{amscd,amssymb,stmaryrd}
\usepackage{amsmath}
\usepackage{amscd,amssymb,stmaryrd}
\usepackage{graphicx}
\usepackage{subfigure,fancybox}
\usepackage{amscd,amstext}
\usepackage{amsmath}
\usepackage{mathtools}
\usepackage{setspace}
\usepackage{amssymb}
\usepackage{slashed}
\usepackage{pstricks,pst-node,pst-plot,pst-coil,pstcol}
\usepackage{graphicx}
\usepackage{ntheorem}
\usepackage{upgreek}
\usepackage{enumerate}
\usepackage{graphicx}
\usepackage{cite}
\usepackage{tikz}
\usepackage{authblk}
\usepackage{txfonts}
\usepackage[colorlinks,linkcolor=black,anchorcolor=black,citecolor=black,CJKbookmarks=True]{hyperref}
\theorembodyfont{\normalfont}
\newtheorem{thm}{Theorem}[section]
\newtheorem{lem}{Lemma}[section]
\newtheorem*{pf}{PROOF:}
\newtheorem{prop}{Proposition}[section]
\newtheorem{defi}{Definition}[section]

\newtheorem{remark}{Remark}[section]
\newtheorem{cor}{Corollary}[section]

\numberwithin{equation}{section}
\title{\large Formation of shifted shock for the 3D compressible Euler equations with damping}
\author[1]{\normalsize Chen Zhendong\footnote{Acknowledgement:
This is part of the Ph. D thesis of the author written under the supervision of Professor Zhouping Xin at the Institute of Mathematical Sciences of The Chinese University of Hong Kong. The research is supported in part by
the Zheng Ge Ru Foundation and by the Hong Kong RGC
Earmarked Research  Grants:  CUHK-14305315, CUHK-14300917, and CUHK-14302819.}}

\affil[1]{Institute of Mathematical Science, The Chinese University of Hong Kong, Shatin, NT, Hong Kong.}
\date{}
\begin{document}           
\def\ra{\rightarrow}\def\Ra{\Rightarrow}\def\lra{\Longleftrightarrow}\def\fai{\varphi}
\def\s{\slashed}\def\Ga{\Gamma}\def\til{\tilde}
\def\de{\Delta}\def\fe{\phi}\def\sg{\slashed{g}}\def\la{\lambda}\def\R{\mathbb{R}}\def\m{\mathbb}
\def\a{\alpha}\def\p{\rho}\def\ga{\gamma}\def\lam{\lambda}\def\ta{\theta}\def\sna{\slashed{\nabla}}
\def\pa{\partial}\def\be{\beta}\def\da{\delta}\def\ep{\epsilon}\def\dc{\underline{\chi}}
\def\si{\sigma}\def\Si{\Sigma}\def\wi{\widetilde}\def\wih{\widehat}\def\beeq{\begin{eqnarray*}}
\def\eeq{\end{eqnarray*}}\def\na{\nabla}\def\lie{{\mathcal{L}\mkern-9mu/}}\def\rie{\mathcal{R}}
\def\ud{\underline}\def\les{\lesssim}\def\ka{\kappa}
\def\dl{\underline{L}}\def\du{\underline{u}}
\def\hs{\hspace*{0.5cm}}\def\bee{\begin{equation}}\def\ee{\end{equation}}\def\been{\begin{enumerate}}
\def\een{\end{enumerate}}\def\bes{\begin{split}}\def\zpi{\prescript{(Z)}{}{\pi}}\def\szpi{\prescript{(Z)}{}{\slashed{\pi}}}
\def\ees{\end{split}}\def\hra{\hookrightarrow}\def\udx{\ud{\xi}_{A}}\def\ude{\ud{\eta}_{A}}
\def\tpi{\prescript{(T)}{}{\pi}}\def\lpi{\prescript{(L)}{}{\pi}}\def\qpi{\prescript{(Q)}{}{\pi}}
\def\stpi{\prescript{(T)}{}{\slashed{\pi}}}\def\sqpi{\prescript{(Q)}{}{\slashed{\pi}}}
\def\zgpi{\prescript{(Z)}{}{\wi{\pi}}}\def\szgpi{\prescript{(Z)}{}{\wi{\slashed{\pi}}}}
\def\tgpi{\prescript{(T)}{}{\wi{\pi}}}\def\lgpi{\prescript{(L)}{}{\wi{\pi}}}\def\qgpi{\prescript{(Q)}{}{\wi{\pi}}}
\def\stgpi{\prescript{(T)}{}{\wi{\slashed{\pi}}}}\def\sqgpi{\prescript{(Q)}{}{\wi{\slashed{\pi}}}}
\def\pre #1 #2 #3{\prescript{#1}{#2}{#3}}\def\spi{\slashed{\pi}}\def\rgpi{\prescript{(R_{i})}{}{\wi{\pi}}}
\def\rpi{\prescript{(R_{i})}{}{\pi}}\def\srpi{\prescript{(R_{i})}{}{\s{\pi}}}
\def\pp{\uprho}\def\vw{\upvarpi}\def\srgpi{\pre {(R_{i})} {} {\wi{\s{\pi}}}}
\def\supnormda{L^{\infty}(\Si_{t}^{\tilde{\da}})}\def\supnormu{L^{\infty}(\Si_{t}^{u})}
\def\normda2{L^{2}(\Si_{t}^{\tilde{\da}})}\def\normu2{L^{2}(\Si_{t}^{u})}
\def\les{\lesssim}\def\sdfai{|\s{d}\fai|}\def\sd{\s{d}}\def\Llra{\Longleftrightarrow}
\def\oli{\overline}\def\vta{\vartheta}\def\ol{\overline}\def\mari{\mathring}
\def\ea+2{E_{\leq|\a|+2}}\def\dea+2{\ud{E}_{\leq|\a|+2}}\def\fa+2{F_{\leq|\a|+2}}\def\dfa+2{\ud{F}_{\leq|\a|+2}}
\def\tea+2{\wi{E}_{\leq|\a|+2}}\def\tdea+2{\ud{\wi{E}}_{\leq|\a|+2}}\def\tfa+2{\wi{F}_{\leq|\a|+2}}\def\tdfa+2{\ud{\wi{F}}_{\leq|\a|+2}}
\def\ntea{\wi{E}_{\leq|\a|+1}}\def\ntdea{\ud{\wi{E}}_{\leq|\a|+1}}\def\ntfa{\wi{F}_{\leq|\a|+1}}\def\ntdfa{\ud{\wi{F}}_{\leq|\a|+1}}
\def\ba2{b_{|\a|+2}}\def\Ka{K_{\leq|\a|+2}}\def\mumba{\mu_m^{-2\ba2}}\def\nba{b_{|\a|+1}}
\def\endpf{\hfill\raisebox{-0.3cm}{\rule{0.1mm}{3mm}\rule{0.5mm}{0.1mm}\rule[-0.6mm]{2.3mm}{0.7mm}\raisebox{3mm}[0pt][0pt]{\makebox[0pt][r]{\rule{2.9mm}{0.1mm}}}\rule[-0.6mm]{0.1mm}{3.5mm}\rule[-0.6mm]{0.5mm}{3mm}} \vspace{0.5cm}\\}

\pagestyle{myheadings} \thispagestyle{empty} \markright{}

\maketitle
\begin{abstract}
In this paper, we show the shock formation of the solutions to the 3-dimensional (3D) compressible isentropic and irrotational Euler equations with damping for the initial short pulse data which was first introduced by D.Christodoulou\cite{christodoulou2007}. Due to the damping effect, the largeness of the initial data is necessary for the shock formation and we will work on the class of large data (in energy sense). Similar to the undamped case, the formation of shock is characterized by the collapse of the characteristic hypersurfaces and the vanishing of the inverse foliation density function $\mu$, at which the first derivatives of the velocity and the density blow up. However, the damping effect changes the asymptotic behavior of the inverse foliation density function $\mu$ and then shifts the time of shock formation compared with the undamped case. The methods in the paper can also be extended to a class of $3D$ quasilinear wave equations for the short pulse initial data.\\
\par\textbf{Keywords: } Compressible Euler equations with damping, delayed shock formation, inverse foliation density function, characteristic null-hypersurfaces
\end{abstract}

\tableofcontents
\section{\textbf{Introduction}}\label{section1}
\hs In this paper, we will consider the following compressible Euler equations with damping in $\mathbb{R}^3$:
\bee\label{cee}
\left\{\bes
&\dfrac{\pa}{\pa t}\p+\nabla\cdot(\p v)=0,\\
&\dfrac{\pa}{\pa t}(\p v)+\nabla\left(\p v\otimes v+p I_3\right)+
a\p v=0,\\
&\dfrac{\pa}{\pa t}S+v\cdot\nabla S=0,
\end{split}\right.
\ee
where $\p$, $v=(v^1,v^2,v^3)$, $p$, $S$ represent the density, the velocity, the pressure and the entropy of the flow, respectively. $I_3$ is a $3\times 3 $ identical matrix and $a$ is the damping constant. These equations describe the motion of a perfect fluid which follow from the conservation of mass, momentum and energy respectively. The enthalpy $h$ of a thermodynamic system is defined as the sum of its internal energy and the product of its pressure and volume:
 \bee
 h=e+pV.
 \ee
 Here $e$ and $V=\frac{1}{\p}$ represent the internal energy and the volume of the system, respectively. It follows from the relation of thermodynamics that
 \bee\label{enthalpy}
 dh=TdS+Vdp,
 \ee
 where $T$ is the absolute temperature. The sound speed $\eta$ is defined as $\eta=\sqrt{\dfrac{\pa p}{\pa\p}}$. The equation of state is given by $p=p(\p,S)$ and we assume that:
   \bee\label{nonlinearcondition}
   p\ \text{is not linear in}\ \dfrac{1}{\p}.
   \ee
  This assumption will be clarified later.
\subsection{\textbf{Brief review of known results}}

\hs If the damping term vanishes, then the system \eqref{cee} is returned to the classical compressible Euler system, which is a prototype of hyperbolic systems of conservation laws:
    \bee\label{conservationlaws}
    \bes
    &\pa_tU+div(F(U))=\pa_tU+A(U)\nabla U=0,\\
    &U(x,0)=U_0(x),
    \end{split}
    \ee
    where $U: \mathbb{R}^m\to \mathbb{R}^n$. Understanding the shock formation and development to \eqref{conservationlaws} is fundamental since the system will develop singularities in general in finite time for initial smooth small data. B.Riemann was the first one to study the nonlinear effects to the $1$D isentropic Euler equations. In\cite{Riemann1860},
    Riemann proved that the wave compression leads to the shock waves and the wave expansion leads to the rarefaction waves.
    Later, P.Lax generalized Riemann's result into the $2\times2$ system (i.e. $n=2$) in one space dimension in\cite{Lax1987HyperbolicSO}. Lax used the Riemann invariants $(\tilde{u},\tilde{v})$, that is, the gradient of each invariant is proportional to the left eigenvector of the coefficient matrix. Then, Lax was able to give a sufficient and necessary condition for the system to admit a shock or not. However, for general system with $n>2$, Lax's method may fail since the system may not admit a coordinate system of Riemann invariants. In 1974, F.John in\cite{doi:10.1002/cpa.3160270307} achieved a remarkable result for the shock formation to general $n\times n$ hyperbolic systems of conservation laws in one space dimension. He considered the following Cauchy problem:
        \bee\label{john}
        \bes
        &\pa_tu+A(u)\pa_xu=0,\\
        &u(x,0)=f(x),
        \end{split}
        \ee
        where $A(u)$ is an $n\times n$ smooth coefficient matrix with real distinct eigenvalues. F.John showed that if the system \eqref{john} is genuinely nonlinear and the initial data $f(x)$ is small with compact support, then any classical solution to \eqref{john} must form shock in finite time, which means that the first order derivatives of $u$ blow up in finite time while itself remains bounded. However, F.John's result failed to apply to the Euler equations in one space dimension since the entropy wave is linearly degenerate, and later T.P.Liu generalized F.John's result in\cite{LIU197992}. T.P.Liu proved that if each characteristic family for the system \eqref{john} is either genuinely nonlinear or linearly degenerate, then under certain condition, any classical solution to \eqref{john} must develop shock in finite time.\\
         \hs In multi-dimensional case, the singularity formation problem for \eqref{conservationlaws} is much more complicated than in $1$D case. The first general result for the singularity formation to the compressible Euler equations in three spatial dimensions was obtained by Sideris\cite{Sideris1985} for polytropic gases. In particular, he exhibited an open set of small initial data for which the corresponding solutions cease to be $C^1$ in finite time by using dissipative energy estimates. However, his results did not provide any information for the mechanism of the breakdown of the solutions. Later, Alinhac studied the two-dimensional compressible isentropic Euler equations with radial symmetry in\cite{Alinhac1993}. He showed that a large class of small radially symmetric data leading to the finite time blow up of the solutions and gave a precise estimate for the blow up time. Later, in a series of works \cite{alinhac1999blowup,alinhac1999blowup1,alinhac2001blowup,alinhac2001blowup2}, Alinhac proved the shock formation to the $2D$ quasilinear wave equations. He constructed a class of initial data which lead to the break down of the solutions in finite time and gave precise estimates of the solutions up to the first singularity. \\
        \hs For the system\eqref{cee}, there are many studies both in $1D$ and multi-dimensions. For the one-dimensional Euler equations with damping, the global existence of smooth solutions with small data was proved by Nishida\cite{1978NonlinearHE} and Slemrod\cite{slemrod} showed that for small data ($L^{\infty}$ sense), the $1D$ Euler equations with damping admit a global smooth solution while for large data, the equations can develop a shock in finite time. Later, these results were generalized by many authors, see \cite{liu1992,NISHIHARA2000191,MARCATI1990,MARCATI2000,MARCATI2005,NISHIHARA2000} and the references therein.\\
\hs In multi-dimensional case, the global existence and $L^2$ estimates for the solutions to the $3D$ isentropic Euler system with damping was obtained by Kawashima\cite{Kawashima1984} and the long time behavior of the solutions was obtained and generalized by many authors, see \cite{WANG2001410,Kawashima2004,PAN2009581,TAN20121546,sideris2003long} and the references therein.
\subsubsection{\textbf{Christodoulou' s theory of shock formation}}
\hs A major breakthrough in understanding the shock mechanism for the hyperbolic systems\eqref{conservationlaws} in multi-dimensions has been made by Christodoulou in a series of works\cite{christodoulou2007,christodoulou2008,christodoulou2014compressible}. In\cite{christodoulou2014compressible}, the classical, non-relativistic, isentropic compressible Euler¡¯s equations in three spatial dimensions with initial irrotational data were studied. Starting from the short pulse data, the authors gave a detailed analysis of the solutions near the singularity (shock) and showed a complete geometric structure for the shock development. Note that under the isentropic and irrotational assumptions, the 3D Euler system can be rewritten as the following quasilinear wave equation:
\begin{equation}\label{christodoulou}
\Box_{g(\fai)}\fai=0,
\end{equation}
 where $\fai$ represent the first order derivatives of the potential function $\fe$ which is defined as $v=-\nabla\fe$ and $g$ is a Lorentzian metric depending on $\fai$.
 Instead of studying the equation\eqref{christodoulou} in the Cartesian coordinates, the authors constructed a system of geometric coordinates (the acoustical coordinates). In the acoustical coordinates, the solution is regular and the shock formation corresponds to the transformation between two coordinates degenerating at which the first order derivatives of $\fai$ blow up in the Cartesian coordinates. In \cite{christodoulou2014compressible}, the authors achieved their results through the following two steps (see also the following picture).
 
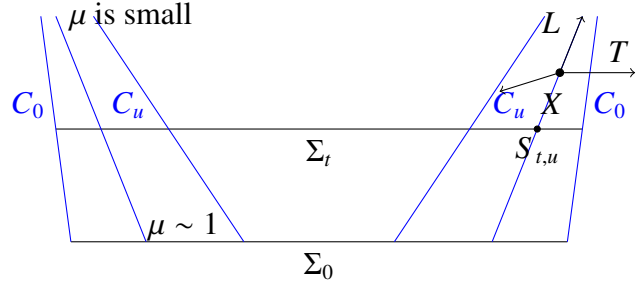
\begin{center}
\begin{tikzpicture}
\draw (-3.3,0)--node[below] {$\Si_0$}(3.3,0);
\draw (-3.5,1.5)--node[below] {$\Si_t$}(3.5,1.5);
\draw[blue] (1,0)--(3,3);
\draw[blue] (2.3,0)--node[above left] {$C_u$}(3.5,3);
\node at (-2.5,3) {$\mu$ is small};
\node at (-1.8,0.2) {$\mu\sim 1$};
\filldraw (2.9,1.5) circle (.04);
\node at (2.9,1.2) {$S_{t,u}$};
\draw[blue] (3.3,0)--node[above right] {$C_0$}(3.7,3);
\draw[blue] (-1,0)--(-3,3);
\draw[blue] (-2.3,0)--node[above right] {$C_u$}(-3.5,3);
\draw[blue] (-3.3,0)--node[above left] {$C_0$}(-3.7,3);
\filldraw (3.2,2.25) circle (.05);
\draw[->] (3.2,2.25)--node[above left] {$L$}(3.5,3);
\draw[->] (3.2,2.25)--node[above right] {$T$}(4.2,2.25);
\draw[->] (3.2,2.25)--node[below right] {$X$}(2.4,2);
\end{tikzpicture}
\end{center}
\begin{remark}
In the above picture, $\Si_{t'}$ represents the hypersurface $\{t=t'\}$ in the space time and $C_u$ represent the characteristic hypersurfaces where $u\in[0,\ep_0]$ with $\ep_0$ being a small constant. Their intersection is the "sphere" $S_{t,u}:=\Si_{t}\cap C_u$ which is a 2-dimensional submanifold (see also the definition (2.22)).
\end{remark}
 \begin{itemize}
 \item The first step is the geometric formulation. Given a solution $\fai$ and the a prior estimates for $\fai$ and its derivatives, they constructed the following objects:
       \been[(1)]
       \item the acoustic ekional function $u$ defined as $g^{\a\be}\pa_{\a}u\pa_{\be}u=0$ whose level sets are the characteristic null-hypersurfaces $C_u$. The ekional function $u$ forms a component of the acoustical coordinates $(t,u,\ta_1,\ta_2)$ where $(\ta_1,\ta_2)$ are the coordinates on $S_{t,u}$;
       \item the inverse foliation density function $\mu$ whose reciprocal measures the density of the hypersurfaces $C_u$. At the blow up point, $\mu\to 0$ and the characteristic hypersurfaces collapse;
       \item as long as $\mu>0$, the authors could construct the frame $\{L,T,X\}$ which are equivalent to $\{\frac{\pa}{\pa x^{\a}} \}_{\a=0,1,2,3}$, where $L$ represents the null vector field which equals $\frac{\pa}{\pa t}$ in the acoustical coordinates while $T$ and $X_i,i=1,2$ (they are all spacelike) represent the vector fields of "radial" derivative and the angular derivatives in the acoustical coordinates respectively.
       \een
 \item The second step is the top order energy estimates, which are based on
       \been[(1)]
       \item deriving fundamental energy estimates to the wave equation $\Box_g\fai=F$ by using the multiplier method where $F$ represent the general inhomogeneous terms;
       \item commuting a list of commutators with $\Box_g\fai=F$ and applying fundamental energy estimates to obtain the top order energy estimates which will be stated in details as follows.
       \een
 \end{itemize}
In the picture above, initially, $\mu\sim1$ and the characteristic hypersurfaces $C_u$ are not compressive while as time increasing, $\mu$ tends to $0$ and $C_u$ become dense so that the shock formation is precisely captured by the vanishing of $\mu$. The vanishing of $\mu$ and the blow up of the first order derivatives of $\fai$ are related as follows. Note the following fact on the transport equation for $\mu$
\bee
L\mu\sim -2T\fai+\text{small terms},
\ee
and when $\mu$ is small, then $L\mu\leq-C$ for $C>0$. It can be shown that $T\fai:=\mu\hat{T}\fai\sim 1$ even $\mu$ is small. As a consequence, $\hat{T}\fai$ blows up like $\frac{1}{\mu}$. This leads to the difficulty that when concerning the higher order energy estimates involving multi-derivatives of $T$ such as $\int |T^m\fai|^2$, the energies may blow up like $\mu^{-n}$. To deal with this difficulty, the authors then derived the exact blow up rate of the higher order energies by considering the weighted energy $\int \mu^{b_{\ast}} T^m\fai$ where the power of $\mu$ should be carefully chosen. To close the top order energy estimates, one has to commute $\{L,T,X\}$ with the transport equations for $\mu$ and $\chi$ to derive the top order estimates for $\mu$ and $\chi$. 
This part is highly involved and requires the regularization of transport equations for $\mu$ and $\chi$, deriving the bounds for $\int \mu^{-b}$ and the bounds for the top order spatial derivatives of $\mu$ and $\chi$. Based on the estimates for $\mu$ and $\chi$, the authors showed the following weighted energy estimates:
\bee\label{christodoulou2}
E_{modified}^{N_{top}}\les \text{initial energy},
\ee
where $E_{modified}^{N_{top}}$ are the top order weighted energies (contain the power of $\mu$). The estimates \eqref{christodoulou2} are not enough to provide the pointwise bounds for the derivatives of $\fai$ since the weights in the energies involve the power of $\mu$ which may vanish as a shock forms. However, by choosing appropriate weights, i.e. $b$, one can lower the order of the derivatives of $\fai$ in the energies to eliminate the power of $\mu$ which leads to
\bee
E^{N_{top}-M}\les\text{initial energy},
\ee
where $E^{N_{top}-M}$ are the desired energies which do not contain the power of $\mu$. By Sobolev embedding, one can show the $L^{\infty}$ bounds for the derivatives of $\fai$ and thus close the a priori assumptions.\\
\hs Later, in\cite{Ontheformationofshocks}, Pin Yu and Shuang Miao applied Christodoulou's framework to the following quasilinear wave equation in three space dimensions:
    \bee\label{yupin}
    \left\{\begin{split}
    &-(1+3G''(0)(\pa_t\fe)^2)\pa^2_t\fe+\de\fe=0,\\
    &(\fe,\pa_t\fe)(-2,x)=(\da^{\frac{3}{2}}\fe_0(\frac{r-2}{\da},\ta),\da^{\frac{1}{2}}\fe_1(\frac{r-2}{\da},\ta)),
    \end{split}\right.
    \ee
    where $G''(0)$ is a non-zero constant, $\ta\in S^2$, $\da>0$ is sufficient small and $\fe_0(s,\ta),\fe_1(s,\ta)\in C_0^{\infty}((0,1]\times S^2)$. Note that the initial data they constructed is of large energy (in $H^s$ sense) since for small initial data, \eqref{yupin} admits a global solution. Then, they proved the break down of the solution is due to the collapse of the characteristic hypersurfaces and gave a sufficient condition leading to the shock formation before $t=-1$.\\
\hs Recently, the $2D$ Euler equations with rotations were studied by J.Speck and J.Luk in\cite{luk2018shock} by applying the same framework.
    They discovered the Euler equations can be rewritten as a system of wave equations with the inhomogeneous terms as follows:
    \bee\label{speck2}
    \begin{split}
    \mu\Box_gv^i&=\mu\mathcal{D}^i+\text{terms of easy to deal with},\\
    \mu\Box_g\varrho&=\mu \mathcal{D},
    \end{split}
    \ee
    where $\varrho$ is the logarithmic density and $\mathcal{D},\mathcal{D}^i$ are the null forms with respect to $g$. Roughly speaking, the null forms $\mathcal{D}$ have the form $\mathcal{D}(f,g)=g^{\a\be}\pa_{\a}f\pa_{\be}g$ for some functions $f,g$ and do not contain $\hat{T}$ derivative in the frame $\{L,T,X\}$ decomposition. Speck and Luk considered the initial plane symmetric data with the perturbation of short pulse and proved the shock formation is due to the collapse of the characteristic hypersurfaces, at which the first order derivatives of $v$ and $\varrho$ blow up while the vorticity remains bounded. Later, they also generalized the formation of shocks to $3D$ case in\cite{Specklukshockformation}.\\
\hs More recently, Buckmaster, Shkoller and Vicol studied the shock formation for the Euler system in multi-dimensional case in a series of works\cite{BSV2Disentropiceuler,BSV3Disentropiceuler,BSV3Dfulleuler,formationofunstableshock}. In\cite{BSV2Disentropiceuler}, they considered the $2D$ isentropic Euler equations under azimuthal symmetry (but different from $1D$ problem) with smooth initial data of finite energy and nontrivial vorticity. By using the modulated self-similar variables, they showed the point shock forms in finite time with explicitly computable blow-up time and location and obtained that the solutions near shock formation point are of cusp type. Later, in\cite{BSV3Disentropiceuler} and \cite{BSV3Dfulleuler}, they generalized the above results into $3D$ isentropic and non-isentropic Euler systems.\\ 
\hs In our work, we found that Christodoulou's framework can be applied to the $3D$ Euler system with damping. Our work is aiming to show the effect of damping on the shock formation. Since for small data, \eqref{cee} admits a global solution and we have to construct a class of short pulse data with large energy. It will be shown that for such data, damping may not be strong to prevent the shock formation. And if a shock forms, the time of shock formation will be shifted compared with the undamped case and the shift depends on the size of the initial data and the damping constant $a$. Although we adapt the Christodoulou's framework and J.Speck's work, their methods can not be adapted directly due to the following reasons:
\begin{itemize}
\item the Lorentian metric $g$ in their work depends only on the variations $\fai$ ($\Box_{g(\fai)}\fai=0$) while our's depends on both the potential function $\fe$ and its variations $\fai$, which will affect the structures of the two transport equations for $\mu$ and $\chi$. Precisely, for the transport equation of $\mu: L\mu=m+\mu e$, there will be one more term $aT\fe$ in $m$, which vanishes in the undamped case. Then, in order to obtain the estimate for $\mu$, it requires the estimate for $T\fe$ and the estimate for $\mu$ will help us to recover the assumption on $T\fe$.\\
     \hs In order to obtain the estimates for the high order derivatives of $\mu$, it requires one to derive the estimates for the high order derivatives of $T\fe$, which include $L^{m-1}T\fe$, $X^{m-1}T\fe$ and $T^m\fe$ $(m\geq 2)$. The estimates for the former two can be derived easily while for $T^m\fe$, they can not be derived trivially. 
    To deal with this difficulty, one has to commute $T^{m-1}$ with the equation $L\mu=m+\mu e$ which gives us $LT^{m-1}\mu=T^{m-1}m+eT^{m-1}\mu+\text{small terms}$. Then, one can derive the estimates for $T^{m-1}\mu$ by using an induction argument and the bootstrap assumptions on $T^m\fe$. The estimates for $T^{m-1}\mu$ will lead us to recover the a priori assumptions on $T^m\fe$. 
\item Christodoulou considered the following covariant wave equation:
      \bee\label{covariantexample}
      \Box_{g}\fai=f,
      \ee
      with $f\equiv0$. 
      In our work, due to the damping, the inhomogeneous terms $f$ contain two parts: $a\frac{\pa\fai}{\pa t}$ and $a\fai\de\fe$, which will blow up as $\mu\to 0$. To deal with the second term, one first notes that $\de\fe\sim \mu^{-1}T\fai+X\fai$ so that this term blows up like $\frac{1}{\mu}$, which means if one multiplies $\mu$ to this term, it can be bounded! Moreover, since $\fai$ is small, then one can obtain that $a\fai\de\fe$ is indeed small by multiplying $\mu$. For the first term, since $\frac{\pa}{\pa t}\sim L+\mu^{-1}T$, then one can also multiply $\mu$ to it but will obtain that $\mu\pa_t\fai\sim1 $, which is not small. This term plays a crucial role in the accurate estimates of $\mu$ and will illustrate the damping effect on the formation of shock.
\end{itemize}
 In this paper, we will adapt the major geometric framework invented by Christodoulou\cite{christodoulou2014compressible} (see also the presentation of J.Speck) and will modify the analysis therein. For definiteness, the construction of the short pulse initial data is similar to that of Christodoulou\cite{christodoulou2007} and Yu-Miao\cite{Ontheformationofshocks}.\\
\subsection{\textbf{Notations}}
Through the whole paper, the following notations will be used unless stated otherwise:
\begin{itemize}
\item Latin indices $\{i,j,k,l,\cdots\}$ take the values $1, 2, 3, $ Greek indices $\{\a,\be,\ga,\cdots\}$
take the values $ 0, 1, 2, 3$ and capital letter $\{A,B,C,\cdots\}$ take values $1,2.$ Repeated indices are meant to be summed.
\item The convention $f\les h$ means that there exists a universal positive constant $C$ such that $f\leq Ch$.
\item The notations l.o.ts (lower order terms) mean the terms are of lower order. Here, the order means the number of total derivatives acting on $\fai$ and we set $\fai$ to be order $0$. For example, one can rewrite $\pa^2\fai+\pa\fai$ as $\pa^2\fai+$l.o.ts. Let $O_b^{\leq a}$ be the terms of order $\leq a$ and with bound $\da^{\frac{b}{2}}$. 
\item For the metric $g_{\a\be}$, $g^{\a\be}$ means its inverse such that $g_{\a\be}g^{\be\ga}=\da^{\ga}_{\be}$ with $\da^{\ga}_{\be}$ being the Kronecker symbol.
\item The box operator $\Box_g:=g^{\a\be}D^2_{\a\be}$ denotes the covariant wave operator corresponding to the spacetime metric g and $\s{\de}:=\sg^{AB}\s{D}^2_{AB}$ denotes the covariant Laplacian corresponding to $\sg$ on $S_{t,u}$, where $D$, $\s{D}$ are the Levi-Civita connections corresponds to $g$, $\sg$ respectively. We also denote $\til{D}$ to be the Levi-Civita connection corresponds to $\til{g}$.
\item For a object $q$, $\s{q}$ means its restriction(projection) on $S_{t,u}$. In particular, $\s{div}$ represents the divergence operator on $S_{t,u}$ such that $\s{div}Y:=\s{D}_AY^A$ for any $S_{t,u}$ vector field $Y$. $\s{d}$, $\lie$ represent the restriction of the standard differentiation $d$ and the lie derivative $\mathcal{L}$ on $S_{t,u}$ respectively.
\item For $q$ being a $(0,2)$ tensor and $Y,Z$ being the vector fields, set the contraction as
      \bee
      q_{YZ}:=q_{\a\be}Y^{\a}Z^{\be},
      \ee
      and similar for the other type tensors.
\item For a spacetime vector field $V$, denote the decomposition of $V$ relative to the frame $\{L,T,X\}$ as
      \bee
      V=V^LL+V^TT+V^AX_A,
      \ee
      with $V^L,V^T$ and $V^A$ being functions.
\item For a $(0,2)$ tensor $\ta$, the following decomposition holds:
      \bee
      \ta=\hat{\ta}+\frac{1}{2}tr\ta\cdot g,
      \ee
      where $\hat{\ta}$ is the trace free part of $\ta$ and $tr\ta=g^{\a\be}\ta_{\a\be}$ is the trace of $\ta$ with respect to $g$.
\end{itemize}
\subsection{\textbf{Geometric blow-up for the Burgers' equation with damping}}
Consider the following Cauchy problem for $1D$ Burger's equation with damping
\bee\label{burgersdamp}
\left\{\bes
&\pa_t\fai+\fai\pa_x\fai=-a\fai,\\
&\fai(x,t=-1)=f(x),
\end{split}\right.
\ee
where $a$ is the damping constant. For simplicity, we assume additionally that $f(0)=0$ and $\min \pa_xf(x)=\pa_xf(0)=-c$. Then for $a=0$, this problem returns to the standard Burgers equation and by standard characteristic method, the solution of \eqref{burgersdamp} will form a shock at time $T_{\ast}=-1+\frac{1}{c}$ for $c>0$ and the location $x_{\ast}=0$ with $\pa_xu(0,t)\ra-\infty$ as $t\to T_{\ast}$.
Following D.Christodoulou, we define the eikonal function for a given solution $\fai$ to \eqref{burgersdamp} as
 \bee
 \left\{\bes
 &Lu:=(\pa_t+\fai\pa_x)u=0,\\
 &u(x,-1)=x,
 \end{split}\right.
 \ee
 where $L$ is the vector filed which generates the characteristics of \eqref{burgersdamp}. This Eikonal function yields a new coordinate system $(t,u)$ in which $L=\dfrac{\pa}{\pa t}|_{(t,u)}$. Moreover, the original equation becomes
 \bee
 L\fai=-a\fai.
 \ee
 Hence, $\fai(t,u)=\fai_0(u)e^{-a(t+1)}$ and if initial $\fai_0$ is smooth, $\fai$ together with it's derivatives will remain smooth for all time in $(t,u)$ coordinates. Thus, the only possibility of formation of singularities is that the transformation between $(t,u)$ coordinates and $(t,x)$ coordinates degenerates.\\
 \hs Define the inverse foliation density function as
 \bee
 \mu=\dfrac{1}{\pa_xu}.
 \ee
 It is clear that $\mu$ has the geometric meaning: $\dfrac{1}{\mu}$ measures the density of level sets of $u$ which are exactly the characteristics in the physical plane and initially, $\mu$ equals to $1$. Note that the Jacobian of the transformation between two coordinates is given as
 \bee
 \de=\dfrac{\pa (t,u)}{\pa (t,x)}=\mu.
 \ee
 Based on above facts, the singularity formation (shock formation) to \eqref{burgersdamp} is equivalent to that the diffeomorphism between $(t,u)$ coordinates and $(t,x)$ coordinates degenerates, that is $\mu\ra 0$, which means that the characteristics become infinitely dense at the singular points.\\
\hs Since $\mu$ satisfies the equation
\bee
\left\{\bes
&L\mu=\mu\pa_x\fai=\dfrac{\pa \fai}{\pa u},\\
&\mu(x,-1)=1,
\end{split}\right.
\ee
it follows that
\bee
\mu(t,x)=\mu_0+\pa_xf(x)\int_{-1}^te^{-a(t+1)} dt=1-\dfrac{\pa_xf(x)}{a}(e^{-a(t+1)}-1).
\ee
In conclusion,
\begin{itemize}
\item for fixed $c$,
\been[(1)]
\item if $a\geq c$, the damping effect is strong enough and $\mu>0$ for all $t>-1$, that is, the characteristics never intersect and we will obtain a global solution for \eqref{burgersdamp};
\item if $a<c$, the damping effect is too weak to prevent the shock formation and $\mu\to 0$ as $t\to T_{\ast}=-\dfrac{1}{a}\ln(1-a)-1$. Furthermore, as $t\to T_{\ast}$, $\pa_x\fai=\dfrac{1}{\mu}\dfrac{\pa \fai}{\pa u}$ blows up like $\dfrac{1}{\mu}$.
    \een
    \item For fixed $a$, the size of the initial data plays a key role for the shock formation to \eqref{burgersdamp} in the following sense:
      \been[(1)]
      \item if $c\leq a$, then $\mu>0$ for all $t>-1$ and we obtain a global solution to \eqref{burgersdamp};
      \item if $c>a$, then the large initial will lead to the shock formation at time $T_{\ast}=-\dfrac{1}{a}\ln(1-\frac{a}{c})-1$.
      \een
\end{itemize}
\subsection{\textbf{The inhomogeneous nonlinear wave equation for the isentropic and irrotational Euler equations with damping}}
\hs For a smooth isentropic flow, the entropy will remain constant from the equation of energy. In this case, it follows from \eqref{enthalpy} that
\bee
dh=\dfrac{1}{\p}dp,\ \dfrac{d\p}{dh}=\dfrac{\p}{\eta^2}.
\ee
\hs Then, the equation of momentum can be written as
\bee\label{momentum2}
\dfrac{\pa}{\pa t}v^i+v\nabla v^i+a v^i=-\dfrac{\pa h}{\pa x^i}.
\ee
Taking $curl$ on both sides of \eqref{momentum2} yields
\bee\label{vorticity}
\dfrac{\pa}{\pa t}\omega+v\nabla\omega=-a\omega+\omega\cdot\nabla v-(\text{div} v)\omega,
\ee
where $\omega=curl v$ is the vorticity of the flow, which implies $\omega$ will vanish everywhere if the initial vorticity is zero. Moreover, for a irrotational flow, there exists a potential function $\fe$ such that $v=-\nabla\fe$. Substitute this into the equation \eqref{momentum2} yields
\[
\dfrac{\pa}{\pa x^i}[-h+\dfrac{\pa}{\pa t}\fe-\dfrac{1}{2}|\nabla\fe|^2+a\fe]=0.
\]
Then, up to a constant, it holds that $h=\dfrac{\pa}{\pa t}\fe-\dfrac{1}{2}|\nabla\fe|^2+a\fe$. Then, the continuity equation can be rewritten as
\bee\label{mass2}
 g^{\a\be}(\fe,\pa\fe)\pa_{\a}\pa_{\be}\fe=\dfrac{1}{\eta^2}a(\dfrac{\pa\fe}{\pa t}-|\nabla\fe|^2),\\
\ee
where $g$ is the acoustical metric given by
\bee
\left(\begin{array}{cccc}
-\eta^2+\sum_{i}(v^i)^2 & -v^1 & -v^2 & -v^3\\
-v^1 & 1 &0 &0\\
-v^2 & 0 & 1 &0\\
-v^3 & 0 & 0 &1
\end{array}\right)
\ee
with it's inverse $g^{-1}$
\bee
\dfrac{1}{\eta^2}\left(\begin{array}{cccc}
-1 & -v^1 & -v^2 & -v^3\\
-v^1 & \eta^2-(v^1)^2 &-v^1v^2 &-v^1v^3\\
-v^2 & -v^1v^2 & \eta^2-(v^2)^2 &-v^2v^3\\
-v^3 & -v^1v^3 & -v^2v^3 &\eta^2-(v^3)^2
\end{array}\right)
\ee
\hs Note that the above equation is invariant under the following symmetries: the space translation, the time translation, and the space rotation. Let $A$ be any one of $\{\pa_{\a},\Omega_{ij}=x^i\pa_j-x^j\pa_i\}$. One cae use the symmetry generated by $A$ to act on the solution $\fe$ to get a $1-$parameter family of solutions $\{\fe_{\tau}|\fe_0=\fe\}$. Then, differentiating with respect to $\tau$ and taking value at $\tau=0$ yield the variations of $\fe$ defined to be
\[
\fai:=A\fe=\dfrac{d\fe}{d\tau}|_{\tau=0}.
\]
The above procedure will produce a equation for $\fai$, which will be called the linearized equation corresponding to \eqref{mass2} with respect to the symmetry $A$.\\
\hs Acting $A$ on the equation \eqref{mass2} directly yields the following inhomogeneous covariant wave equations for $\fai$
\bee\label{nonlinearwave}
\square_{\tilde{g}(h)}\fai=-\dfrac{2\eta'}{\eta^2}a\fai\de\fe+\dfrac{a}{\p\eta}\left(\dfrac{\pa\fai}{\pa t}-\dfrac{\pa\fe}{\pa x^i}\dfrac{\pa \fai}{\pa x^i}\right),
\ee
where $\til{g}$ is the conformal acoustical metric
\bee
\til{g}=\Omega g=\dfrac{\p}{\eta} g.
\ee
\begin{remark}
Note that the enthalpy $h$ here is different from the enthalpy in the Christodoulou's work, so that the metric $g$ here depends not only on $\pa\fe$, but also on $\fe$ itself. Moreover, due to the damping effect, the inhomogeneous terms of \eqref{nonlinearwave} will affect the behavior of $\fai$ and the inverse foliation density $\mu$ which need to be carefully studied.
\end{remark}
\begin{remark}
Note that for the metric $g$, the corresponding Christoffel symbols $\Ga_{\a\be\ga}$ satisfy
\bee\label{christoffelsymbol}
\Ga_{0ij}=\Ga_{ijk}=0,
\ee
which will be used later.
\end{remark}
\begin{remark}
Let $H=-2h-\eta^2$. Then for isentropic and irrotational flows, the condition \eqref{nonlinearcondition} is equivalent to
\bee
\dfrac{d H}{dh}\neq0,
\ee
which plays a key role in the shock formation.
\end{remark}
\subsection{\textbf{Initial data}}
\hs Let the initial data for \eqref{cee} be given on the hyperplane $\Si_{-2}$ which is isentropic and irrotational. Denote $S_{-2,a}$ to be the sphere on $\Si_{-2}$ centered at origin with radius $2+a$. Inside the sphere $S_{-2,0}$, set
\bee
\p=\p_0,\ s=s_0,\ v=0,\ \eta=\eta_0.
\ee
Then the initial data for the nonlinear wave equation \eqref{mass2} inside $S_{-2,0}$ is given by
\bee
\fe=0,\ \pa_t\fe=h_0,
\ee
where $h_0$ is the initial enthalpy. For simplicity, one can set
\bee
\p_0=1,\eta_0=1,h_0=0.
\ee
Let the initial data be only non-trivial in the $\da$-annulus region:
\bee
\Si_{-2}^{\da}:=\left\{x\in\Si_{-2} |2\leq r(x)\leq 2+\da\right\}.
\ee
\hs Following the work in \cite{Ontheformationofshocks}, we construct the short pulse data for \eqref{mass2} in $\Si_{-2}^{\da}$ as follows.
\begin{lem}
For any given seed data $(\fe_1,\fe_2)$, there exists a $\da'>0$ depending only on $(\fe_1,\fe_2)$ such that for all $\da<\da'$, there exists another function $\fe_0\in C^{\infty}((0,1]\times S^2)$ smoothly depending on $(\fe_1,\fe_2)$ with the property that if one sets the initial data of \eqref{mass2} as follows:\\
\hs for $r\leq 2$ on $\Si_{-2}$, $(\fe(-2,x),\pa_t\fe(-2,x))=(0,0)$; while for $2\leq r\leq 2+\da$,
\bee\label{initialdata}
\fe(-2,x)=\da^2\fe_0\left(\dfrac{r-2}{\da},\ta\right),\hs \pa_t\fe(-2,x)=
\da\fe_1\left(\dfrac{r-2}{\da},\ta\right).
\ee
Then, it holds that on $\Si_{-2}$,
\bee\label{pat+parfe}
|(\pa_t-\pa_r)^2\fe|\les\da.
\ee
As a result, it holds that on $\Si_{-2}$, $|(\pa_t-\pa_r)\pa_{\a}\fe|\les\da$ for $\a=0,1,2,3$. 
\end{lem}
\begin{remark}
Indeed, \eqref{pat+parfe} can be generalized for higher orders of $(\pa_t-\pa_r)$. The point is that although the high order radial derivatives of $\fe$ as well as the time derivatives are large, their suitable summations may be small (this means that the derivatives along the null hypersurfaces are small).
\end{remark}
\begin{pf}
First, take arbitrary $\fe_1$ and fix it. Since $\eta\sim1$ on $\Si_{-2}$, it follows that
\bee\label{pat2fe}
\bes
\dfrac{\pa^2\fe}{\pa t^2}&=\de\fe+2\dfrac{\pa\fe}{\pa x^i}\dfrac{\pa^2\fe}{\pa t\pa x^i}-\dfrac{\pa\fe}{\pa x^i}\dfrac{\pa\fe}{\pa x^j}\dfrac{\pa^2\fe}{\pa x^i\pa x^j}-a\left(\dfrac{\pa\fe}{\pa t}-|\nabla\fe|^2\right)\\
&=\pa^2_r\fe+\dfrac{2}{r}\pa_r\fe+\dfrac{1}{r^2}\s{\de}_{S^2}\fe+2\pa_r\fe\pa_r\pa
_t\fe-r(\pa_r\fe)^2\pa^2_r\fe+(\pa_r\fe)^3\\
&-a\left(\pa_t\fe-\sum_{i,j}\dfrac{x^ix^j}{r^2}(\pa_r\fe)^2\right).
\end{split}
\ee
Then, using \eqref{initialdata} yields that
\[
|\pa_t^2\fe-2\pa_r\pa_t\fe+\pa_r^2\fe|\les\da,
\]
which is equivalent to
\bee
\bes
&\pa^2_r\fe+\dfrac{2}{r}\pa_r\fe+\dfrac{1}{r^2}\s{\de}_{S^2}\fe+2\pa_r\fe\pa_r\pa
_t\fe-r(\pa_r\fe)^2\pa^2_r\fe+(\pa_r\fe)^3\\
&-a\left(\pa_t\fe-\sum_{i,j}\dfrac{x^ix^j}{r^2}(\pa_r\fe)^2\right)-2\pa_r\pa_t\fe
+\pa_r^2\fe=O(\da),
\end{split}
\ee
or
\bee
\bes\label{equi}
 &2\pa_s^2\fe_0+\dfrac{2}{r}\da\pa_s\fe_0+\dfrac{1}{r^2}\da^2\s{\de}_{S^2}\fe_0\\
&+2\da\pa_s\fe_0\pa_s\fe_1-r\da^2\pa_s^2\fe_0(\pa_s\fe_0)^2+\da^3(\pa_s\fe_0)^3\\
&-a\da\fe_1+\sum_{i,j}\dfrac{x^ix^j}{r^2}\da^2(\pa_s\fe_0)^2-2\pa_s\fe_1=O(\da).
\end{split}
\ee
One can assume the following ansatz for $\fe_0$:
\bee\label{ansatz}
|\pa_s\fe_0|+|\pa^2_{\ta}\fe_0|+\da|\pa_s^2\fe_0|\leq C,
\ee
where $C$ is a constant independent of $\da$. Based on \eqref{ansatz}, one can drop all the $\da$-terms or higher order terms on the left hand side of \eqref{equi} to get
\bee\label{equi2}
 2\pa_s^2\fe_0-2\pa_s\fe_1=O(\da).
\ee
This can be obtained by solving the following Cauchy problem:
\bee\label{ode1}
\left\{\bes
&\pa_s^2\fe_0-\pa_s\fe_1=\da\fe_2,\\
&\fe_0(0,\ta)=0,\hs \pa_t\fe_0(0,\ta)=0.
\end{split}\right.
\ee
By choosing $\da$ sufficiently small, one could recover \eqref{ansatz}. Hence, it holds that on $\Si_{-2}$
\bee\label{Lfai}
\bes
(\pa_t-\pa_r)\pa_r\fe&=\pa_s^2\fe_0-\pa_s\fe_1\les\da,\\
2(\pa_t-\pa_r)\pa_t\fe&=\pa_t^2\fe-2\pa_r\pa_t\fe+\pa_r^2\fe+O(\da)\\
&=(\pa_t-\pa_r)^2\fe+O(\da)\les\da.
\end{split}
\ee
\end{pf}
\section{\textbf{The Geometric formulation, basic structure equations and the main results}}\label{section2}
\subsection{\textbf{The maximal development}}
\hs Considering the following Cauchy problem on a Lorentzian manifold $(M,g)$ in the Minkowski space-time:
\bee\label{eg}
\left\{\begin{split}
&\square_{g(\fai)}\fai=0,\\
&\fai(0,x),\pa_{t}\fai(0,x)\ \text{are given on}\ \Si_{0}.
\end{split}\right.
\ee
\hs A development of a given initial data is defined as:
\been[(1)]
\item a domain $D$ in $M$, whose past boundary is $\Si_{0}$, that is, for any point $a\in D$, there exists a smooth curve $\be$ starting from $a$ to the past such that $\be$ will meet $\Si_0$;
\item a smooth solution $\fai$ to \eqref{eg} defined on $D$ with initial data on $\Si_{0}$ with following property.\\
    As long as $\fai$ is well-defined on $D$, the characteristic hypersurfaces to \eqref{eg} are well-defined. The past null cone at any point $p\in D$ is defined as the characteristic hypersurfaces emanating from $p$ to the past. Then, for any inextendible curve $\ga :[0,\tau)\ra D$ with $\ga(0)=p\in D$ towards the past and $\ga'(s)$ being inside the past null cone at $\ga(s)$ (that is, the vector $\ga'(s)$ is causal, i.e. $g(\ga'(s),\ga'(s))\leq 0$), $\ga$ must intersect $\Si_0$.
\een
\hs Indeed, by the terminology of Lorentz geometry, the above is equivalent to say $\Si_0$ is a Cauchy hypersurface of $D$, that is, every inextendible timelike curve in $D$ will intersect $\Si_0$ exactly once.\\
\hs The local well-posedness theory implies that if $(D_{1},\fai_{1})$,$(D_{2},\fai_{2})$ are two developments on same initial data, then $\fai_{1}=\fai_{2}$ on $D_{1}\cap D_{2}$. Therefore, the union of all developments of a given initial data is also a development, which is called \textbf{the maximal development}, the corresponding solution and domain are \textbf{the maximal solution} and \textbf{the domain of the the maximal solution}, respectively.\\
\hs In the following, we restrict the initial data on following the annular region:
\bee
\Si_{-2}^{\da^2}:=\{x\in\Si_{-2}:\ 2\leq r(x)\leq 2+\da^2\},
\ee
where $r=\sqrt{(x^1)^2+(x^2)^2+(x^3)^2}$.
Define the function $u:=r-2$ on $\Si_{-2}$. For each value of $u$, the corresponding level set $S_{-2,u}$ is a sphere of radius $2+u$ and as a consequence,
\bee
\Si_{-2}^{\tilde{\da}}=\bigcup_{u\in[0,\tilde{\da}]}S_{-2,u},\ \tilde{\da}:=\da^2.
\ee
\hs Consider, in the domain of maximal solution, the family of incoming characteristic null-hypersurfaces emanated from $S_{-2,u}$, denoted to be $C_u$. Obviously, $C_u\cap\Si_{-2}=S_{-2,u}$, for any $u\in[0,\tilde{\da}]$.
By the domain of dependence, the solution inside $C_{0}$ is completely trivial and it's natural to consider the solution in the region foliated by $C_u$
\bee
W_{\tilde{\da}}=\bigcup_{u\in[0,\tilde{\da}]}C_u,
\ee
which is a subset of the domain of the maximal solution. Next, the function $u$ will be extended to $W_{\tilde{\da}}$ in the following definition, which implies that its level sets are precisely the incoming characteristic null-hypersurfaces $C_u$.
\begin{defi}
The eikonal function $u$ satisfies
\bee\label{ekinoal}
\left\{\bes
&g^{\a\be}\dfrac{\pa u}{\pa x^{\a}}\dfrac{\pa u}{\pa x^{\be}}=0,\\
& u|_{\Si_{-2}}=r-2,
\end{split}\right.
\ee
and $\pa_tu>0$.
\end{defi}
\hs See the following picture for the above geometric notations.
\begin{center}
\begin{tikzpicture}
\draw (-4,0)--node[below] {$\Sigma_{-2}$}(4,0);
\draw[blue] (-4,0)--node[above left] {$C_{\tilde{\da}}$}(-2,3);
\draw (-3,1.5)--node[below] {$\Sigma_t$}(3,1.5);
\draw[blue] (-2.7,0)--node[above left] {$C_u$}(-1.5,3);
\filldraw (2.1,1.5) circle (.04);
\node at (2.1,1.2) {$S_{t,u}$};
\draw[blue] (-1.7,0)--node[above right] {$C_0$}(-1.3,3);
\draw[blue] (4,0)--node[above right] {$C_{\tilde{\da}}$}(2,3);
\draw[blue] (2.7,0)--node[above right] {$C_u$}(1.5,3);
\draw[blue] (1.7,0)--node[above left] {$C_0$}(1.3,3);
\filldraw (1.8,2.25) circle (.05);
\draw[->] (1.8,2.25)--node[above right] {$L$}(1.5,3);
\draw[->] (1.8,2.25)--node[above right] {$T$}(2.7,2.25);
\draw[->] (1.8,2.25)--node[below right] {$X$}(1.3,1.6);
\end{tikzpicture}
\end{center}
\subsection{\textbf{Frames and coordinates}}
\begin{defi}
The incoming null geodesic vector field is defined as
\bee
\hat{L}=-g^{\a\be}\pa_{\a}u\pa_{\be}.
\ee
\end{defi}
\hs It follows that $\hat{L}$ is $g-$null (since $g(\hat{L},\hat{L})=g^{\a\be}\pa_{\a}u\pa_{\be}u=0$) and $g-$orthogonal to $C_u$. Using $\hat{L}$, one can define the important inverse foliation density function as follows.
\begin{defi}
The inverse foliation density function is set to be
\bee
\mu:=\dfrac{1}{-g^{\a\be}\pa_{\a}t\pa_{\be}u}=\dfrac{1}{\hat{L}^0},
\ee
whose reciprocal measures the density of the foliation $\bigcup_{u\in[0,\tilde{\da}]}C_u\cap\Si_t$.
\end{defi}
\hs It will be shown later that $\hat{L}^{\a}$ blows up like $\frac{1}{\mu}$ when $\mu$ vanishes. So, one needs to rescale the vector filed $\hat{L}$ as
\bee
L=\mu\hat{L}.
\ee
Since $L$ is proportional to $\hat{L}$, then $L$ is also $g-$orthogonal to $C_u$ and $g-$null. Moreover, it holds that $Lt=1$.\\
\hs Define for each $u\in[0,\tilde{\da}]$ and fixed $t$, the min of $\mu$ on the set $S_{t,u}$ to be $\mu(t,u)$. Then, define
\bee\label{defmumu}
\mu_m^u=\min\{\inf_{u'\in[0,u]}\mu(t,u'),1\},
\ee
and
\bee
s_{\ast}=\sup\{t|t\geq-2\ \text{and}\ \mu_m^{\tilde{\da}}(t)>0\}.
\ee
\hs For each $u\in[0,\tilde{\da}]$, one can define $t_{\ast}(u)$ to be the lifespan of the solution to \eqref{ekinoal} and define $t_{\ast}$ to be:
\bee
t_{\ast}=\inf_{u\in[0,\tilde{\da}]}t_{\ast}(u)=\sup\{\tau|\text{smooth solution exists for all}\ (t,u)\in[-2,\tau)\times[0,\tilde{\da}]\}.
\ee
We finally restrict time on $[-2,t^{\ast})$ with
\bee\label{defsasttast}
s^{\ast}=\min\{s_{\ast},\sigma\},\quad t^{\ast}=\min\{t_{\ast},s^{\ast}\},
\ee
where $\sigma<0$ is a \textbf{fixed small constant}. In the following, we work on $W_{\tilde{\da}}^{\ast}\subset W_{\tilde{\da}}$, where
\bee
W_{\tilde{\da}}^{\ast}=\bigcup_{(t,u)\in[-2,t^{\ast})\times[0,\tilde{\da}]}S_{t,u},
\ee
where $S_{t,u}=C_u\cap\Si_t$ is a $2-$dimensional topological sphere.\\
\hs Define on $W_{\tilde{\da}}^{\ast}$ a vector filed $T$ such that it is tangential to $\Si_t$ and $g-$orthogonal to $\{S_{t,u}\}$ for $u\in[0,\tilde{\da}]$, and normalized by
\bee
Tu=1.
\ee
Let $\Lambda=[L,T]$ be the commutator of $L,T$ which is $S_{t,u}$ tangent since $\Lambda t=[L,T]t=0$ and $\Lambda u=[L,T]u=0$.
\begin{prop}\label{LTrelation}
It holds that
\bee
\bes
&g(L,T)=-\mu,\quad g(T,T)=\kappa^2=(\eta^{-1}\mu)^2\ \text{where}\ \kappa>0,\\
&L=\pa_t-\left(\eta\hat{T}^i+\dfrac{\pa\fe}{\pa x^i}\right)\dfrac{\pa}{\pa x^i},\ \text{where}\ \hat{T}=\kappa^{-1}T,\ |\hat{T}|=1.
\end{split}
\ee
Let $\dl=\eta^{-1}\kappa L+2T$. Then $\dl$ is an outgoing null vector field and $g-$orthogonal to $S_{t,u}$ verifying $g(L,\dl)=-2\mu$.
\end{prop}
\begin{pf}
\bee
g(L,T)=\mu g(\hat{L},T)=-\mu g_{\a\be}g^{\a\ga}\pa_{\ga}uT^{\be}=-\mu\da^{\ga}_{\be}T^{\be}\pa_{\ga}u=-\mu Tu=-\mu.
\ee
\hs Consider the vector field
\bee
N=-\eta^2 g^{\a\be}\pa_{\a}t\pa_{\be}.
\ee
Then, $N$ is $g-$orthogonal to $\Si_t$, future directed and timelike verifying $Nt=1$. Indeed,
\bee
N=\dfrac{\pa}{\pa t}+v^i\dfrac{\pa}{\pa x^i}.
\ee
Since $N$ is $g-$orthogonal to $S_{t,u}$, then
\bee\label{Ndecompose}
N=L+f T,
\ee
where $f$ is a function. In \eqref{Ndecompose}, taking inner product with $T$ and noticing the fact that $N$ is $g-$orthogonal to $\Si_t$ and then $g(N,T)=0$ yield $\mu=f\kappa^2$. On the other hand, since
\bee
-\eta^2=g(N,N)=g(L+fT,L+fT)=-f^2\kappa^2,
\ee
 it follows that $\mu=\eta\kappa.$ Hence,
\bee\label{Ldecomposition}
\bes
L&=N-\dfrac{\eta}{\kappa}T=N-\eta\hat{T}=\dfrac{\pa}{\pa t}-\left(\eta\hat{T}^i+\dfrac{\pa\fe}{\pa x^i}\right)\dfrac{\pa}{\pa x^i}.
\end{split}
\ee
\end{pf}
\hs One can now define a coordinate system on $S_{t,u}$ as follows.\\
\hs Since $S_{-2,0}$ is the standard Euclidean sphere and each $S_{-2,u}$ is diffeomorphic to $S_{-2,0}$, then if the local coordinates $(\til{\ta}^1,\til{\ta}^2)$ are chosen on $S_{-2,0}$, the diffeomorphism defines a local coordinate system on $S_{-2,u}$ for every $u\in[0,\tilde{\da}]$ which is denoted to be $(\ta_1,\ta_2)$. Extend $(\ta_1,\ta_2)$ to $S_{t,u}$ for $(t,u)\in[-2,t_{\ast})\times[0,\tilde{\da}]$ by requiring that
\bee
L\ta^A=0,\ \text{for}\ A=1,2.
\ee
\hs The local coordinates $(\ta^1,\ta^2)$ together with $(t,u)$ define a complete local coordinate system $(t,u,\ta^1,\ta^2)$ on $W_{\tilde{\da}}^{\ast}$, which is called the acoustical coordinates.\\
\hs Set:
\bee
\bes
\Si_t^u&:=\{(t,u',\ta)\in\Si_t|\ 0\leq u'\leq u\},\\
C_u^t&:=\{(t',u,\ta)\in C_u|\ -2\leq t'\leq t\},\\
W_t^u&:=\bigcup_{(t',u')\in[-2,t)\times[0,u]} S_{t',u'}.
\end{split}
\ee
\hs Define two $S_{t,u}-$tangent vector fields as $X_A=\dfrac{\pa}{\pa\ta^A}$ for $A=1,2$. Then, $\{X_A\}$ together with $\{L,T\}$ and $\{L,\dl\}$ form the frames $\{L,T,X_1,X_2\}$ and $\{L,\dl,X_1,X_2\}$, respectively.\\
\hs Now we turn to represent the metric $g$ in the acoustical coordinates and in the null frame $\{L,\dl,X_1,X_2\}$.\\
\hs Since $Lu=0$ and $L\ta^A=0$, it holds that
\bee\label{Lacoustical}
L=\dfrac{\pa}{\pa t}\left|_{(t,u,\ta)}\right..
\ee
Since $T$ is $\Si_t$-tangential, then one can decompose it as $T=\dfrac{\pa}{\pa u}-\Xi$, where $\Xi$ is an $S_{t,u}$ tangential vector field which can be written as $\Xi=\Xi^A\dfrac{\pa}{\pa \ta^A}$. Moreover, it holds that $[L,\Xi]=-[L,T]=-\Lambda.$
\begin{remark}
Indeed, one can set $\Xi\equiv 0$ on any one of the hypersurface $\Si_t$ for fixed $t$. 
\end{remark}
\hs With the notation
\bee
\sg_{AB}=g\left(\dfrac{\pa}{\pa\ta^A},\dfrac{\pa}{\pa\ta^B}\right),
\ee
 \eqref{Lacoustical} and Proposition\ref{LTrelation}, $g$ can be represented in $(t,u,\ta)$ as:
\bee
g=-2\mu dudt+\kappa^2 du^2+\sg_{AB}(d\ta^A+\Xi^Adu)(d\ta^B+\Xi^Bdu).
\ee
Due to $g(L,L)=g(\dl,\dl)=0$ and Proposition\ref{LTrelation}, $g$ can be represented in the frame $\{L,\dl,X_1,X_2\}$ as:
\bee\label{decomposeg2}
g^{\a\be}=-\dfrac{1}{2\mu}\left(L^{\a}\dl^{\be}+L^{\be}\dl^{\a}\right)+\sg^{AB}X_A^{\a}
X_B^{\be}.
\ee
\hs Consider the Jacobian of the change of coordinates
\bee
(t,u,\ta^1,\ta^2)\ra(x_0,x_1,x_2,x_3).
\ee
Then, it follows that:
\begin{thm}(Fundamental theorem)
\bee
\de=\det\dfrac{\pa(x_0,\cdots x_3)}{\pa(t,u,\ta^1,\ta^2)}=\eta^{-1}\mu\sqrt{\det\sg}.
\ee
Since $\eta\sim1$ (see Lemma\ref{hrhoeta}), the above relation implies that the diffeomorphism between the acoustical coordinates and the rectangular coordinates never degenerates \textbf{as long as $\mu>0$}, which implies that the frames $\{L,\dl,X_1,X_2\}$ and $\{L,T,X_1,X_2\}$ are equivalent to $\{\pa_{\a}\}_{\a=0,1,2,3}$ as long as $\mu>0$.
\end{thm}
\begin{pf}
Since $x_0=t$, it follows that
\bee
\dfrac{\pa x_0}{\pa t}=1,\ \dfrac{\pa x_0}{\pa u}=\dfrac{\pa x_0}{\pa \ta^A}=0.
\ee
By \eqref{Lacoustical}, $\dfrac{\pa x_i}{\pa t}=L^i.$ It follows from the definition of $T$ that $\dfrac{\pa x_i}{\pa u}=T^i+\Xi^AX_A^i:=T^i+\ga^i$. Set $\bar{g}$ to be the restriction of $g$ on $\Si_t$, which is the Euclidean metric. Then,
\bee
\bes
\de&=\left|\begin{array}{cccc}
1&0&0&0\\
L^1&T^1+\ga^1&X_1^1&X_2^1\\
L^2&T^2+\ga^2&X_1^2&X_2^2\\
L^3&T^3+\ga^3&X_1^3&X_2^3
\end{array}\right|
=\left|\begin{array}{ccc}
T^1+\ga^1&X_1^1&X_2^1\\
T^2+\ga^2&X_1^2&X_2^2\\
T^3+\ga^3&X_1^3&X_2^3
\end{array}\right|\\
&=\left|\begin{array}{ccc}
T^1&X_1^1&X_2^1\\
T^2&X_1^2&X_2^2\\
T^3&X_1^3&X_2^3
\end{array}\right|=|T|\cdot|X_1\wedge X_2|=\eta^{-1}\mu\sqrt{\det\sg}.
\end{split}
\ee
\end{pf}
\begin{remark}
This theorem 
implies that the regular solutions in the acoustical coordinates can be brought back to the rectangular coordinates as long as $\mu>0$. As it will be shown that there is no singularity in the acoustical coordinates, the only possibility for the singularity(shock) formation is the Jacobian of the change of the coordinates vanishing. As a consequence, the shock formation is equivalent to that the transformation between the acoustical coordinates and the rectangular coordinates degenerates, that is $\mu\to 0$, which means that the collapse of characteristic null hypersurfaces $C_u$.
\end{remark}
\subsection{\textbf{Connection coefficients and $2^{nd}$ fundamental forms $k,\chi$ and $\ta$}}
\begin{defi}
Let $\bar{\mathcal{L}}$ and $\lie$ be the restriction of the lie derivative $\mathcal{L}$ on $\Si_t$ and $S_{t,u}$ respectively.
\been
\item The $2^{nd}$ fundamental form $k$ for the embedding $\Si_t^u\hookrightarrow W_{t}^{u}$ is defined as:
    \bee
    2\eta k=\bar{\mathcal{L}}_{N}g.
    \ee
\item The $2^{nd}$ fundamental form $\chi_{AB}$ for the embedding $S_{t,u}\hookrightarrow C_u^t$ is defined as:
    \bee
    2\chi_{AB}=\lie_L\sg_{AB}=2g(D_{A}L,X_B).
    \ee
\item The $2^{nd}$ fundamental form $\ta_{AB}$ for the embedding $S_{t,u}\hookrightarrow \Si_t^u$ is defined as:
    \bee
    2\eta^{-1}\mu\ta_{AB}=\lie_T\sg_{AB}=2g(D_AT,X_B).
    \ee
\een
\end{defi}
\hs It follows from the definitions that $\chi=\eta(\s{k}-\ta)$ where $\s{k}$ is the restriction of $k$ on $S_{t,u}$. Indeed, one can compute $k$ as follows:
\begin{align}\label{kij}
2\eta k_{ij}&=(\bar{\mathcal{L}}_Ng)_{ij}=N(g(\pa_i,\pa_j))-g([N,\pa_i],\pa_j)-g(\pa_i,[N,\pa_j])\\
&=g_{mj}\pa_iv^m+g_{im}\pa_jv^m=-2\pa_i\fai_j.
\end{align}
Set the following three $1-$forms $\ep,\zeta,\eta$ on $S_{t,u}$ as:
\bee
\bes
\eta^{-1}\mu\ep_A=k(X_A,T),\ \zeta_A=g(D_AL,T),\ \eta_A=-g(D_AT,L).
\end{split}
\ee
Note that
\bee\label{etaAcompute}
\zeta_A=\eta\kappa\ep_A-\kappa X_A\eta,\quad \eta_A=\zeta_A+X_A\mu.
\ee
Then, the connection coefficients can be computed by in the frame $\{L,T,X_1,X_2\}$ as follows, which can be verified directly.
\begin{lem}\label{connection1}
\begin{align*}
D_LL&=\mu^{-1}(L\mu)L,\quad D_TL=\sg^{AB}\eta_BX_A-\eta^{-1}L(\kappa)L,\\
D_LT&=-\sg^{AB}\zeta_BX_A-\eta^{-1}L(\kappa)L,\quad D_AL=-\mu^{-1}\zeta_AL+\chi_A^BX_B=D_LX_A,\\
D_TT&=\eta^{-2}\kappa(T\eta+L\kappa)L+(\eta^{-1}L\kappa+\mu^{-1}T\mu)T-\dfrac{1}{2}\sg^{AB}X_B(\kappa^2)X_A,\\
D_AT&=\eta^{-1}\kappa\ep_AL+\mu^{-1}\eta_AT+\kappa\ta_A^CX_C,\quad D_AX_B=\s{D}_AX_B+\eta^{-1}\s{k}_{AB}L+\mu^{-1}\chi_{AB}T.
\end{align*}
\end{lem}
\hs We will use the following structure equations whose proofs are given in\cite{christodoulou2014compressible}.
\been[(1)]
\item The Gauss equation for the $2^{nd}$ fundamental form $\ta$ is:
      \bee\label{gaussta}
      \dfrac{1}{2}(tr\ta)^2-\dfrac{1}{2}|\ta|^2_{\sg}=K,
      \ee
      where $K$ is the Gauss curvature of $S_{t,u}$.
\item The Codazzi equations for the $2^{nd}$ fundamental form $\chi$ are:
      \bee\label{Codazzi1}
      \s{D}_A\chi_{BC}-\s{D}_B\chi_{AC}=R_{ABCL}-\mu^{-1}(\zeta_A\chi_{BC}-
      \zeta_B\chi_{AC})=\beta_C\varepsilon_{AB}-\mu^{-1}(\zeta_A\chi_{BC}-
      \zeta_B\chi_{AC}),
      \ee
      where $R$ is the Riemannian curvature tensor, $\beta$ is the $1-$form and $\varepsilon_{AB}$ is the area form of $\sqrt{\det\sg}$. Contracting $\sg^{AC}$ on both sides of \eqref{Codazzi1} leads to
      \bee\label{Codazzi2}
      \s{div}\chi-\sd tr\chi=\beta^{\ast}-\mu^{-1}(\zeta\cdot\chi-\zeta\cdot
      tr\chi),
      \ee
      where $\beta^{\ast}$ is the $\sg$-dual of $\beta.$
\item \bee\label{Ttrchi}
      Ttr\chi=\s{D}^B\eta_B+\mu^{-1}\zeta\cdot\eta-\eta^{-1}L(\eta^{-1}\mu)tr\chi-\eta^{-1}\mu\ta\cdot\chi-\sg^{AB}R_{ATBL}.
      \ee
      It follows from \eqref{etaAcompute} and \eqref{Ttrchi} that $Ttr\chi-\s{\de}\mu\sim\mu\sd^2\fai+l.o.ts.$
\een
\hs In the null frame $\{L,\dl,X_1,X_2\}$, the connection coefficients are given as follows.
\begin{lem}\label{connection2}
\begin{align*}
D_LL&=\mu^{-1}(L\mu)L,\quad D_{\dl}L=-L(\eta^{-1}\kappa)L+2\eta^AX_A,\\
D_AL&=-\mu^{-1}\zeta_AL+\chi_A^B X_B=D_LX_A,\quad D_L\dl=-2\zeta^AX_A,\\
D_{\dl}\dl&=(\mu^{-1}\dl\mu+L(\eta^{-1}\kappa))\dl-2\mu\sg^{AB}X_B(\eta^{-2}\kappa)
X_A,\\
D_A\dl&=\mu^{-1}\eta_AL+\dc_A^BX_B,\quad D_AX_B=\s{D}_AX_B+\dfrac{1}{2}\mu^{-1}\dc_{AB}L+\dfrac{1}{2}\mu^{-1}\chi_{AB}\dl,
\end{align*}
where $\dc$ is the $2^{nd}$ fundamental form for the embedding $S_{t,u}\hookrightarrow C_u^t$ associated to $\dl$, that is, $2\dc=\lie_{\dl}\sg.$
\end{lem}
\subsection{\textbf{Decomposition of wave equation in the null frame}}
\hs Let $\s{\de}$ be the Laplacian operator on $S_{t,u}$. 
It follows from \eqref{decomposeg2} and Lemma\ref{connection2} that for any function $f$,
\begin{align}
\square_gf&=g^{\a\be}D^2_{\a\be}f=-\mu^{-1}(L(\dl f)-(D_L\dl)f)+\sg^{AB}[(\s{D}^2_{AB})f+(\s{D}_AX_B)f-(D_AX_B)f]\\
&=\s{\de}f-\frac{1}{2}\mu^{-1}tr\dc(Lf)-\frac{1}{2}\mu^{-1}tr\chi(\dl f)-\mu^{-1}L(\dl f)-2\mu^{-1}\zeta\cdot\sd f.\label{decomposewavenonconformal}
\end{align}
\hs Direct computations yield
\bee\label{relationconfromal}
\bes
\square_{\til{g}}f&=\Omega^{-1}\square_gf+\Omega^{-2}\dfrac{d\Omega}{dh}g^{\a\be}\pa_{\a}h\pa_{\be}f.
\end{split}
\ee
Then, it follows from\eqref{decomposewavenonconformal} and \eqref{relationconfromal} that
\bee\label{decomposewaveeq}
\bes
\square_{\til{g}}f&=\Omega^{-1}\s{\de}f-\Omega^{-1}\mu^{-1}L(\dl f)+2\Omega^{-1}\mu^{-1}\zeta\cdot\sd f\\
&-\dfrac{1}{2}\mu^{-1}\Omega^{-1}\left((tr\dc+\dfrac{d\log\Omega}{dh}\dl h)Lf+
(tr\chi+\dfrac{d\log\Omega}{dh}Lh)\dl f\right).
\end{split}
\ee
\subsection{\textbf{Curvature tensor}}
The Riemannian curvature tensor is defined as:
\bee
\bes
R_{\a\be\ga\da}&=\dfrac{1}{2}\left(\pa_{\a}\pa_{\da}g_{\be\ga}+\pa_{\be}\pa_{\ga}
g_{\a\da}-\pa_{\a}\pa_{\ga}g_{\be\da}-\pa_{\be}\pa_{\da}g_{\a\ga}\right)\\
&+g^{\kappa\lambda}\left(\Ga_{\a\da\kappa}\Ga_{\be\ga\lambda}-
\Ga_{\a\ga\kappa}\Ga_{\be\da\lambda}\right),
\end{split}
\ee
where the first line is the principle part since it contains the second derivatives of $g$. In our case, the non-vanishing principle part of the component of the curvature tensor is $R_{0i0j}$ with its principle part given by
\bee
R_{0i0j}^{[P]}=-\dfrac{1}{2}\dfrac{dH}{dh}D^2_{ij}h.
\ee
\hs Define
\bee
\a_{AB}=R_{\a\be\ga\da}X_A^{\a}L^{\be}X_B^{\ga}L^{\da}=R_{ijkl}X_A^iL^jX_B^kL^l+R_{i0jk}X_A^iX_B^jL^k+R_{ijk0}X_A^iL^jX_B^k+
R_{0i0j}X_A^iX_B^j.
\ee
\hs Hence,
\bee
\bes
\a_{AB}^{[P]}&=R_{0i0j}^{[P]}X_A^iX_B^j=-\dfrac{1}{2}\dfrac{dH}{dh}D^2_{AB}h=-\dfrac{1}{2}\dfrac{dH}{dh}\left(\s{D}^2_{AB}h-\eta^{-1}\s{k}_{AB}Lh-\mu^{-1}
\chi_{AB}Th\right).
\end{split}
\ee
Note that the term
\bee
\dfrac{1}{2}\mu^{-1}\dfrac{dH}{dh}\chi_{AB}Th
\ee
in $\a_{AB}^{[P]}$ is singular in $\mu$. Then, let $\a_{AB}'=\a_{AB}-\dfrac{1}{2}\mu^{-1}\dfrac{dH}{dh}\chi_{AB}Th$, which is not singular in $\mu$. Precisely,
\bee
\a_{AB}'=-\dfrac{1}{2}\dfrac{dH}{dh}\left(\s{D}^2_{AB}h-\eta^{-1}\s{k}_{AB}Lh\right)+\a_{AB}^{[N]},
\ee
with $\a_{AB}^{[N]}$ given by
\bee\label{a_AB^N}
X_A^iX_B\fai_iL^{\a}L\fai_{\a}-L^{\a}L^{\be}X_A\fai_{\a}X_B\fai_{\be}+X_AhX_Bh+2LhX_A^iX_B\fai_i-X_AhL^{\a}X_B\fai_{\a}-X_BhL^{\a}X_A\fai_{\a},
\ee
which contains at most first order derivatives of $g$.\\
\hs As it will be shown that $\chi_{AB}$ and $Th$ are bounded, it follows that the shock formation is equivalent to that the curvature tensor $\a_{AB}$ blows up (since $\frac{dH}{dh}\neq 0$). Since the curvature tensor is purely a geometric entity (tensorial), the shock formation mechanism doesn't depend on the choice of the coordinates and the frames.
\subsection{\textbf{The rotation vector fields}}
\begin{defi}
Let $\mari{R}_i$ be the standard space rotation in the Euclidean space such that $\mari{R}_i=\ep_{ijk}x^j\frac{\pa}{\pa x^k}=\frac{1}{2}\ep_{ijk}(x^j\pa_k-x^k\pa_j)$ where $\ep_{ijk}$ is the skew-symmetric symbol. Then, define the rotation vector fields as follows:
\bee
R_i=\Pi\mari{R}_i, \ i=1,2,3,
\ee
where $\Pi$ is the projection from $\Si_t$ to $S_{t,u}$. Precisely,
\bee
\Pi_{i}^j=\da_i^j-\bar{g}_{ik}\hat{T}^k\hat{T}^j=\da_i^j-\hat{T}^i\hat{T}^j.
\ee
\end{defi}
Let $\{Z_i\}$ with $i=1,2,3,4,5$ be the set of the commutation vectorfields such that
\bee
Z_1=Q=tL,\ Z_2=T,\ Z_{i+2}=R_i,\ \text{for}\ i=1,2,3.
\ee
\subsection{\textbf{Transport equations for $\mu$ and $\chi$}}
\hs Here we derive the governing equations for the inverse foliation density function $\mu$ and the $2^{nd}$ fundamental form $\chi$, which will play a key role in the analysis later.
\begin{prop}
$\mu$ and $\chi$ satisfy the following transport equations respectively.
\begin{align}
L\mu&=m+\mu e,\label{transportmu}\\
\begin{split}\label{transportchi}
L\chi_{AB}&=(\mu^{-1}L\mu)\chi_{AB}+\chi_A^C\chi_{BC}-\a_{AB}\\
&=e\chi_{AB}+\chi_{A}^C\chi_{BC}-\a'_{AB}+a\eta^{-1}\hat{T}\fe\chi_{AB},
\end{split}\\
Ltr\chi&=(\mu^{-1}L\mu)tr\chi-|\chi|^2_{\sg}-tr\a.\label{transporttrchi}
\end{align}
with
\begin{align}
m&=\dfrac{1}{2}\dfrac{dH}{dh}Th+aT\fe,\ H=-2h-\eta^2,\label{mcomputation}\\
e&=\dfrac{1}{2\eta^2}\left(\dfrac{\p}{\p'}\right)'Lh+\eta^{-1}\hat{T}^i(L\fai_i).\label{edefinition}
\end{align}
\end{prop}
\begin{pf}
Since $g(N,T)=0$, $g(N,X_A)=0$ and $g(N,L)=-\eta^2$ due to Proposition\ref{LTrelation}, it follows from Lemma\ref{connection1} that
\bee
-\kappa(T\eta+L\kappa)=g(D_TT,N)=-g(T,D_TN)=-\eta k(T,T),
\ee
which implies
\bee\label{lmudefi}
L\mu=\kappa L\eta+\eta L\kappa=\kappa L\eta-\eta T\eta+\eta\mu k(\hat{T},\hat{T}).
\ee
Since $\eta\mu k(\hat{T},\hat{T})=-\mu\hat{T}^i\hat{T}^j\pa_i\fai_j=-\eta \hat{T}^i(T\fai_i)$ due to \eqref{kij} and Proposition\ref{LTrelation} implies that
\begin{align*}
-\eta\hat{T}^i(T\fai_i)&=L^i(T\fai_i)+\fai_iT\fai_i\\
&=T^jL^i\pa_j\fai_i+L^0T^i\pa_0\fai_i-L^0T^i\pa_0\fai_i+\fai_iT\fai_i\\
&=T^i(L\fai_i)-Th+aT\fe,
\end{align*}
it follows that the right hand side of \eqref{lmudefi} can be computed as
\bee
\bes
&\dfrac{d\eta}{dh}\dfrac{1}{\eta}\mu Lh-\dfrac{d\eta}{dh}\eta Th-\eta\hat{T}^i(T\fai_i)\\
&=\mu\left[\dfrac{d\eta}{dh}\dfrac{1}{\eta}Lh+\dfrac{1}{\eta}\hat{T}^i(L\fai_i)
\right]-\dfrac{d\eta}{dh}\eta Th-Th+aT\fe\\
&=\mu e+m,
\end{split}
\ee
with
\[
m=\dfrac{1}{2}\dfrac{dH}{dh}Th+aT\fe,\quad e=\dfrac{1}{2\eta^2}\left(\dfrac{\p}{\p'}\right)'Lh+\eta^{-1}\hat{T}^i(L\fai_i).
\]
Next, it follows from the definition of $\chi_{AB}$ that
\bee
\lie_L\chi_{AB}=L(g(D_AL,X_B))=g(D_AD_LL+D_{[L,A]}L-R_{LA}L,X_B)+g(D_AL,D_LX_B),
\ee
where $R$ is the curvature tensor. It follows from Lemma\ref{connection1} and $[L,X_A]=0$ that
\bee
\lie_L\chi_{AB}=(\mu^{-1}L\mu)\chi_{AB}-R_{LALB}+\sg_{CD}\chi_A^C\chi_B^D=(\mu^{-1}L\mu)\chi_{AB}+\chi_A^C\chi_{BC}-\a_{AB}.
\ee
\hs Contracting $\sg^{AB}$ on both sides and using the fact that $\lie_L\sg^{AB}=-2
\sg^{AC}\sg^{BD}\chi_{CD}$ yield the desired equation for $tr\chi$, \eqref{transporttrchi}.
\end{pf}
\subsection{\textbf{The main results}}
\hs In the end of this section, we state the main result in this paper.
\begin{thm}
 For the nonlinear wave equation \eqref{mass2}, 
 let the initial data $(\fe,\pa_t\fe)$ be given as in Section 1.6.\\
\hs Then, there exist a positive integer $N_{top}$ and a large positive constant $M$ depending only on the initial data such that the following estimates hold on $W_{\tilde{\da}}^{\ast}$ for $2\leq|\a|\leq [\dfrac{N_{top}}{2}]+3$ and $t\in[-2,t^{\ast})$ provided that $\da$ is sufficiently small:
    \bee
    \max_{\a}\max_{i_1,...i_n}\da^m||R_{i_n}\cdots R_{i_1}(T)^mQ^p\fai||_{\supnormda}+||\fai||_{\supnormda}+||L\fai||_{\supnormda}+\da||T\fai||_{\supnormda}+||\sd\fai||_{\supnormda}\les\da M,
    \ee
    with $p+m+n=|\a|$. In particular, $\fai,\mu,\chi$ and $g$ are regular on $\Si_{t^{\ast}}$ in the acoustical coordinates, as well their derivatives.\\
\hs In the rectangular coordinates, there are following two possibilities :
\been[(1)]
\item if the following largeness condition is satisfied:
      \begin{center}
      \textbf{$\max\pa_{s}\fe_1=c\geq-\frac{2}{a^{\ast}}>0$},
      \end{center}
      where $a^{\ast}=4\int_{-2}^{\sigma}\dfrac{e^{-a(\tau+2)}}{\tau}\ d\tau$ is a constant depending only on $a$, then there exists a $\da_0>0$ such that for all $\da<\da_0$, shock formation to the nonlinear wave equations \eqref{nonlinearwave} occurs before $t=T_{\ast}$, where $T_{\ast}\in[-2,t^{\ast})$ is defined as
      \bee
      c\int_{-2}^{T_{\ast}}\dfrac{e^{-a(\tau+2)}}{\tau}\ d\tau=-1+O(\da).
      \ee
      \hs Moreover, in this case, the inverse foliation density $\mu$ vanishes at some points on $\Si_{T_{\ast}}$, at which the Jacobian $\de$ of the transformation between the acoustical coordinates and the rectangular coordinates becomes $0$. 
  Moreover, the incoming characteristic null hypersurfaces $C_u$ become infinitely dense and $\pa_i\pa_j\fe$ blows up for some $i,j\in\{1,2,3\}$ at these points while $\fe$ and $\pa\fe$ remain bounded. That is, for the compressible Euler equations, $\pa_iv^j$ and $\pa_i\p$ blow up before $t=T_{\ast}$ while $v^i$ and $\p$ remain bounded. Therefore, 
  if the largeness condition holds, the damping term can not prevent the shock formation for the short pulse data. Furthermore, $T_{\ast}$ is an increasing function of $a$, then the damping effect will shift the shock time in the following sense compared with the undamped case:
\begin{itemize}
\item if $a>0$, then the damping effect will delay the formation of shock;
\item if $a<0$, then the anti-damping term will lead to the formation of shock in advance.
\end{itemize}
\item If the initial data satisfies:
      \begin{center}
      \textbf{$\max\pa_{s}\fe_1=c'\leq-\frac{1}{2a^{\ast}}$ },
      \end{center}
      then there exists a $\da_0>0$ such that for all $\da<\da_0$, smooth solutions to the nonlinear wave equations \eqref{nonlinearwave} exist on $[-2,t^{\ast}]$.
\een

\end{thm}
\begin{remark}
For simplicity, here, $a$ is a finite number (in the sense that there exists a universal constant $C$ such that $|a|\leq C$) so that the method and result in the paper can be generalized to the general nonlinear bounded damping. 
\end{remark}
\begin{remark}
For the case (2) above, we only establish the global existence on $[-2,t^{\ast}]$ while for $t>t^{\ast}$, there is no conclusion about the finite time shock formation or the smooth global existence of the solutions.
\end{remark}
\section{\textbf{Bootstrap assumptions and preliminary estimates}}\label{section3}
\hs Let $\fe$ be the solution to \eqref{mass2} and $\fai$ be its variations which are guaranteed by the local well-posedness theory. By the choice of the initial data, it holds that on $\Si_{-2}$
\begin{align}
&|L\fai|=|(\pa_t-\pa_r)\fai|\les\da,\ |T\fai|=|\pa_r\fai|\les 1,\ |\sd\fai|\les\da,\\
&\da^m|R_{i_n}\cdots R_{i_1}(T)^mQ^p\fai|\les \da,
\end{align}
where $p\leq2$. Then, assume that the following bootstrap assumptions hold for all $t\in[-2,t^{\ast})$.
\been
\item
 For $2\leq|\a|\leq N_{\infty}:=[\dfrac{N_{top}}{2}]+3$,
\bee\label{bs}
\max_{\a}\max_{i_1,...i_n}\da^m||R_{i_n}\cdots R_{i_1}(T)^mQ^p\fai||_{\supnormda}+||\fai||_{\supnormda}+||L\fai||_{\supnormda}+\da||T\fai||_{\supnormda}+||\sd\fai||_{\supnormda}\les\da M,
\ee
for $p+m+n=|\a|$ with $p\leq 2$;
\item
\bee\label{bs2}
||T\fe||_{\supnormda}\les M,
\ee
\een
where $M$ is a large constant independent of $\da$ to be determined later and $N_{top}$ is an integer measuring the totally derivatives acting on $\fai$. 
For $M$ large enough, it follows immediately that $|\fe|\leq|\int_{-2}^t\dfrac{\pa \fe}{\pa t} ds|+|O(\da^2)|\les\da M$.
\begin{lem}\label{hrhoeta}
Under the bootstrap assumptions \eqref{bs} and \eqref{bs2}, the following estimates hold for sufficiently small $\da$:
\beeq
&&||h||_{\supnormda}\les\da M,\hs ||\p-\p_0||_{\supnormda}\les\da M,\hs ||\eta-\eta_0||_{\supnormda}\les\da M,\\
&&||L\fe||_{\supnormda}\les\da M,\hs ||\sd\fe||_{\supnormda}\les\da M.
\eeq
\end{lem}
\hs As a consequence, $C^{-1}\leq\Omega\leq C$ for some constant $C>1$.
\begin{pf}
The estimate for $h$ follows immediately from the definition $h=\fai_0-\dfrac{1}{2}(\fai_i)^2+a\fe$. Recall that $\p$ is a smooth function of $h$ and $|\p-\p_0|\leq|h|$. Thus, $||\p-\p_0||_{\supnormda}\les\da M$. It follows from $\eta(h)=\left(\dfrac{\p'(h)}{\p(h)}\right)^{-\frac{1}{2}}$ that $||\eta-\eta_0||_{\supnormda}\les\da M$. The bounds for $\p-\p_0$ and $\eta-\eta_0$ imply that $\p\sim1$ and $\eta\sim1$ for sufficiently small $\da$.\\
\hs \eqref{Ldecomposition} implies that
\bee
|L\fe|=|\dfrac{\pa}{\pa t}\fe-(\eta\hat{T}^i+\fai_i)\dfrac{\pa}{\pa x^i}\fe|=|\fai_0-\eta\hat{T}^i\fai_i-\sum_i(\fai_i)^2|\les\da M.
\ee
\hs It follows from \eqref{bs} that
\bee
|\sd\fe|^2=\sg^{AB}X_A\fe X_B\fe=\sg^{AB}X_A^iX_B^j\fai_i\fai_j\leq\bar{g}^{ij}\fai_i\fai_j=\sum_i(\fai_i)^2\les \da^2 M^2.
\ee
\end{pf}
\hs In the following, the bounds $\p\sim 1$ and $\eta,\eta'\sim 1$ will be used without mention.
\subsection{\textbf{Preliminary estimates for the metric, $\mu$ and the second fundamental forms under the bootstrap assumptions}}
\begin{lem}\label{esthme}
Under the bootstrap assumptions \eqref{bs} and \eqref{bs2}, the following estimates hold for sufficiently small $\da$:
\beeq
&&||L(h,\p,\eta)||_{\supnormda}+||\sd (h,\p,\eta)||_{\supnormda}+\da||T(h,\p,\eta)||_{\supnormda}\les\da M,\\
&&||m||_{\supnormda}\les M, \hs||e||_{\supnormda}\les\da M.
\eeq
\end{lem}
\hs The proof is trivial due to the definitions of $h,\p,\eta,m,e$.
\begin{lem}\label{estmu}
Under the same assumptions as in Lemma\ref{esthme}, it holds that for sufficiently small $\da$
\[
||\mu||_{\supnormda}+||L\mu||_{\supnormda}\les M.
\]
As a consequence, one can obtain $||T\fe||_{\supnormda}\les \da M^2$, which recovers the bootstrap assumption \eqref{bs2} for sufficiently small $\da$.
\end{lem}
\begin{pf}
Since $L=\dfrac{\pa}{\pa t}|_{(t,u,\ta)}$, integrating \eqref{transportmu} along integral curves of $L$ yields
\beeq
\mu(t,u,\ta)&=&e^{\int_{-2}^te(\tau)d\tau}\left(\mu(-2,u,\ta)+
\int_{-2}^te^{\int_{-2}^{\tau}-e(s)ds}m(\tau,u,\ta)d\tau\right)\\
&\les& M.
\eeq
The estimate for $L\mu$ follows from $|L\mu|\leq|m|+|\mu||e|$ and the estimate for $\mu$.\\
\hs 
  It follows from the definition of $\hat{T}$ that $|T\fe|=|\eta^{-1}\mu\hat{T}^i\fai_i|\les\da M^2$. Hence, one can recover the bootstrap assumption \eqref{bs2} for sufficiently small $\da$. From now on, the bootstrap assumptions mean the assumptions \eqref{bs} and we will use the refined estimate $|T\fe|\les\da M^2$.
\end{pf}
Let $\dfrac{\sg_{AB}}{u-t}$ be the null $2^{nd}$ fundamental form in Minkowski space (c.f.\cite{Ontheformationofshocks} (3.17)). Then set
\begin{align*}
\chi'_{AB}&=\chi_{AB}+\dfrac{\sg_{AB}}{u-t},\ \ta'_{AB}=\ta_{AB}-\dfrac{\sg_{AB}}{u-t}.
\end{align*}
\begin{lem}\label{2ndff}
Under the bootstrap assumptions, it holds that for sufficiently small $\da$
\bee
||\s{k}||_{\supnormda}\les\da M,\quad ||\chi'||_{\supnormda}\les\da M, \quad ||\ta'||_{\supnormda}\les\da M.
\ee
\end{lem}
\begin{remark}
As a corollary, one can estimate $||\chi_{AB}||_{\supnormda}$ and $||\ta_{AB}||_{\supnormda}$, which are bounded by some constant $C$ (Since $\sg^{AB}X_A^iX_B^j+\hat{T}^i\hat{T}^j=\da^{ij}$). Since $\chi'$ has desired estimate, then when dealing with $\chi$, one can decompose it to be a $\da$-term and the term with $\dfrac{1}{u-t}$, where the latter one can be solved by Gronwall inequality. Moreover, when considering the $2^{nd}$ fundamental form with respect to $\til{g}$: $\til{\chi}_{AB}=\til{g}(\til{D}_{A}L, X_{B}),$ one can derive the same estimate, since $\til{\chi}_{AB}-\chi_{AB}=\dfrac{1}{2}L(\Omega)\sg_{AB}.$
\end{remark}
\begin{pf}
Since $\s{k}_{AB}=X_A^iX_B^jk_{ij}=-\eta^{-1}X_A^iX_B^j\pa_i\fai_j$, it follows that
\beeq
|\s{k}|^2&=&\sg^{AC}\sg^{BD}\s{k}_{AB}\s{k}_{CD}=\eta^{-2}\sg^{AC}\sg^{BD}X_A^i
X_B^j\pa_i\fai_jX_C^kX_D^l\pa_k\fai_l\\
&\leq&\eta^{-2}\bar{g}^{ik}\sg^{BD}X_B^j\pa_j\fai_iX_D^l\pa_l\fai_k=
\eta^{-2}\sum_i\sg^{BD}X_B^j\pa_j\fai_iX_D^l\pa_l\fai_i\\
&\leq&\eta^{-2}\sum_i|\sd\fai_i|^2\les \da^2 M^2.
\eeq
As a consequence, it holds that
\[
|\sd^2\fe|^2=\sg^{AC}\sg^{BD}X_A^i
X_B^j\pa_i\fai_jX_C^kX_D^l\pa_k\fai_l\les|\s{k}|^2\les\da^2 M^2,
\]
which will be used to estimate the curvature tensor $\a'_{AB}$.\\
 \hs It follows from \eqref{transportchi} and $L\dfrac{\sg_{AB}}{u-t}=\dfrac{2\chi'_{AB}}{u-t}-\dfrac{\sg_{AB}}{(u-t)^2}$ that
\bee\label{Lchi'AB}
L\chi'_{AB}=e\chi'_{AB}+\chi_A^{'\ C}\chi'_{BC}-\a'_{AB}+a\eta^{-1}\hat{T}\fe\chi'_{AB}-\dfrac{(e+a\eta^{-1}\hat{T}^i\fai_i)\sg_{AB}}{u-t}.
\ee
Therefore,
\bee
L(|\chi'|^2_{\sg})=2(e+a\eta^{-1}\hat{T}\fe)|\chi'|^2_{\sg}+\frac{2}{u-t}|\chi'|_{\sg}^2-2\chi^{'\ AC}\chi_A^{'\ D}\chi'_{CD}
-2\chi'^{AB}\a'_{AB}-\dfrac{2(e+a\eta^{-1}\hat{T}^i\fai_i)}{u-t}tr\chi',
\ee
which implies that
\bee\label{lt-uchi'}
L((t-u)^2|\chi'|_{\sg})\leq(t-u)^2\left((|\chi'|+|e|+|a\eta^{-1}\hat{T}^i\fai_i|)|\chi'|+|\a'|+\frac{|e|+|a\eta^{-1}\hat{T}^i\fai_i|}{u-t}\right).
\ee
It remains to estimate $\a'_{AB}$. Due to the expression $\a'_{AB}=-\dfrac{1}{2}\dfrac{dH}{dh}\s{D}^2_{A,B}h+\dfrac{1}{2}\dfrac{dH}{dh}
\eta^{-1}\s{k}_{AB}Lh+\a_{AB}^{[N]}$ and \eqref{a_AB^N}, where the second term and the last term are bounded by $\da M$, it remains to estimate $\s{D}^2_{A,B}h$. It follows from Proposition\ref{relationangular} that
\bee
|\s{D}^2_{A,B}h|\les|\s{d}^2 h|+|\s{d}h|\les|\sd^2\fai_0|+|\sd(\fai_i\sd\fai_i)|+|\sd^2\fe|\les\da M,
\ee
 which implies that $|\a'|\les\da M$.\\
\hs Let $P(t)$ be the set of $t'$ such that $||\chi'||_{L^{\infty}(\Si_{t'}^{\tilde{\da}})}\leq C_0\da M$ holds for all $t'\in[-2,t]$ with sufficiently small $\da$, where $C_0$ is a constant to be determined later. Initially, $\chi_{AB}=\eta(\s{k}_{AB}-\ta_{AB})=\eta(\sd\fai-\dfrac{\sg_{AB}}{u+2})$ and $|\chi'_{AB}|=|\chi_{AB}+\dfrac{\sg_{AB}}{u+2}|\leq|(1-\eta)||\frac{\sg_{AB}}{u+2}|+|\eta\sd\fai|$. Then, it follows from the construction of the initial data and choosing $C_0$ large enough that $||\chi'||_{L^{\infty}(\Si_{-2}^{\tilde{\da}})}<C_0\da<\frac{M}{4}C_0\da$. That is, $-2\in P(t)$.\\
  \hs Let $t_0$ be the upper bound of $P(t)$. For $t\in[-2, t_0]$, it holds that $|\chi'|+|e|+|a\eta^{-1}\hat{T}^i\fai_i|\leq (C_0+C_1)\da M$ and $|\a'|\leq C_2\da M$ for some constants $C_1$ and $C_2$. Thus, it follow from \eqref{lt-uchi'} that
\bee\label{L(t-u)^2chi}
L((t-u)^2|\chi'|_{\sg})\leq(t-u)^2\left[ (C_0+C_1)\da M|\chi'|_{\sg}+(C_1+C_2)\da M\right].
\ee
Let $x=(t-u)^2|\chi'|$. Then, along integral curves of $L$, one can rewrite \eqref{L(t-u)^2chi} as
\bee
\dfrac{dx}{dt}\leq fx+g,
\ee
where $f=(C_0+C_1)\da M$ and $g=C_3(C_1+C_2)\da M$. Integrating from $-2$ to $t$ yields
\bee
x(t)\leq e^{\int_{-2}^tf\ dt'}\left(x(-2)+\int_{-2}^tge^{\int_{-2}^{t'}-f\ ds} dt'\right).
\ee
Since $(t-u)^2\sim 1$, then
\[
|\chi'|_{\sg}\leq e^{2(C_0+C_1)\da M}\left(||\chi'||_{L^{\infty}(\Si_{-2}^{\da})}+C_3(C_1+C_2)\da M\right).
\]
Choose $C_0$ such that $C_0>4C_3(C_1+C_2)$. Then, it holds that $||\chi'||_{\supnormda}<C_0\da M$ for $\da<\dfrac{\log 2}{2(C_0+C_1)M}$. Hence, by a continuity argument, one shows that $P(t)$ holds for all $t\in[-2,t^{\ast})$.\\
\hs The estimate for $\ta'$ follows from $|\ta'_{AB}|=|\s{k}_{AB}-\eta^{-1}\chi_{AB}-\dfrac{\sg_{AB}}{u-t}
|\les\da M$.
\end{pf}
\begin{lem}\label{LTidTi}
It holds that
\begin{align}
L\hat{T}&=(\hat{T}^j(X^A\fai_j)+X^A(\eta))X_A,\label{LTi}\\
X_A\hat{T}&=\ta_A^BX_B,\\
T\hat{T}&=-X^A(\kappa)X_A.
\end{align}
Therefore, the following bounds hold for sufficiently small $\da$:
\begin{align}
|L\hat{T}^i|&\les|\hat{T}^j||X^A\fai_j|+|X^A\eta|\les\da M,\\
|\sd\hat{T}^i|&\les|\ta|\les M.\label{sdhatTi}
\end{align}
\end{lem}
\begin{pf}
The computation for \eqref{LTi} will be given and the proofs for rest are the same and thus omitted. Since $\hat{T}$ is $\Si_t$-tangential, one can expand $L\hat{T}$ as $L\hat{T}=p_L\hat{T}+q_L^AX_A$, for some functions $p_L$, $q_L^A$ to be determined later. Then since, $L\hat{T}=D_L\hat{T}-\Ga_{\a\be}^{\ga}
L^aT^{\be}\pa_{\ga}$, it follows from Remark\ref{christoffelsymbol} that
\bee
p_L=g(L\hat{T},\hat{T})=g(D_L\hat{T},\hat{T})-\Ga_{\a\be\ga}L^{\a}\hat{T}^{\be}\hat{T}^{\ga}=0,
\ee
and
\begin{align}
q_L^A\sg_{AB}&=g(L\hat{T},X_B)=g(D_L\hat{T},\hat{T})-\Ga_{\a\be\ga}L^{\a}\hat{T}^{\be}X_B^{\ga}\\
&=\kappa^{-1} g(D_LT,X_B)=-\kappa^{-1}\zeta_B\\
&=-\eta\ep_B+X_B(\eta)=\hat{T}^i(X_B(\fai_i))+X_B(\eta).
\end{align}
\end{pf}
\begin{lem}\label{ltfe}
The following estimates hold for sufficiently small $\da$:
\begin{align}
||LT\fe||_{\supnormda}&\les \da M^2,\\
||T^2\fe||_{\supnormda}&\les\da^{-1}M,\label{T^2fe}\\
||\sd T\fe||_{\supnormda}&\les M.
\end{align}
Meanwhile, the following estimates on $\mu$ hold:
\bee
||\sd\mu||_{\supnormda}\les M,\ ||T\mu||_{\supnormda}\les\da^{-1}M.
\ee
\end{lem}
\begin{pf}
Since $LT\fe=L(\eta^{-1}\mu)\hat{T}^i\fai_i+\eta^{-1}\mu L(\hat{T}^i\fai_i)$, it holds that
\beeq
|LT\fe|\les|L(\eta^{-1}\mu)\fai_i|+|\mu||\fai_i L(\hat{T}^i)+L(\fai_i)|\les\da M^2.
\eeq
\hs Suppose that $|\sd T\fe|\les M$. Then $|\sd Th|=|\sd T\fai_0-\sd(\fai_iT\fai_i)
+a\sd T\fe|\les M$. Since in the acoustical coordinates, $\s{d}_Af=X_Af$ for any function $f$, which implies
 \bee
 \s{d}_ALf=LX_Af+\s{d}f\cdot [X_A,L]=L\sd_Af.
 \ee
 Then, Commuting $\sd$ with \eqref{transportmu} yields
 \bee\label{lsdmu}
 L\sd\mu=\sd m+\mu\sd e+e\sd\mu.
 \ee
 \hs For $\sd m$, it holds that $|\sd m|=|\dfrac{1}{2}\sd\left(\dfrac{dH}{dh}\right)Th+\dfrac{1}{2}\dfrac{dH}{dh}\sd Th|+|a\sd T\fe|\les M$.\\ 
 \hs Note that $|\sd e|\les|\sd Lh|+|\sd T^i||L\fai_i|+|\sd L\fai_i|$. Due to \eqref{sdhatTi}, $|\sd L\fe|=|\sd\fai_0-[\sd(\eta^{-1}\hat{T}^i)\fai_i+\eta^{-1}\hat{T}^i\sd\fai_i]
-\fai_i\sd\fai_i|\leq\da M^2$, then $|\sd Lh|\les \da M^2$, which implies $|\sd e|\les\da M^2$. Thus, integrating \eqref{lsdmu} along integral curves of $L$ yields $||\sd\mu||_{\supnormda}\les M$. As a consequence,
\beeq
|\sd T\fe|&=&|\sd(\eta^{-1}\mu\hat{T}^i\fai_i)|=|\sd(\eta^{-1}\mu)\hat{T}^i\fai_i+\eta^{-1}\mu\sd\left(\hat{T}^i\fai_i\right)|\\
&\les&\da M^2<M,
\eeq
for $\da$ sufficiently small, which recovers the assumption $|\sd T\fe|\les M$. As a corollary, one obtains the bound $||\sd\mu||_{\supnormda}\les M$. The same argument can be applied to prove $||T^2\fe||_{\supnormda}\les\da^{-1}M$ and at the same time, one can obtain the estimate $||T\mu||_{\supnormda}\les\da^{-1} M$. To see this, commuting $T$ with \eqref{transportmu} yields
\bee\label{LTmu}
LT\mu=Tm+(Te)\mu+(T\mu)e-(\eta^A+\zeta^A)X_A\mu=Tm+eT\mu+O(\da M^2).
\ee
Under the assumption $|T^2\fe|\les\da^{-1}M$, it holds that $|Tm|\les |T^2\fai_0|+|aT^2\fe|\les\da^{-1}M$. Then, integrating \eqref{LTmu} along integral curves of $L$ yields $|T\mu|\les\da^{-1}M$ and as a consequence, $|T^2\fe|=|T(\eta^{-1}\mu \hat{T}^i\fai_i)|\les M^2<\da^{-1}M$.
\end{pf}
\subsection{\textbf{Accurate estimates for $\mu$ and its derivatives}}
\hs Considering the inhomogeneous covariant wave equation $\square_{\til{g}}\fai_{\a}=-\dfrac{2\eta'}{\eta^2}a\fai_{\a}\de\fe+\dfrac{a}
{\p\eta}\left(\dfrac{\pa\fai_{\a}}{\pa t}-\fai_i\dfrac{\pa\fai_{\a}}{\pa x^i}\right)$ and noticing that \eqref{decomposewaveeq} can be written as
\bee\label{decompose}
\bes
\mu\Omega\square_{\til{g}}f&=\mu\s{\de}f-L(\dl f)-\dfrac{1}{2}\left(
(\til{tr}\chi+L\log\Omega)\dl f+(\til{tr}\dc+\dl\log\Omega) Lf\right)\\
&-2\zeta\cdot\sd f+\Omega^{-1}\dfrac{d\Omega}{dh}\sd h\cdot\sd f
\end{split}
\ee
yield the following transport equation for $\dl\fai_{\a}$:
\bee\label{transportdlfai}
\bes
L(\dl \fai_{\a})+\dfrac{1}{2}\til{tr}\chi\dl \fai_{\a}&=\underbrace{\mu\s{\de}\fai_{\a}-\frac{1}{2}L(\log\Omega)\dl\fai_{\a}-\dfrac{1}{2}(\til{tr}\dc+\dl\log\Omega)L\fai_{\a}-2\zeta\cdot\sd \fai_{\a}+\Omega^{-1}\dfrac{d\Omega}{dh}\sd{h}\cdot\sd \fai_{\a}}_{I}\\
&\underbrace{-\dfrac{\mu a}{\eta^2}\left(\dfrac{\pa\fai_{\a}}{\pa t}-\fai_i\dfrac{\pa\fai_{\a}}{\pa x^i}\right)+\dfrac{2\p\mu\eta'}{\eta^3}a\fai_{\a}\de\fe}_{II}.
\end{split}
\ee
It follows from Lemma\ref{estmu},\ref{2ndff}, the bootstrap assumptions \eqref{bs}, Proposition\ref{relationangular} and the relation $tr\dc=\eta^{-2}\mu tr\chi+2\eta^{-1}\mu tr\ta$ that $I$ is bounded by $\da M^2$. It remains to estimate $II$. For $\mu\de\fe$, it holds that
\bee
\bes
|\mu\de\fe|&=|\mu\da^{ij}\frac{\pa\fai_j}{\pa x^i}|=|\mu(\hat{T}^i\hat{T}^j+\sg^{AB}X_A^iX_B^j)\frac{\pa\fai_j}{\pa x^i}|\\
&\les|\hat{T}^iT\fai_i|+\mu|\sd\fai|\les M.
\end{split}
\ee
\hs Since $\mu\dfrac{\pa\fai_{\a}}{\pa t}-\mu\fai_i\dfrac{\pa\fai_{\a}}{\pa x^i}=\mu L\fai_{\a}+\eta^{2}T\fai_{\a}+\mu\fai_i\dfrac{\pa\fai_{\a}}{\pa x^i}-\mu\fai_i\dfrac{\pa\fai_{\a}}{\pa x^i}=\frac{1}{2}\eta^{2}\dl\fai_{\a}+\frac{1}{2}\mu L\fai_{\a}$, then \eqref{transportdlfai} can be written as
\bee\label{transportdlfai2}
L(\dl\fai_{\a})+\dfrac{1}{2}\left(\til{tr}\chi+a\right)\dl\fai_{\a}=O(\da M^2).
\ee
For convenience, we may replace $\a$ by $2a$ in the analysis.
\begin{prop}\label{accruatemu1}
For sufficiently small $\da$, it holds that
\bee\label{aestlmu}
|e^{at}|t| L\mu(t,u,\ta)-2e^{-2a} L\mu(-2,u,\ta)|\les\da M^2,
\ee
\bee\label{aestmu}
\arrowvert\mu(t,u,\ta)-1-2A_1(t)L\mu(-2,
u,\ta)\arrowvert\les\da M^2,
\ee
where
\[
A_1(t)=\int_{-2}^t\dfrac{1}{e^{a(\tau+2)}(-\tau)}\ d\tau.
\]
\end{prop}

\begin{pf}
\hs Since $\til{tr}\chi=\til{tr}\chi'-\dfrac{2}{u-t}$, where $\til{tr}\chi'$ is bounded by $\da M$ due to Lemma\ref{2ndff}, it follows from\eqref{transportdlfai2} that
\bee
L(\dl\fai_{\a})+(a-\dfrac{1}{u-t})\dl\fai_{\a}=O(\da M^2),\\
\ee
which implies
\bee
L(e^{a t}(u-t)\dl\fai_{\a})=O(\da M^2).\label{Leatu-tdlfai}
\ee
\hs Integrating \eqref{Leatu-tdlfai} along integral curves of $L$ and using the bootstrap assumptions yield
\bee\label{estdlfai}
\bes
&|e^{at}(u-t)\dl\fai_{\a}(t,u,\ta)-e^{-2a}(u+2)\dl\fai_{\a}(-2,u,\ta)|\les\da M^2,\\
&|e^{at}|t|\dl\fai_{\a}(t,u,\ta)-2e^{-2a}\dl\fai_{\a}(-2,u,\ta)|\les\da M^2,\\
&|e^{at}|t| T\fai_{\a}(t,u,\ta)-2e^{-2a} T\fai_{\a}(-2,u,\ta)|\les\da M^2,\\
&|e^{at}|t|\fai_{\a}(t,u,\ta)-2e^{-2a}\fai_{\a}(-2,u,\ta)|
\les\da M^2,
\end{split}
\ee
where the last inequality follows by integrating third one along integral curves of $T$ from $0$ to $u$. It follows from $|T\fe|\les\da M^2$ that
\bee\label{esttfe}
|e^{at}|t|T\fe(t,u,\ta)-2e^{-2a} T\fe(-2,u,\ta)|\les |T\fe(t,u,\ta)|+|T\fe(-2,u,\ta)|\les\da M^2.
\ee
\hs It follows from \eqref{transportmu}, Lemma\ref{esthme} and \ref{estmu} that
\bee\label{lmu-2t}
\bes
e^{at}(-t)L\mu(t,u,\ta)-2e^{-2a} L\mu(-2,u,\ta)&=e^{at}(-t)m(t,u,\ta)-2e^{-2a}m(-2,u,\ta)+O(\da M^2),
\end{split}
\ee
It follows from \eqref{mcomputation}, Lemma\ref{hrhoeta}
, \eqref{estdlfai} and \eqref{esttfe} that
\bee\label{m-2t}
\bes
&|e^{at}(-t)m(t,u,\ta)-2e^{-2a}m(-2,u,\ta)|=|e^{a\tau}(-\tau)[(-1-\eta\eta')Th+aT\fe]|^t_{-2}|\\
=&|e^{a\tau}(-\tau)[(-1-\eta'\eta)(T\fai_0-\fai_iT\fai_i+aT\fe)+aT\fe]|^t_{-2}|\les\da M^2.
\end{split}
\ee
The estimate \eqref{m-2t} and \eqref{lmu-2t} complete the proof for \eqref{aestlmu}. \eqref{aestmu} will follow from \eqref{aestlmu} as follows.
\bee
\bes
\mu(t,u,\ta)&=\mu(-2,u,\ta)+\int_{-2}^t\dfrac{e^{a\tau}(-\tau)}{e^{a\tau}(-\tau)}L\mu(\tau,u,\ta)d\tau\\
&=1+2L\mu(-2,u,\ta)\int_{-2}^t\dfrac{1}{e^{a(\tau+2)}(-\tau)}d\tau+O(\da M^2).\\
\end{split}
\ee
\end{pf}
\subsubsection{\textbf{Two key properties of $\mu$ near shock}}
\hs Define the shock region as: $W_{shock}=\{(t,u,\ta)|\mu(t,u,\ta)\leq\dfrac{1}{10}\}$. The following proposition implies once a point $p$ lies in $W_{shock}$, then all the points after $p$ along integral curves of $L$ will stay in $W_{shock}$ as well.
\begin{prop}\label{keymu1}
For sufficiently small $\da$ and all $(t,u,\ta)\in W_{shock}$, it holds that
\bee
L\mu(t,u,\ta)\leq-\dfrac{1}{C(a)}\dfrac{1}{e^{at}(-t)}\les-1,
\ee
where $C>0$ is a constant which may depend on $a$.
\end{prop}
\begin{remark}
Since $||\mu e||_{\supnormda}\les\da M^2$, it follows from Proposition\ref{keymu1} that in the shock region $m\les-1$, which implies that $|Th|\geq c>0$ for some constant $c$ due to \eqref{mcomputation} and the estimate $|T\fe|\les\da M^2$.
On the other hand,
\bee
\bes
|Th|&=|T\fai_0-\fai_iT\fai_i+a\eta^{-1}\mu\hat{T}^i\fai_i|=\eta^{-1}\mu|\hat{T}\fai_0-\fai_i\hat{T}\fai_i+a\hat{T}^i\fai_i|,
\end{split}
\ee
which implies that
\bee
|\hat{T}\fai_0-\fai_i\hat{T}\fai_i+a\hat{T}^i\fai_i|\geq\mu^{-1}c.
\ee
Since $|\hat{T}|=1$ and $|\fai_i|\les\da M$, we deduce that
\bee
\text{when a shock forms,}\ \hat{T}\fai_{\a}\ \text{blows up for some $\a\in\{0,1,2,3\}$}.
\ee
\hs More precisely, since $\hat{T}\fai_0=\hat{T}^i\pa_i\fai_0=\hat{T}^i\pa_0\fai_i=\hat{T}^i\left(L\fai_i+
\eta\hat{T}^j\pa_j\fai_i+\fai_j\pa_j\fai_i\right)$ and $|L\fai_i|\les\da M$, it follows that \textbf{when a shock forms, $\pa_i\pa_j\fe$ blows up, for some $i,j\in\{1,2,3\}$}. That is, for the compressible Euler equations, when a shock forms, $\nabla v$ and $\nabla\p$ blow up.
\end{remark}
\begin{remark}
Here we illustrate how the largeness of the initial data leads to the shock formation in finite time. Since $m=\dfrac{1}{2}\dfrac{dH}{dh}Th+aT\fe=(-1-\eta\eta')(T\fai_0-\fai_iT\fai_i+aT\fe)+aT\fe$,
and by the estimate \eqref{estdlfai}, \eqref{esttfe}, and $|\eta-\eta_0|+|\eta'-\eta'_0|\les\da M$, it is clear that up to the order $O(\da M^2)$,
$T\fai_{\a}(t,u,\ta)$, $\fai_i(t,u,\ta)$ and $T\fe(t,u,\ta)$ can be replaced by $\dfrac{2e^{-2a}}{|t|e^{at}}T\fai_{\a}(-2,u,\ta)$, $\dfrac{2e^{-2a}}{|t|e^{at}}\fai_i(-2,u,\ta)$ and $\dfrac{2e^{-2a}}{|t|e^{at}}T\fe(-2,u,\ta)$ respectively, so that
\bee
\bes
m&=2\left(\dfrac{2e^{-2a}}{|t|e^{at}}\right)^2\fai_i T\fai_i(-2,u,\ta)-2\dfrac{2e^{-2a}}{|t|e^{at}}T\fai_0(-2,u,\ta)
-3\dfrac{2e^{-2a}}{|t|e^{at}}aT\fe(-2,u,\ta)+O(\da)\\
&=-\dfrac{4e^{-2a}}{|t|e^{at}}\left(\pa_r\fe_1\left(\dfrac{r-2}{\da},\ta\right)\right)+O(\da).
\end{split}
\ee
Hence, $L\mu=m+\mu e=-\dfrac{4e^{-2a}}{|t|e^{at}}\left(\pa_r\fe_1\left(\dfrac{r-2}{\da},\ta\right)\right)+O(\da)$, which can be integrated to obtain
\bee\label{mushock}
\bes
\mu(t,u,\ta)-\mu(-2,u,\ta)&=\int_{-2}^tL\mu(\tau,u,\ta)d\tau=4\pa_r\fe_1\left(\frac{r-2}{\da},\ta\right)\int_{-2}^t\dfrac{e^{-a(\tau+2)}}{\tau}d\tau+O(\da).
\end{split}
\ee
Let
\bee\label{a^ast}
a^{\ast}=4\int_{-2}^{\sigma}\dfrac{e^{-a(\tau+2)}}{\tau}d\tau.
 \ee
 Then, it follows that for fixed constant $a$:
\been
\item if $\max\pa_s\fe_1=c\geq-\frac{2}{a^{\ast}}>0$, then for all $a\neq 0$, there exists a $T_{\ast}(a)\in[-2,\sigma)$ such that $\mu(T_{\ast})=0$, which means that a shock forms at $t=T_{\ast}$ for sufficiently small $\da$. Moreover, $T_{\ast}(a)$ can be computed as
\bee\label{Tast}
c\int_{-2}^{T_{\ast}}\dfrac{e^{-a(\tau+2)}}{\tau}\ d\tau=-1+O(\da).
\ee
It follows that $T_{\ast}(a)$ is a increasing function of $a$. That is, compared with the undamped case, if $a>0$, then the damping term will delay the shock formation, while for $a<0$, the anti-damping term will lead to the formation of shock in advance.
\item If $\max\pa_s\fe_1=c'\leq-\frac{1}{2a^{\ast}}$, then it follows from \eqref{mushock} that $\mu>0$ on $t\in[-2,t^{\ast}]$, so that one can obtain a smooth solution on $[-2,t^{\ast}]$.
\een
Moreover, if $\dfrac{dH}{dh}=0$, i.e. the flow is the Chaplygin gas, then $\mu>0$ on $t\in[-2,t^{\ast}]$.
\end{remark}
\begin{pf}
By Proposition\ref{accruatemu1}, it holds that
\begin{align}
-2e^{-2a}\int_{-2}^t\dfrac{e^{-a\tau}}{(-\tau)}d \tau L\mu(-2,u,\ta)&\geq(1-\mu(t,u,\ta)
-O(\da M^2)),
\end{align}
which implies that
\begin{align}
2e^{-2a}L\mu(-2,u,\ta)&\leq-\dfrac{1}{C(a)}\underbrace{(1-\mu(t,u,\ta)-O(\da M^2))}_{I},\label{proofkeymu2}
\end{align}
where the term $I$ in the shock region is $\geq\dfrac{1}{2}$ for sufficiently small $\da$ and $C(a)>0$. It follows from this and Proposition\ref{accruatemu1} that
\bee
\bes
e^{at}(-t)L\mu(t,u,\ta)&\leq 2e^{-2a}L\mu(-2,u,\ta)+O(\da M^2)\\
&\les-\dfrac{1}{C}.
\end{split}
\ee
\end{pf}
\hs In the energy estimates later on, we will encounter the integral $\int_{-2}^t\mu^{-1}T\mu$ which may not be bounded. However, the positive part of $\mu^{-1}T\mu$ is integrable in time as shown in the following proposition.
\begin{prop}\label{keymu2}
Let $s^{\ast}$ be given by \eqref{defsasttast}. Then for sufficiently small $\da$ and all $(t,u,\ta)\in W_{shock}$, it holds that
\bee\label{mu-1tmu+}
(\mu^{-1}T\mu)_+\les\dfrac{\da^{-1}}{|e^{-at}-e^{-as^{\ast}}|^{\frac{1}{2}}}.
\ee
\end{prop}
\begin{pf}
 Choose a new coordinate system $(\wi{\ta^1},\wi{\ta^2})$ on $S_{t,u}$ such that $T=\dfrac{\pa}{\pa u}\mid_{t,u,\wi{\ta}}$. Then $\wi{\ta}$ will be constant along integral curves of $T$. Consider the following function $f$ along an integral curve of $T$: $f(u)=T(\log\mu)(u),u\in[0,\tilde{\da}].$ Then $f$ will attain its maximum at $u=u^{\ast}\in[0,\tilde{\da}].$ Without loss of generality, one may assume $f(u^{\ast})>0,$ otherwise $(\mu^{-1}T\mu)_+=0$ and then \eqref{mu-1tmu+} holds trivially. Since when $u=0$, $\mu$ is constant, which implies $f(u^{\ast})=0$.  Then, $\dfrac{d}{du}f(u^{\ast})\geq 0$, i.e. $\mu^{-1}T^2\mu-(\mu^{-1}T\mu)^2\geq 0,$ which implies 
\[
(\mu^{-1}T\mu)_+\leq\sqrt{\dfrac{||T^2\mu||_{L^{\infty}([0,\tilde{\da}])}}{\inf_{[0,\tilde{\da}]}\mu}}.
\]
It remains to bound $T^2\mu$ and $\mu.$ Following the same argument in the proof of 
Proposition\ref{accruatemu1}, one can obtain
\bee\label{T^3fai}
||t|e^{at}T^{3}\fai_0(t,u,\ta)-2e^{-2a}T^3\fai_0(-2,u,\ta)|\leq\da^{-1}M^2.
\ee
In the proof of Lemma\ref{ltfe}, we have verified actually that $|T^2\fe|\les M^2$, and the similar argument for \eqref{T^2fe} implies 
\bee\label{T^3fe}
|T^3\fe|\les\da^{-1}M^2.
\ee
 Then commuting $T^2$ with $L\mu=m+\mu e$ yields
\bee\label{lt2muequation}
LT^2\mu=T^2m+(T^2\mu)e+2(T\mu)(Te)+\mu(T^2e)+[L,T^2]\mu.
\ee
Then, it follows from \eqref{T^3fai}-\eqref{lt2muequation} and the similar argument in proving Proposition\ref{accruatemu1} that, 
\bee\label{LT^2mu}
|e^{at}|t|LT^2\mu(t,u,\ta)-2e^{-2a}LT^2\mu(-2,u,\ta)|=O(\da^{-1}M^2),
\ee
which implies
\bee
|T^2\mu|\les\da^{-2},
\ee
by integrating \eqref{LT^2mu} along integral curves of $L$.\\
\hs It follows from  Proposition\ref{keymu1} and \ref{accruatemu1} that for all $(t,u,\ta)\in W_{shock}$
\bee
 2e^{2a}L\mu(-2,u,\ta)\les|t|e^{at}L\mu(t,u,\ta)+O(\da M^2)\les-1/C+O(\da M^2)\les -1/C,
 \ee
 for sufficiently small $\da$.\\
 Hence, for $t<s^{\ast}$,
\beeq
\mu(t,u,\ta)&\geq&-\int_{t}^{s^{\ast}}\dfrac{e^{a\tau}(-\tau)}{e^{a\tau}(-\tau)}L\mu(\tau,u,\ta)d\tau\\
&=&-\int_{t}^{s^{\ast}}\dfrac{1}{e^{a\tau}(-\tau)}(2e^{2a}L\mu(-2,u,\ta)+O(\da M^2))d\tau\\
&\geq&(\dfrac{1}{C}+O(\da M^2))\int_t^{s^{\ast}}\dfrac{1}{e^{a\tau}(-\tau)}\\
&\geq&\left(\dfrac{1}{2C}+O(\da M^2)\right)\dfrac{1}{a}\left(e^{-at}-e^{-as^{\ast}}\right)\gtrsim|e^{-at}-e^{-as^{\ast}}|,
\eeq
for sufficiently small $\da$. Therefore, $(\mu^{-1}T\mu)_+\les\dfrac{\da^{-1}}{|e^{-at}-e^{-as^{\ast}}|^{\frac{1}{2}}}$.
\end{pf}
\section{\textbf{Estimates for deformation tensors}}\label{section4}
\hs The deformation tensor associated with a vector filed $Z$ with respect to $g$ is defined as: $\zpi_{\a\be}=D_{\a}Z_{\be}+D_{\be}Z_{\a}$ or equivalently $\zpi_{\a\be}=\mathcal{L}_{Z}g_{\a\be}$ and with respect to the conformal metric $\til{g}$ is given by
\bee\label{relationdeformationtensor}
\zgpi_{\a\be}=\mathcal{L}_{Z}\til{g}_{\a\be}=\Omega\zpi_{\a\be}+Z(\Omega)g_{\a\be}.
\ee
Note that for any vector fields $X$ and $Y$, it holds that
\bee
\zpi_{XY}=g(D_XZ,Y)+g(D_YZ,X).
\ee
Therefore, by Lemma\ref{connection1} and \ref{connection2}, one can compute the deformation tensors associated with the commutation vectorfields $\{Z_i\}$ with $i=1,2,3,4,5$.\\
 \hs For $Z_1=Q=tL$, the deformation tensor $\qpi$ can be computed in the null frame $(L,\dl,X_1,X_2)$ as:
\bee\label{computationqpi}
\bes
&\qpi_{LL}=0,\hs \qpi_{\dl\dl}=4t\mu L(\eta^{-2}\mu)-4\mu\eta^{-2}\mu,\\
&\qpi_{L\dl}=-2tL\mu-2\mu,\hs \qpi_{LA}=0,\\
&\qpi_{\dl A}=2t(\eta_A+\zeta_{A}),\hs \sqpi_{AB}=2t\chi_{AB}.
\end{split}
\ee
Hence, the corresponding deformation tensor with respect to $\til{g}$ is given by:
\beeq
&&\qgpi_{LL}=0,\hs \qgpi_{\dl\dl}=4\Omega\mu(tL(\eta^{-2}\mu)-\eta^{-2}\mu),\\
&&\qgpi_{L\dl}=-2(\Omega t L\mu+\Omega\mu+\mu tL\Omega),\hs
\qgpi_{LA}=0,\\
&&\qgpi_{\dl A}=2\Omega t(\eta_A+\zeta_A),\hs \pre {(Q)} {} {\hat{\til{\s{\pi}}}}_{AB}=
2\Omega t\hat{\chi}_{AB},\\
&&tr\sqgpi=2\Omega t\cdot\til{tr}\chi.
\eeq
\hs For $Z_2=T$, the deformation tensor $\tpi$ with respect to $g$ in the null frame $(L,\dl,X_1,X_2)$ is given by:
\bee\label{computationtpi}
\bes
&\tpi_{LL}=0,\hs \tpi_{\dl\dl}=2\mu T(\eta^{-2}\mu),\\
&\tpi_{L\dl}=-2 T\mu,\hs \tpi_{LA}=-(\eta_A+\zeta_A),\\
&\tpi_{\dl A}=-\eta^{-2}\mu(\eta_A+\zeta_A),\hs \stpi_{AB}=2\eta^{-1}\mu\ta_{AB},
\end{split}
\ee
and with respect to $\til{g}$ is given by:
\beeq
&&\tgpi_{LL}=0,\hs \tgpi_{\dl\dl}=2\Omega\mu T(\eta^{-2}\mu),\\
&&\tgpi_{L\dl}=-2(\Omega T\mu+\mu T\Omega),\hs
\tgpi_{LA}=-\Omega(\eta_A+\zeta_A),\\
&&\tgpi_{\dl A}=-\Omega\eta^{-2}\mu(\eta_A+\zeta_A),\hs \pre {(T)} {} {\hat{\til{\s{\pi}}}}_{AB}=
2\Omega\eta^{-1}\mu\hat{\ta}_{AB},\\
&&tr\stgpi=2\Omega\eta^{-1}\mu\til{tr}\ta,
\eeq
where $\eta^{-1}\mu\til{tr}\ta=\Omega^{-1}(\Omega\eta^{-1}\mu tr\ta+T\Omega)$.\\
\hs Since $R_i$ is $\Si_t$ tangent, it's convenient to compute $\rpi$ in the frame $(L,T,X_1,X_2)$ (since $T$ is $\Si_t$ tangent) instead of in the null frame. It follows from the fact that $L$ and $T$ are $g$-orthogonal to $S_{t,u}$ and Lemma\ref{connection1} that
\bee\label{computationrpi1}
\bes
&\rpi_{LL}=2g(D_LR_i,L)=-2g(R_i,D_LL)=0,\\
&\rpi_{TT}=2g(D_TR_i,T)=-2g(R_i,D_TT)=2\eta^{-1}\mu R_i(\eta^{-1}\mu),\\
&\rpi_{LT}=g(D_LR_i,T)+g(D_TR_i,L)=-g(R_i,D_LT)-g(R_i,D_TL)=-R_i\mu.
\end{split}
\ee
\hs Due to the definition $R_i=\Pi \mathring{R}_i$ and the fact that for any vector field $Z$: $g(\Pi Z,X_A)=g(Z,X_A)$, it holds that
\bee\label{rpiLATAABcompute}
\bes
\rpi_{LA}&=g((D_L\Pi)\mari{R}_i,X_A)+g(D_L\mari{R}_i,X_A)-g(R_i,D_AL),\\
\rpi_{TA}&=g((D_T\Pi)\mari{R}_i,X_A)+g(D_T\mari{R}_i,X_A)-g(R_i,D_AT),\\
\rpi_{AB}&=g((D_A\Pi)\mari{R}_i,X_B)+g((D_B\Pi)\mari{R}_i,X_A)\\
&+g(D_A\mari{R}_i,X_B)+g(D_B\mari{R}_i,X_A).
\end{split}
\ee
 Decompose $\mari{R}_i$ as
\[
\mari{R}_i=R_i^AX_A+\lam_i\hat{T},
\]
for some function $\lam_i$. Then, it follows that
\[
g((D_Z\Pi)\mari{R}_i,X_A)=\lam_i\eta\mu^{-1} g((D_Z\Pi)T,X_A)=
-\lam_i\eta\mu^{-1} g(D_ZT,X_A).
\]
Hence,
\bee\label{gDLTAPi}
g((D_L\Pi)\mari{R}_i,X_A)=\lam_i\eta\mu^{-1}\zeta_A,\ g((D_T\Pi)\mari{R}_i,X_A)=\lam_iX_A(\eta^{-1}\mu),\ g((D_A\Pi)\mari{R}_i,X_B)=-\lam_i\ta_{AB}.
\ee
\hs For $D_L\mari{R}_i$, it holds that
\bee\label{gDLRi}
\bes
g(D_L\mari{R}_i,X_A)&=L^{\lam}\left(\dfrac{\pa\mari{R}_i^{\a}}{\pa x^{\lam}}+\Ga^{\a}_{\lam\gamma}\mari{R}_i^{\gamma}\right)X_A^{\be}g_{\a\be}=L^{\lam}\left(g_{\a\be}\dfrac{\pa\mari{R}_i^{\a}}{\pa x^{\lam}}+\Ga_{\lam \be\gamma}\mari{R}_i^{\gamma}\right)X_A^{\be}\\
&=L^{\lam}\left(\bar{g}_{kl}\dfrac{\pa\mari{R}_i^k}{\pa x^{\lam}}\right)X_A^l=L^{\lam}\da_{kl}\ep_{ijk}\da_{j\lam}X_A^l=L^{k}\ep_{ik l}X_A^l.
\end{split}
\ee
Similarly, it holds that
\bee\label{gDTARi}
g(D_T\mari{R}_i,X_A)=T^{k}\ep_{ik l}X_A^l,\  g(D_A\mari{R}_i,X_B)=X_A^{k}\ep_{ik l}X_B^l.
 \ee
 Collecting \eqref{rpiLATAABcompute}-\eqref{gDTARi} yields
\bee\label{rpiLATAAB}
\bes
\rpi_{LA}&=\lam_i\eta\mu^{-1}\zeta_A+L^{\lam}\ep_{i\lam l}X_A^l-R_i^B\chi_{AB},\\
\rpi_{TA}&=\lam_i X_A(\eta^{-1}\mu)+T^{\lam}\ep_{i\lam l}X_A^l-\eta^{-1}\mu R_i^B\ta_{AB},\\
\rpi_{AB}&=-2\lam_i\ta_{AB}+X_A^{\lam}\ep_{i\lam l}X_B^l+X_B^{\lam}\ep_{i\lam l}X_A^l=-2\lam_i\ta_{AB}.
\end{split}
\ee
\hs Since 
$|\ep_A|\les|\s{k}|\les\da M$ and
\bee\label{estzetaetaA}
|\eta\mu^{-1}\zeta_A|=|\eta\ep_A-\sd_A\eta|\les\da M,\ |\eta_A|=|\zeta_A+\sd_A\mu|\les M,
\ee
due to Lemma\ref{ltfe}, the following estimates hold for the deformation tensor $\qgpi$:
\bee\label{estqgpi}
\bes
&\qgpi_{LL}=0,\hs ||\mu^{-1}\qgpi_{\dl\dl}||_{\supnormda}\les M,\hs
||\qgpi_{L\dl}||_{\supnormda}\les M,\\
&\qgpi_{LA}=0,\hs ||\mu^{-1}\qgpi_{\dl A}||_{\supnormda}\les M,\\
&||\pre {(Q)} {} {\hat{\til{\s{\pi}}}}_{AB}||_{\supnormda}\les\da M,\hs ||tr\sqgpi||_{\supnormda}\les 1.
\end{split}
\ee
\hs The following estimates for the deformation tensor $\tgpi$ hold:
\bee\label{esttgpi}
\bes
&\tgpi_{LL}=0,\hs ||\mu^{-1}\tgpi_{\dl\dl}||_{\supnormda}\les\da^{-1} M,\hs
||\tgpi_{L\dl}||_{\supnormda}\les\da^{-1} M,\\
&||\tgpi_{LA}||_{\supnormda}\les M,\hs ||\mu^{-1}\tgpi_{\dl A}||_{\supnormda}\les M,\\
&||\pre {(T)} {} {\hat{\til{\s{\pi}}}}_{AB}||_{\supnormda}\les\da M,\hs ||tr\stgpi||_{\supnormda}\les 1.
\end{split}
\ee
\hs Define
\[
y^i=\hat{T}^i-\dfrac{x^i}{u-t},
\]
which measure the difference between $\hat{T}^i$ and $\pa_r^i$ in the acoustical coordinates. Since $x^{k}\ep_{ik l}X_A^l=\sum_l\mari{R}_i^lX_A^l=\bar{g}(\mari{R}_i,X_A)=\sg_{AB}R_i^B$, then
\bee
\rpi_{TA}=\lam_i X_A(\eta^{-1}\mu)+\eta^{-1}\mu y^{k}\ep_{ik l}X_A^l-\eta^{-1}\mu R_i^B\ta'_{AB}.
\ee
Similarly, define
\[
z^i=L^i+\dfrac{x^i}{u-t},
\]
which measure the difference between $L^i$ and $\pa_t^i-\pa_r^i$ in the acoustical coordinates. Then, 
\[
z^i=-\eta y^i-\dfrac{(\eta-1)x^i}{u-t}-\fai_i.
\]
One can rewrite $\rpi_{LA}$ as
\[
\rpi_{LA}=\lam_i\eta\mu^{-1}\zeta_A+z^{k}\ep_{ik l}X_A^l-R_i^B\chi'_{AB}.
\]
It suffices to estimate $\lam_i$ and $y^i$.\\
\hs Since $\lam_i=\bar{g}(\mari{R}_i,\hat{T})=\sum_k\mari{R}_i^k\hat{T}^k$, it follows that
\[
L\lam_i=\sum_kL(\mari{R}_i^k)\hat{T}^k+\sum_kL(\hat{T}^k)\mari{R}_i^k.
\]
Then, it follows from $L(\mari{R}_i^k)=\ep_{ijk}L^j=-\ep_{ijk}(\eta\hat{T}^j+\fai_j)$, Lemma\ref{LTidTi} and \eqref{estzetaetaA} 
 that
\bee\label{Llami}
|L\lam_i|=|-\ep_{ijk}\fai_j\hat{T}^k-\eta\mu^{-1}\zeta^Ax^k\ep_{ikl}X_A^l|\les\da M.
\ee
Integrating \eqref{Llami} along integral curves of $L$ and noticing that $\lam_i=0$ on $\Si_{-2}$ yield
\bee\label{estlami}
||\lam_i||_{\supnormda}\les\da M.
\ee
\hs Define
\bee\label{defy'i}
y^{' i}=\hat{T}^i-\pa_r^i,
\ee
which measure the difference between $\hat{T}$ and the Euclidean radial derivative and will be estimated as follows.
Let $\Si$ be the projection from $\Si_t$ to the unit Euclidean sphere. Then, decompose $\hat{T}$ as
\[
\hat{T}=\bar{g}(\pa_r,\hat{T})\pa_r+\Si\hat{T}.
\]
It follows that
\begin{align}
1=|\hat{T}|^2&=|\Si\hat{T}|^2+|\bar{g}(\pa_r,\hat{T})|^2=r^{-2}\sum_i|\bar{g}(
\mari{R}_i,\hat{T})|^2+|\bar{g}(\pa_r,\hat{T})|^2=r^{-2}\sum_i\lam_i^2+
|\bar{g}(\pa_r,\hat{T})|^2,\label{1-g(parT)^2}\\
|y'|^2&=|\hat{T}|^2+1-2\bar{g}(\pa_r,\hat{T})=r^{-2}\sum_i\lam_i^2+|1-\bar{g}
(\pa_r,\hat{T})|^2.\label{1-g(parT)}
\end{align}
Since $|Tr|=|\eta^{-1}\mu\sum_i\dfrac{x^i\hat{T}^i}{r}|\les M$, then integrating $Tr$ along integral curves of $T$ from $0$ to $u$ and noticing that $r=-t$ when $u=0$ yield $|r+t|=|\int_0^uTrdu'|\les\da M$.
On $S_{t,0}$, $g(\pa_r,\hat{T})=1$, which implies that the angle between $\pa_r$ and $\hat{T}$ is less than $\dfrac{\pi}{2}$ due to the continuity of $\hat{T}$ in $u$ and the bootstrap assumptions. Then, $1-\bar{g}(\pa_r,\hat{T})\les1-|\bar{g}(\pa_r,\hat{T})|^2\les\da^2M^2$ due to \eqref{1-g(parT)^2} and \eqref{estlami}. 
 Thus, it follows from\eqref{1-g(parT)} that
\bee\label{esty'i}
|y'|\les|\lam_i|+|1-g(\pa_r,\hat{T})|\les\da M.
\ee
  Therefore,
 \bee
 |y^i-y^{' i}|\leq\arrowvert\dfrac{1}{r}-\dfrac{1}{u-t}\arrowvert|x^i|\les\arrowvert\dfrac{u-(r+t)}{r(u-t)}\arrowvert|x^i|\les\da M,
 \ee
  which implies
  \bee\label{estyi}
  |y^i|\les\da M.
  \ee
   Thus, it follows from \eqref{relationdeformationtensor}, \eqref{estlami} and \eqref{estyi} that
\bee\label{estrgpi}
\bes
&\rgpi_{LL}=0,\hs ||\mu^{-1}\rgpi_{\dl \dl}||_{\supnormda}\les M,\hs ||\rgpi_{L\dl}||_{
\supnormda}\les M,\\
&||\rgpi_{LA}||_{\supnormda}\les\da M,\hs ||\rgpi_{\dl A}||_{\supnormda}\les\da M^2,\\
&||\pre {(R_i)} {} {\hat{\til{\s{\pi}}}}_{AB}||_{\supnormda}\les\da^2 M^2,\hs ||tr\srgpi||_{\supnormda}\les \da M.
\end{split}
\ee
\hs At the end of this section, we introduce the following Proposition, whose proof can be found in\cite{christodoulou2014compressible}.
\begin{prop}\label{relationangular}
Under the bootstrap assumptions, there exists a numerical constant $C$ such that for any $S_{t,u}$ $1-$form $\xi$, the following bounds hold:
\bee\label{h0}
|\xi|^2\leq C\sum_i|\xi(R_i)|^2,
\ee
\bee\label{h1}
|\s{D}\xi|^2+|\xi|^2\leq C\sum_i|\lie_{R_i}\xi|^2,
\ee
for sufficiently small $\da$.
\end{prop}
\begin{remark}
In particular, taking $\xi=\sd f$ for any smooth function $f$ into \eqref{h0} and \eqref{h1} yields
\beeq
|\sd f|^2\leq C\sum_i|R_i f|^2,\quad |\s{D}^2f|^2+|\sd f|^2\leq  C\sum_i|\lie_{R_i}\sd f|^2.
\eeq
\end{remark}
\begin{remark}
Similarly, for any $S_{t,u}$ $(0,2)$ tensor $\lam$, \eqref{h0} and \eqref{h1} hold with $\xi$ replaced by $\lam$.
\end{remark}
\section{\textbf{Multiplier and commuting vector fields methods, estimates on lower order terms}}\label{section6}
\hs One can derive the energy estimates for the high order derivatives of $\fai$ as follows.
\been[i)]
\item The first part is based on the classical multiplier method. 
    A multiplier should be a non-space like vector field, which yields positive energy. Here, the following two vector fields will be used as multipliers:
     \bee\label{defmultipliers}
     K_{0}=\dl, K_{1}=L,
      \ee
      which are slightly different from the multipliers in \cite{christodoulou2014compressible}(since there is no need to control the growth in time). Then, contracting these multipliers with the energy-momentum tensor $\til{T}$ (see \eqref{energymomentum} later) and applying the divergence theorem (Lemma\ref{divergencethm}) yield the following fundamental energy estimates for the solution of $\Box_{\tilde{g}}\fai=f$:
    \bee\label{fee1}
    \left.\bes
    &E_{0}(t,u)+F_{0}(t,u)\les E_{0}(-2)+\int_{W^{u}_{t}}\Omega^{2}(-f\cdot K_{0}\fai-\dfrac{1}{2}\wi{T}^{\a\be}\wi{\pi}_{0,\a\be}):=E_{0}(-2)+Q_0,\\
    &E_{1}(t,u)+F_{1}(t,u)\les E_{1}(-2)+\int_{W^{u}_{t}}\Omega^{2}(-f\cdot K_{1}\fai-\dfrac{1}{2}\wi{T}^{\a\be}\wi{\pi}_{1,\a\be}):=E_1(-2)+Q_1,
    \end{split}\right.
    \ee
   where $E_{i}$ and $F_{i}$ $(i=0,1)$ are the energies and fluxes associated with $K_{i}$ respectively. However, \eqref{fee1} provides only the estimates for the first order derivatives of $\fai$, which is not enough to close the bootstrap assumptions. 
\item For high order energy estimates, another method will be adapted, which is called commuting vector fields method. 
    Commuting $Z_i^{\a}$ with the wave equation \eqref{nonlinearwave} where $Z_i$ is any one of the commutators $\{Q=tL, T,R_{i}\}$ and then applying \eqref{fee1} yield the energy (in)equality for $Z_i^{\a}\fai$. However, to obtain the estimates for $Z_i^{\a}\fai$, one has to deal with two additional terms in $Q_0$ and $Q_1$.
    \been[(1)]
    \item One major term arises from the product of the energy momentum tensor and the deformation tensors associated with $K_{0}$ and $K_{1}$, which can be bounded by energies and fluxes.  
    \item Another term arises from differentiating with the metric. That is, commuting $Z_i^{\a}$ ($|\a|=n-1$) with the wave equation \eqref{nonlinearwave} leads to $\square_{\wi{g}}Z_i^{\a}\fai=\p_{n}$, where $\p_n$ contains several terms due to the nonlinear structure of \eqref{nonlinearwave} and is given by:
        \begin{equation}\label{commutation1}
        \p_{n}=\dfrac{1}{\Omega^2}div_{g}\pre Z {} {J_{n-1}}+Z\p_{n-1}+\dfrac{1}{2}tr_{\wi{g}}\zgpi\cdot\p_{n-1},
        \end{equation}
        where $J_{n-1}$ is given by
        \[
        \pre Z {} {J_{n-1}^{\nu}}=(\zgpi^{\mu\nu}-\dfrac{1}{2}tr_{g}\zgpi\cdot g^{\mu\nu})\pa_{\mu}\fai_{n-1}.
        \]
        The divergence part in \eqref{commutation1} will give rise to the top order acoustical terms (that is, the terms of top order spatial derivatives for $\mu$ and $\chi$), which are the most difficult terms to deal with.
    \een
\een
\hs Let $\mu_{\sg}$ be the area form of $S_{t,u}$. Since $|\det g|=\mu^2\det \sg$, which implies that the volume form on $W_{t}^{ u}$ equals to $\mu dt d ud\mu_{\sg}$. Therefore, for convenience, the following notations will be adopted:
\beeq
&&\int_{\Si_{t}^{ u}}f=\int_{0}^{ u}\int_{S_{t, u}}f d\mu_{\sg}d u',\\
&&\int_{C_{ u}^{t}}f=\int_{-2}^{t}\int_{S_{t, u}}f d\mu_{\sg}dt',\\
&&\int_{W_{t}^{ u}}f=\int_{0}^{ u}\int_{-2}^t\int_{S_{t, u}}\mu\cdot f(t', u',\ta) d\mu_{\sg}d u'dt'.\\
\eeq
\subsection{\textbf{Multiplier method, fundamental energy estimates}}
\hs Let $f$ be the right hand side of\eqref{nonlinearwave} and then \eqref{nonlinearwave} becomes: $\square_{\til{g}}\fai=f$. 
Define the energy-momentum tensor associated with $\fai$ as:
\begin{equation}\label{energymomentum}
\begin{split}
\til{T}_{\mu\nu}:&=\pa_{\mu}\fai\pa_{\nu}\fai-\dfrac{1}{2}\til{g}_{\mu\nu}\til{g}^{\a\be}\pa_{\a}\fai\pa_{\be}\fai,\\
&=\pa_{\mu}\fai\pa_{\nu}\fai-\dfrac{1}{2}g_{\mu\nu}g^{\a\be}\pa_{\a}\fai\pa_{\be}\fai:=T_{\mu\nu},
\end{split}
\end{equation}
which describes the density of energy and momentum in spacetime. Direct computations yield
\beeq
\til{D}^{\mu}\til{T}_{\mu\nu}&=&\square_{\til{g}}\fai\pa_{\nu}\fai+\pa_{\mu}\fai\til{g}^{\mu\lam}\pa_{\lam}\pa_{\nu}\fai-\til{g}^{\a\be}\pa_{\nu}\pa_{\a}\fai\pa_{\be}\fai\\
&=&\square_{\til{g}}\fai\pa_{\nu}\fai=f(\pa_{\nu}\fai).
\eeq
In the null frame, due to \eqref{decomposeg2}, $\til{T}_{\a\be}$ can be computed as
\bee\label{tilTdecompose}
\bes
&\til{T}_{\dl\dl}=(\dl\fai)^2,\til{T}_{LL}=(L\fai)^2,\til{T}_{L\dl}=\mu|\s{d}\fai|^2,\til{T}_{LA}=L\fai\cdot X_{A}\fai,\\
&\til{T}_{\dl A}=\dl\fai\cdot X_{A}\fai,\til{T}_{AB}=X_{A}\fai\cdot X_{B}\fai-\dfrac{1}{2}\sg_{AB}(-\mu^{-1}\dl\fai L\fai+|\s{d}\fai|^2).
\end{split}
\ee
Define
\bee
\til{P}^{\mu}_{\a}:=-\til{T}^{\mu}_{\nu}K_{\a}^{\nu},\ \text{for}\ \a=0,1,\\
\ee
where $K_0$ and $K_1$ are the two multipliers defined by \eqref{defmultipliers}. 
Then,
\bee\label{DmuPmu}
\til{D}_{\mu}\til{P}^{\mu}_{\a}=-\left(f K_{\a}\fai+\dfrac{1}{2}\til{T}^{\mu\nu}\til{\pi}_{\a,\mu\nu}\right):=\til{Q}_{\a},\ \text{for}\ \a=0,1,\\
\ee
where $\til{\pi}_{\a,\mu\nu}$ are the deformation tensors associated with $K_{\a}$ $(\a=0,1)$ with respect to $\til{g}$.
Note that
\bee\label{tilDtilPtilQ}
\bes
\til{D}_{\mu}\til{P}_{\a}^{\mu}&=\dfrac{1}{\sqrt{|\det\til{g}|}}\pa_{\mu}\left(\sqrt{|\det\til{g}|}\til{P}^{\mu}\right)=\dfrac{1}{\Omega\sqrt{|\det g|}}\pa_{\mu}\left(\Omega^2\sqrt{|\det g|}\til{P}^{\mu}\right)=\Omega^{-2}D_{\mu}P^{\mu},
\end{split}
\ee
where $P^{\mu}=\Omega^{2}\til{P}^{\mu}$ and $Q=\Omega^2\til{Q}$. Then, it follows from\eqref{DmuPmu} that
\bee\label{DPQ}
D_{\mu}P_{\a}^{\mu}=Q_{\a},\ \a=0,1.
\ee
The following divergence theorem holds:
\begin{lem}\label{divergencethm}
For any spacetime vector field $J=J^{t}\dfrac{\pa}{\pa t}+J^{u}\dfrac{\pa}{\pa u}+\s{J}^A\dfrac{\pa}{\pa \ta^A}$, it holds that
\beeq
\int_{W_t^u}\mu D_{\a}J^{\a}dt'du'd\mu_{\sg}&=&\int_{\Si_t^u}\mu J^t du'd\mu_{\sg}-\int_{\Si_{-2}^u}\mu J^t du'd\mu_{\sg}\\
&+&\int_{C_u^t}\mu J^u dt'd\mu_{\sg}-\int_{C_{0}^t}
\mu J^u dt'd\mu_{\sg}.
\eeq
\end{lem}
 Since $P_{\a}^{\mu}=-\Omega^2\til{T}^{\mu}_{\nu}K_{\a}^{\nu}=- \Omega T^{\mu}_{\nu}K_{\a}^{\nu}, \a=0,1$, then
\bee\label{P0tP0uP1tP1u}
\bes
\mu P_0^t&=-\dfrac{1}{2}\left(P_{0,\dl}+\eta^{-2}\mu P_{0,L}\right)\\
&=-\dfrac{\Omega}{2}(-\til{T}_{\dl\dl}-\til{T}_{\dl L})=\dfrac{\Omega}{2}\left(\eta^{-2}\mu^2\sdfai^2+(\dl\fai)^2\right),\\
\mu P_0^u&=-P_{0,L}=\Omega T_{L\dl}=\Omega\mu\sdfai^2,\\
\mu P_1^t&=-\dfrac{1}{2}\left(P_{1,\dl}+\eta^{-2}\mu P_{1,L}\right)\\
&=-\dfrac{\Omega}{2}(-\til{T}_{L\dl}-\til{T}_{L L})=\dfrac{\Omega}{2}\left(\eta^{-2}\mu(L\fai)^2+\mu\sdfai^2\right),\\
\mu P_1^u&=-P_{1,L}=\Omega T_{LL}=\Omega(L\fai)^2.
\end{split}
\ee
Since $\Omega\sim1$ and $\eta\sim1$, one may define the energies and fluxes as:
\beeq
E_0(t)&=&\int_{\Si_t^u}(\dl\fai)^2+\mu^2\sdfai^2 du'd\mu_{\sg},\\
E_1(t)&=&\int_{\Si_t^u}\mu(L\fai)^2+\mu\sdfai^2 du'd\mu_{\sg},\\
F_0(u)&=&\int_{C_u^t}\mu\sdfai^2 dt'd\mu_{\sg},\\
F_1(u)&=&\int_{C_u^t}(L\fai)^2 dt'd\mu_{\sg}.
\eeq
Applying Lemma\ref{divergencethm} to \eqref{DPQ}, using \eqref{P0tP0uP1tP1u} and noticing that $\fai$ vanishes on $C_0^u$ yield 
\bee
\int_{\Si_t^u}\mu P_{\a}^t-\int_{\Si_{-2}^u}\mu P_{\a}^t+\int_{C_u^t}\mu P_{\a}^u=\int_{W_t^u}Q_{\a},\ \a=0,1,
\ee
which implies
\bee\label{energy1}
E_{\a}(t)+F_{\a}(t)\les E_{\a}(-2)+\int_{W_t^u}Q_{\a},\ \a=0,1,\\
\ee
where $Q_{\a}=-\left(\Omega^2 f K_{\a}\fai+\dfrac{1}{2}\til{T}^{\mu\nu}\til{\pi}_{\a,\mu\nu}\right)$. Due to \eqref{relationdeformationtensor} and Lemma\ref{connection2}, the deformation tensors associated with $K_0$ and $K_1$ can be computed as 
\begin{align*}
&\til{\pi}_{0,LL}=0,\hs \til{\pi}_{0,\dl\dl}=0,\hs \til{\pi}_{0,L\dl}=-2\mu(\mu^{-1}\Omega\dl\mu+\Omega L(\eta^{-2}\mu)+\dl\Omega),\\
&\til{\pi}_{0,LA}=-2\Omega(\eta_A+\zeta_A),\hs \til{\pi}_{0,\dl A}=-2\Omega\mu X_A(\eta^{-2}\mu),\\
&\hat{\til{\s{\pi}}}_{0,AB}=2\Omega\hat{\dc}_{AB},\hs tr\til{\s{\pi}}_0=2\Omega\til{tr}\dc,\\
&\til{\pi}_{1,LL}=0,\hs \til{\pi}_{1,\dl\dl}=4\Omega\mu L(\eta^{-2}\mu),\hs \til{\pi}_{1,L\dl}=-2\Omega L\mu-2\mu(L\Omega),\\
&\til{\pi}_{1,LA}=0,\hs \til{\pi}_{1,\dl A}=2\Omega(\eta_A+\zeta_A),\\
&\hat{\til{\s{\pi}}}_{1,AB}=2\Omega\hat{\chi}_{AB},\hs tr\til{\s{\pi}}_1=2\Omega\til{tr}\chi.
\end{align*}
Rewrite $Q_0$ as $Q_0=-\Omega^2f K_0\fai-\dfrac{1}{2}T^{\mu\nu}\til{\pi}_{0,\mu\nu}=\sum_{\a=0}^7Q_{0,\a}$, where
\begin{align*}
Q_{0,0}&=-\Omega^2f K_0\fai=\Omega^2\mu^{-1}\left(\dfrac{2\eta'}{\eta^2}a\fai\mu\de\fe-\dfrac{a\p}{\eta}
\dl\fai\right)\dl\fai,\\
Q_{0,1}&=-\dfrac{1}{2}T^{LL}\til{\pi}_{0,LL}=-\dfrac{1}{8}\mu^{-2}(\dl\fai)^2\til{\pi}_{0,LL}=0,\\
Q_{0,2}&=-\dfrac{1}{2}T^{\dl\dl}\til{\pi}_{0,\dl\dl}=-\dfrac{1}{8}\mu^{-2}(L\fai)^2\til{\pi}_{0,\dl\dl}=0,\\
Q_{0,3}&=-T^{L\dl}\til{\pi}_{0,L\dl}=
-\dfrac{1}{4}\mu^{-1}\sdfai^2\til{\pi}_{0,L\dl}=
\boxed{
\dfrac{1}{2}\sdfai^2\left(\mu^{-1}\Omega\dl\mu+\dl\Omega+\Omega L(\eta^{-2}\mu)\right)},\\
Q_{0,4}&=-T^{LA}\til{\pi}_{0,LA}=\dfrac{1}{2}\mu^{-1}(\dl\fai)(\sd^A\fai)\til{\pi}_{0,LA}=-\Omega\mu^{-1}(\dl\fai)(\sd^A\fai)(\zeta_A+\eta_A),\\
Q_{0,5}&=-T^{\dl A}\til{\pi}_{0,\dl A}=\dfrac{1}{2}\mu^{-1}(L\fai)(\sd^A\fai)\til{\pi}_{0,\dl A}=\Omega(L\fai)(\sd^A\fai)X_A(\eta^{-2}\mu).
\end{align*}

 Note that
\beeq
T^{AB}\til{\pi}_{0,AB}&=&\left[\underbrace{((\sd^A\fai)(\sd^B\fai)-\dfrac{1}{2}\sg^{AB}\sdfai^2)}_{\text{teace-free}}+\dfrac{1}{2\mu}\sg^{AB}(L\fai)(\dl\fai)\right]\cdot\left(\hat{\til{\s{\pi}}}_{0,AB}+\dfrac{1}{2}tr\til{\s{\pi}}_0\sg_{AB}\right)\\
&=&\left[(\sd^A\fai)(\sd^B\fai)-\dfrac{1}{2}\sg^{AB}\sdfai^2\right]\cdot\hat{\til{\s{\pi}}}_{0,AB}+\dfrac{1}{2\mu}tr\til{\s{\pi}}_0(L\fai)(\dl\fai).
\eeq
Thus, $-\dfrac{1}{2}T^{AB}\til{\pi}_{0,AB}=Q_{0,6}+Q_{0,7}$, where
\begin{align*}
Q_{0,6}&=-\dfrac{1}{2}\left[(\sd^A\fai)(\sd^B\fai)-\dfrac{1}{2}\sg^{AB}\sdfai^2\right]\cdot\hat{\til{\s{\pi}}}_{0,AB}=
-\Omega\hat{\dc}_{AB}(\sd^A\fai)(\sd^B\fai),\\
Q_{0,7}&=-\dfrac{1}{4\mu}tr\til{\s{\pi}}_0(L\fai)(\dl\fai)=-\dfrac{1}{2\mu}\Omega\til{tr}\dc(L\fai)(\dl\fai).
\end{align*}
Similar computations lead to the decomposition of $Q_1=\sum_{\a=0}^7Q_{1,\a}$, where
\begin{align*}
Q_{1,0}&=-\Omega^2f K_1\fai=\Omega^2\mu^{-1}\left(\dfrac{2\eta'}{\eta^2}a\fai\mu\de\fe-\dfrac{a\p}{\eta}
\dl\fai\right)L\fai,\\
Q_{1,1}&=-\dfrac{1}{2}T^{LL}\til{\pi}_{1,LL}=-\dfrac{1}{8}\mu^{-2}(\dl\fai)^2
\til{\pi}_{1,LL}=0,\\
Q_{1,2}&=-\dfrac{1}{2}T^{\dl\dl}\til{\pi}_{1,\dl\dl}=-\dfrac{1}{8}\mu^{-2}(L\fai)^2\til{\pi}_{1,\dl\dl}=-\dfrac{1}{2\mu}\Omega(L\fai)^2L(\eta^{-2}\mu),\\
Q_{1,3}&=-T^{L\dl}\til{\pi}_{1,L\dl}=
-\dfrac{1}{4}\mu^{-1}\sdfai^2\til{\pi}_{1,L\dl}=\boxed{\dfrac{1}{2\mu}\sdfai^2\left(\mu L\Omega+\Omega L\mu\right)},\\
Q_{1,4}&=-T^{LA}\til{\pi}_{1,LA}=\dfrac{1}{2}\mu^{-1}(\dl\fai)(\sd^A\fai)\til{\pi}_{1,LA}=0,\\
Q_{1,5}&=-T^{\dl A}\til{\pi}_{1,\dl A}=\dfrac{1}{2}\mu^{-1}(L\fai)(\sd^A\fai)\til{\pi}_{1,\dl A}=\Omega\mu^{-1}(L\fai)(\sd^A\fai)(\eta_A+\zeta_A),\\
Q_{1,6}&=-\dfrac{1}{2}\left[(\sd^A\fai)(\sd^B\fai)-\dfrac{1}{2}\sg^{AB}\sdfai^2\right]\cdot\hat{\til{\s{\pi}}}_{1,AB}=
-\Omega\hat{\chi}_{AB}(\sd^A\fai)(\sd^B\fai),\\
Q_{1,7}&=-\dfrac{1}{4\mu}tr\til{\s{\pi}}_1(L\fai)(\dl\fai)=-\dfrac{1}{2\mu}\Omega\til{tr}\chi(L\fai)(\dl\fai).
\end{align*}
Due to the definitions of $E$ and $F$, it turns out that the boxed terms are most difficult to handle since they involve $\mu^{-1}\sdfai^2$ and possibly shock formation ($\mu\to0$).\\
\subsubsection{\textbf{Estimates for the error integrals $Q_i$}}
\hs We start with the most difficult terms. First,
\beeq
\int_{W_t^u}Q_{1,3}&=&\int_{W_t^u}\dfrac{1}{2}L\Omega
\sdfai^2+\Omega\int_{W_t^u}\dfrac{1}{2}(\mu^{-1}L\mu)\sdfai^2\\
&\les&\int_{-2}^t E_1(t') dt'+\int_{W_t^u}\dfrac{1}{2}(\mu^{-1}L\mu)\sdfai^2,
\eeq
due to Lemma\ref{hrhoeta}. \textbf{To deal with any integral involving $\mu^{-1}L\mu$ or $\mu^{-1}\dl\mu$, one can split it into shock part (in $W_{shock}$) and non-shock part. The difficult shock part can be bounded by using two key propositions of $\mu$: Proposition\ref{keymu1} and \ref{keymu2} as follows.}
\[
\int_{W_t^u}\dfrac{1}{2}(\mu^{-1}L\mu)\sdfai^2=
\left(\int_{W_t^u\cap W_{shock}}+\int_{W_t^u\cap
W_{n-s}}\right)\dfrac{1}{2}(\mu^{-1}L\mu)\sdfai^2,
\]
where $W_{n-s}=W_t^u\setminus W_{shock}$. 
The integral in the non-shock region is easy to handle. Since $\mu$ has a positive lower bound on $W_t^u\cap W_{n-s}$, it follows from Proposition\ref{accruatemu1} that
\bee\label{q13nonshock}
\int_{W_t^u\cap W_{n-s}}\mu^{-1}L\mu\sdfai^2\les\int_{-2}^tE_1(t')dt'.
\ee
\hs In the shock region, it follows from Proposition\ref{keymu1} that
\bee\label{q13shock}
\int_{W_t^u\cap W_{shock}}(\mu^{-1}L\mu)\sdfai^2\les\int_{W_t^u\cap W_{shock}}-\mu^{-1}\sdfai^2:=-K(t,u),
\ee
where $K(t,u)$ is a non-negative integral defined by $K(t,u)=\int_{W_t^u\cap W_{shock}}\mu^{-1}\sdfai^2$, which plays a key role to control $\int\mu^{-1}\sdfai^2$ in the shock region. Combining \eqref{q13nonshock} with \eqref{q13shock} yields
\bee\label{estq13}
\int_{W_t^u}Q_{1,3}\les\int_{-2}^t E_1(t') dt'-K(t,u).
\ee
\hs Next, we treat $Q_{0,3}$.
\beeq
\int_{W_t^u}Q_{0,3}&=&\int_{W_t^u}\mu^{-1}T\mu\sdfai^2+\dfrac{1}{2}\sdfai^2
(\dl\Omega+L(\eta^{-2}\mu)+\eta^{-2}L\mu)\\
&\les&\int_{W_t^u}\mu^{-1}T\mu\sdfai^2+M\int_0^uF_0(u')du',
\eeq
due to Lemma\ref{hrhoeta}. As before,
\[
\int_{W_t^u}\mu^{-1}T\mu\sdfai^2=\left(\int_{W_t^u\cap W_{shock}}+\int
_{W_t^u\cap W_{n-s}}\right) \mu^{-1}T\mu\sdfai^2.
\]
The integral over non-shock region can be bounded as
\[
\int
_{W_t^u\cap W_{n-s}}\mu^{-1}T\mu\sdfai^2\les\da^{-1}\int_{-2}^tE_1(t')dt',
\]
due to Lemma\ref{ltfe}; while for the integral over shock region, it follows from Proposition\ref{keymu2} that
\[
\int_{W_t^u\cap W_{shock}}\mu^{-1}T\mu\sdfai^2\les\da^{-1}\int_{-2}^t\dfrac{1}{|e^{-a\tau}-e^{-as^{\ast}}|^{\frac{1}
{2}}}E_1(\tau)d\tau.
\]
Hence, it holds that
\[
\int_{W_t^u}Q_{0,3}\les\da^{-1}\int_{-2}^t\dfrac{1}{|e^{-a\tau}-e^{-as^{\ast}}|^{\frac{1}
{2}}}E_1(\tau)d\tau+\da^{-1}\int_{-2}^tE_1(t')dt'+M\int_0^uF_0(u')du'.
\]
This finishes the estimates for the two boxed terms. Next, $Q_{1,2},Q_{1,6},Q_{1,7}$, $Q_{0,5},Q_{0,6}$ and $Q_{0,7}$ are estimated directly by the definition of $E_{\a},F_{\a}$ ($\a=0,1$) as follows.
\beeq
&&\int_{W_t^u}Q_{1,2}\les M\int_{0}^uF_1(u')du',\hs \int_{W_t^u}Q_{1,6}\les
\da M\int_0^uF_0(u')du',\\
&&\int_{W_t^u}Q_{1,7}\les\int_{-2}^t\da E_0(t')dt'+\da^{-1}\int_0^uF_1(u')du',
\hs \int_{W_t^u}Q_{0,5}\les\int_{-2}^tE_1(t')dt',\\
&&\int_{W_t^u}Q_{0,6}\les\da\int_{-2}^tE_0(t')dt',\hs \int_{W_t^u}Q_{0,7}\les
\da\int_{-2}^tE_0(t')dt'+\da^{-1}\int_0^uF_1(u')du'.
\eeq
\hs Next, we estimate $Q_{1,0}$. Since
\bee\label{mudefe}
\bes
\mu\de\fe&=\mu\da^{ij}\dfrac{\pa\fai_j}{\pa x^i}=\mu(\hat{T}^i\hat{T}^j+\sg^{AB}X_A^iX_B^j)\dfrac{\pa\fai_j}{\pa x^i}\\
&=\eta\hat{T}^iT\fai_j+\mu\sg^{AB}X_A^jX_B\fai_j,
\end{split}
\end{equation}
then
\begin{align*}
\int_{W_t^u}Q_{1,0}&\les\int_{W_t^u}\mu^{-1}\fai T\fai L\fai+\fai\sd\fai L\fai+\mu^{-1}\dl\fai L\fai\\
&\les\da\int_{-2}^tE_0(t')dt'+\da^{-1}\int_{0}^uF_1(u')du'.
\end{align*}
Similarly,
\[
\int_{W_t^u}Q_{0,0}\les\int_{-2}^tE_0(t')dt'+\int_0^uF_1(u')du'.
\]
Finally, we treat $Q_{1,5}$ and $Q_{0,4}$. First, due to \eqref{estzetaetaA}, it holds that
\[
\int_{W_t^u}Q_{1,5}\les\int_{W_t^u}\mu^{-1}(L\fai)\cdot\sd\fai(\eta_A+\zeta_A)\les\left(\int_{W_t^u\cap W_{shcok}}+\int_{W_t^u\cap W_{n-s}}\right)\mu^{-1}(L\fai)\cdot\sd\fai.
\]
As before, one can get
\[
\int_{W_t^u\cap W_{n-s}}\mu^{-1}(L\fai)\cdot\sd\fai\les\int_{W_t^u}(L\fai)^2+
\sdfai^2=\int_{-2}^tE_1(t')dt';
\]
while for the integral over shock region, it holds that
\beeq
\int_{W_t^u\cap W_{shock}}\mu^{-1}(L\fai)\cdot\sd\fai&\les&\left(
\int_{W_t^u\cap W_{shock}}\mu^{-1}(L\fai)^2\right)^{\frac{1}{2}}\cdot\left(\int_{W_t^u\cap
W_{shock}}\mu^{-1}\sdfai^2\right)^{\frac{1}{2}}\\
&\les&\da K(t,u)+\da^{-1}\int_0^uF_1(u')du'.
\eeq
The same approach applying to $Q_{0,4}$ yields
\[
\int_{W_t^u}Q_{0,4}\les\int_{-2}^tE_0(t')dt'+\int_0^uF_0(u')du'
+K(t,u).
\]
It follows from \eqref{energy1} that 
\bee\label{energy2}
\bes
&E_0(t,u)+F_0(t,u)+\da^{-1}(E_1(t,u)+F_1(t,u))\les (E_0+\da^{-1}E_1)(-2,u)\\
&\left|\int_{W_t^u}\sum_{i=1}^7Q_{0,i}\right|+
\left|\da^{-1}\int_{W_t^u}\sum_{i=1}^7Q_{1,i}\right|+\underbrace
{\left|\int_{W_t^u}f\cdot\dl \fai\right|}_{Q_{0,0}}+\underbrace{
\left|\int_{W_t^u}\da^{-1}f\cdot L \fai\right|}_{Q_{1,0}}.
\end{split}
\ee
It follows from the estimates for $Q_{\a}(\a=0,1)$ above that for sufficiently small $\da$, the right hand side of\eqref{energy2} is bounded by
\bee\label{Q0Q1bound2}
\bes
&\int_{-2}^tE_0(t')dt'+\da^{-1}\int_0^uF_0(u')du'+\da^{-1}\int_{-2}^tE_1(t')dt'
+\da^{-1}\int_0^u\da^{-1}F_1(u')du'\\
&-\da^{-1}K(t,u)+\da^{-1}\int_{-2}^t\dfrac{1}{|e^{-a\tau}-e^{-as^{\ast}}|^{\frac{1}
{2}}}E_1(\tau)d\tau.
\end{split}
\ee
Applying Gronwall inequality to \eqref{energy2} and using \eqref{Q0Q1bound2} yield
\bee\label{energy3}
\bes
&E_0(t,u)+F_0(t,u)+\da^{-1}(E_1(t,u)+F_1(t,u))\les (E_0+\da^{-1}E_1)(-2,u)\\
&-\da^{-1}K(t,u)+\da^{-1}\int_{-2}^t\dfrac{1}{|e^{-a\tau}-e^{-as^{\ast}}|^{\frac{1}
{2}}}E_1(\tau)d\tau.
\end{split}
\ee
Since $\dfrac{1}{|e^{-a\tau}-e^{-as^{\ast}}|^{\frac{1}
{2}}}$ is integrable in $\tau$, applying Gronwall inequality to \eqref{energy3} yields
\bee\label{energyfinal}
E_0(t,u)+F_0(t,u)+\da^{-1}(E_1(t,u)+F_1(t,u)+K(t,u))\les (E_0+\da^{-1}E_1)(-2,u).
\ee
\begin{remark}\label{energyremark}
\textbf{
To derive the energy estimates for  high order derivatives of $\fai$, the corresponding error integrals $Q_{0,i},Q_{1,i}$ (for $i\ge 1$) can be estimated as above exactly by high order energies and fluxes directly. However, the treatment for $Q_{0,0},Q_{1,0}$ will involve complicated high order acoustical terms and will be given details later.}
\end{remark}
\subsection{\textbf{Some top order acoustical terms}}
\hs Let $\fai$ be a solution to $\square_{\til{g}}\fai=f$ and $Z$ be any one of $\{Q,T,R_i\}$. Then, direct computations yield (c.f.\cite{christodoulou2014compressible}Prop7.1):
\bee\label{commute1}
\square_{\til{g}}(Z\fai)=Zf+\dfrac{1}{2}tr_{\til{g}}\zgpi\cdot f+\Omega^{-2} div_g\pre Z {} J,
\ee
where $J$ is given by:
\[
\pre Z {} J^{\nu}=(\zgpi^{\mu\nu}-\dfrac{1}{2}tr_g\zgpi\cdot g^{\mu\nu})\pa_{\mu}\fai.\]
Here, $\zgpi^{\mu\nu}=g^{\a \mu}g^{\be \nu}\zgpi_{\a\be}$. For the wave equation \eqref{nonlinearwave}: $\Box_{\til{g}}\fai=f$, let $\fai_1=\fai,\p_1=f, J_1=J$ and $\fai_n=Z\fai_{n-1}$. Then, the following relations hold:
\bee\label{recursion1}
\bes
\square_{\til{g}}\fai_n&=\p_n,\fai_n=Z\fai_{n-1},\p_1=-\dfrac{2\eta'}{\eta^2}a\fai\de
\fe+\dfrac{a}{\p\eta}\left(\dfrac{\pa\fai}{\pa t}-\fai_i\dfrac{\pa\fai}{\pa x^i}\right),\\
\p_n&=Z\p_{n-1}+\dfrac{1}{2}tr_{\til{g}}\zgpi\cdot\p_{n-1}+\Omega^{-2} div_g\pre Z {} J_{n-1},\\
\pre Z {} J^{\nu}_{n-1}&=(\zgpi^{\mu\nu}-\dfrac{1}{2}tr_g\zgpi\cdot g^{\mu\nu})\pa_{\mu}\fai_{n-1}.
\end{split}
\ee
We rescale $\p_n$ as:
\bee\label{recursion1'}
\bes
\til{\p}_n&=\Omega^2\mu\p_n=Z\til{\p}_{n-1}+\left(\dfrac{1}{2}tr_{\til{g}}\zgpi
-\mu^{-1}Z\mu-2Z(\log \Omega)\right)\cdot\til{\p}_{n-1}+\mu div_g\pre Z {} J_{n-1}\\
&=Z\til{\p}_{n-1}+\pre {Z} {} {\da}\cdot\til{\p}_{n-1}+\pre Z {} {\sigma_{n-1}},\\
\til{\p}_1&=\Omega\left(-\dfrac{2\eta'}{\p\eta}a\fai\mu\de\fe+\dfrac{a\eta}{\p}
\dl\fai\right),
\end{split}
\ee
where
\bee
\pre {Z} {} {\da}=\dfrac{1}{2}tr_{\til{g}}\zgpi
-\mu^{-1}Z\mu-2Z(\log \Omega),\quad \pre Z {} {\sigma_{n-1}}=\mu div_g\pre Z {} J_{n-1}.
\ee
\hs 
It follows from the estimates for the deformation tensors \eqref{estqgpi}, \eqref{esttgpi} and \eqref{estrgpi} and the relation $tr_{\til{g}}\zgpi=-\mu^{-1}\Omega^{-1}\zgpi_{L\dl}+tr_{\til{g}}\szgpi$ that $|\pre {Z} {} {\da}|\les\da M$ for $Z=R_i, Q$, and $|\pre {Z} {} {\da}|\les M$ for $Z=T$.\\
\hs To proceed, one needs the following elementary lemma, which can be verified inductively.
\begin{lem}\label{recur}
Let $y_n (n=1,2,...)$ be a sequence 
and $A_n$ be a given sequence of operators. Suppose that $x_n$ is a sequence satisfying: $x_n=A_{n-1}x_{n-1}+y_{n-1}.$ Then,
\[
x_n=A_{n-1}A_{n-2}\cdots A_1 x_{1}+\sum_{m=0}^{n-2}A_{n-1}\cdots A_{n-m}y_{n-m-1}.\]
\end{lem}
Applying this lemma to \eqref{recursion1'} yields
\bee\label{recursion2}
\bes
\til{\p}_n&=(Z_{n-1}+\pre {Z_{n-1}} {} {\da})\cdots(Z_{1}+\pre {Z_{1}} {} {\da})\til{\p}_{1}+\sum_{k=0}^{n-2}(Z_{n-1}+\pre {Z_{n-1}} {} {\da})\cdots(Z_{n-k}+\pre {Z_{n-k}} {} {\da})\pre {Z_{n-k-1}} {} {\sigma_{n-k-1}}.
\end{split}
\ee
We will investigate the terms in $\pre Z {} {\sigma_{n-k-1}}$ by using the following lemma.
\begin{lem}\label{divdecom}
For any space time vector field $J$, 
it holds that
\bee\label{JdecomposeLTX}
J=(-\eta^{-2}J_L-\mu^{-1}J_{T})L-\mu^{-1}J_LT+\s{J}.
\ee
Moreover, div$J$ can be expressed as
\beeq
\mu \text{div} J&=&-L(\eta^{-2}\mu J_{L})-L(J_T)-T(J_{L})+\s{div}(\mu\s{J})-(\eta^{-2}\mu tr\chi+\eta^{-1}\mu tr\ta) J_L-tr\chi J_T\\
&=&-\dfrac{1}{2}L({J_{\dl}})-\dfrac{1}{2}\dl(J_{L})+\s{div}(\mu\s{J})
-\dfrac{1}{2}L(\eta^{-2}\mu) J_{L}-\dfrac{1}{2}tr\chi J_{\dl}-\dfrac{1}{2}tr\dc J_L.
\eeq
\end{lem}
 Note that
\bee
\bes
&\pre Z {} {J_{n-1,L}}=-\dfrac{1}{2}tr_{\sg}\szgpi(L\fai_{n-1})+\zgpi_{L}\cdot\s{d}\fai_{n-1},\\
&\pre Z {} {J_{n-1,\dl}}=-\dfrac{1}{2}tr_{\sg}\szgpi(\dl\fai_{n-1})+\zgpi_{\dl}\cdot\s{d}\fai_{n-1}
-\dfrac{1}{2}\mu^{-1}\zgpi_{\dl\dl}(\dl\fai_{n-1}),\\
&\mu\s{J}_{n-1}^A=-\dfrac{1}{2}\zgpi_L^A(\dl\fai_{n-1})-\dfrac{1}{2}\zgpi_{\dl}^A(L\fai_{n-1})+\dfrac{1}{2}(\zgpi_{L\dl}-\mu tr_{\sg}\szgpi)\s{d}\fai_{n-1}^A+\mu\szgpi_B^A\s{d}\fai_{n-1}^B.
\end{split}
\ee
Applying Lemma\ref{divdecom} to $\pre Z {} {\sigma_{n-1}}$ yields the following decomposition:
\bee\label{decomsi}
\pre Z {} {\sigma_{n-1}}=\pre Z {} {\sigma_{1,n-1}}+\pre Z {} {\sigma_{2,n-1}}+\pre Z {} {\sigma_{3,n-1}},
\ee
where $\pre Z {} {\sigma_{i,n-1}} (i=1,2,3)$ are given separately as follows:
\beeq
\pre Z {} {\sigma_{1,n-1}}&&=\dfrac{1}{2}tr\szgpi\left(L\dl\fai_{n-1}
+\dfrac{1}{2}\til{tr}\chi\dl\fai_{n-1}\right)\\
&&+\dfrac{1}{4}\mu^{-1}\zgpi_{L\dl}L^2\fai_{n-1}-\dfrac{1}{2}\zgpi_{\dl}\cdot
\lie_{L}\sd\fai_{n-1}-\dfrac{1}{2}\zgpi_{L}\lie_{\dl}\sd\fai_{n-1}\\
&&-\dfrac{1}{2}\zgpi_{L}\cdot\sd\dl\fai_{n-1}-\dfrac{1}{2}\zgpi_{\dl}\cdot\sd L\fai_{n-1}+\mu\pre {(Z)} {} {\hat{\til{\s{\pi}}}}\cdot\s{D}^2\fai_{n-1}+\dfrac{1}{2}\zgpi_{L\dl}\s{\de}\fai
_{n-1},
\eeq
\beeq
\pre Z {} {\sigma_{2,n-1}}&&=\dfrac{1}{4}(\dl tr\szgpi)L\fai_{n-1}+\dfrac{1}{4}(L tr\szgpi)\dl\fai_{n-1}+\dfrac{1}{4}L(\mu^{-1}\zgpi_{\dl\dl})\dl\fai_{n-1}\\
&&-\dfrac{1}{2}L(\zgpi_{\dl})\sd\fai_{n-1}-\dfrac{1}{2}\dl(\zgpi_L)\sd\fai_{n-1}+\dfrac{1}{2}\s{d}\zgpi_{\dl L}\s{d}\fai
\\
&&-\dfrac{1}{2}\s{div}(\zgpi_{L})\dl\fai_{n-1}-\dfrac{1}{2}\s{div}(\zgpi_{\dl})
L\fai_{n-1}+\s{div}(\mu\pre {(Z)} {} {\hat{\til{\s{\pi}}}})\cdot\sd\fai_{n-1},
\eeq
\beeq
\pre Z {} {\sigma_{3,n-1}}&&=\pre Z {} {\sigma^L_{3,n-1}}L\fai_{n-1}+
\pre Z {} {\sigma^{\dl}_{3,n-1}}\dl\fai_{n-1}+\pre Z {} {\s{\sigma}_{3,n-1}}\sd\fai_{n-1},\\
\pre Z {} {\sigma^L_{3,n-1}}&&=\dfrac{1}{4}L(\eta^{-2}\mu)tr\szgpi+\dfrac{1}{4}\mu^{-1}
\zgpi_{\dl\dl}+\dfrac{1}{4}tr\dc tr\szgpi,\\
\pre Z {} {\sigma^{\dl}_{3,n-1}}&&=-\dfrac{1}{4}L(\log\Omega)tr\szgpi,\\
\pre Z {} {\s{\sigma}_{3,n-1}}&&=-\dfrac{1}{2}tr\szgpi\cdot\Lambda-\dfrac{1}{2}
(L(\eta^{-2}\mu)+tr\dc)\zgpi_{L}-\dfrac{1}{2}tr\chi\zgpi_{\dl}.
\eeq
\textbf{Note that $\pre Z {} {\sigma_{2,n-1}}$ is the major term which contains the product of $1^{st}$ order derivatives of both deformation tensors of $Z$ and $\fai_{n-1}$. $\pre Z {} {\sigma_{1,n-1}}$ contains the product of deformation tensors of $Z$ and $2^{nd}$ order derivatives of $\fai_{n-1}$ and $\pre Z {} {\sigma_{3,n-1}}$ contains the rest which is easy to handle.}\\
\hs In the paper, the highest order of an object will be $N_{top}+1$ (which is called the top order) and we consider the case of $n=N_{top}+1$ in \eqref{recursion1} with $\fai_n=Z_{n-1}\cdots Z_1\fai$.\\
\hs Then we will consider the top order acoustical terms in $\til{\p}_n$ associated with the different choices of $\fai_n$ and the lower order terms will be estimated in the next subsection. Since the Lie-brackets $[R_i,T]$, $[R_i,Q]$ and $[T,Q]$ are one order lower than their product, then up to the commutation terms, there are following 3 possibilities of $\fai_n$ :
\[
\fai_n=R_i^{\a+1}\fai,\quad \fai_n=R_i^{\a'}T^{l+1}\fai,\quad \fai_n=R_i^{\a'}T^{l'}Q^{\be'+1}\fai,
\]
where $N_{top}=|\a|+1=|\a'|+l+1=|\a'|+l'+\be'+1$. \textbf{Since the top order derivatives of $\fai$ can be bounded directly by the associated energies and fluxes, it suffices to consider the top order acoustical terms. That is, we consider the top order spatial derivatives for $\mu$ and $\chi$.} If $\be'>0$, then by the transport equations for $\mu$ and $\chi$, $L \mu $ and $Ltr\chi$ consist terms of derivatives of $\fai$ and lower order acoustical terms. Thus, it suffices to consider the case $\be'=0.$ \\
\hs  It follows from \eqref{recursion2} that there are at most $k$ derivatives acting on $\pre {Z_{n-k-1}} {} {\sigma_{n-k-1}}$, at most $n-1$ derivatives acting on $\til{\p}_1$ and at most $n-2$ derivatives acting on $\pre {Z_{n-1}} {} {\da}$. Since $\til{\p}_1$ contains $\mu$ and $\pre {Z_{1}} {} {\da}$ contains the first order derivatives of $\mu$, then $(Z_{n-1}+\pre {Z_{n-1}} {} {\da})\cdots(Z_{1}+\pre {Z_{1}} {} {\da})\til{\p}_{1}$ in \eqref{recursion2} will only produce the acoustical terms with order of $\leq N_{top}$, which are not top order. While by the decomposition of $\pre Z {} {\sigma_{n-1}}$, the highest order acoustical terms in $(Z_{n-1}+\pre {Z_{n-1}} {} {\da})\cdots(Z_{n-k}+\pre {Z_{n-k}} {} {\da})\pre {Z_{n-k-1}} {} {\sigma_{n-k-1}}$ are at most of order $k+2$, which are contained in $Z_{n-1}\cdots Z_{n-k}\pre {Z_{n-k-1}} {} {\sigma_{n-k-1}}$. Therefore, it follows that the possible top order acoustical terms will be contained in $Z_{n-1}\cdots Z_2\pre {Z_1} {} {\sigma_{2,1}}$.\\
\hs Since $L\mu$ and $Ltr\chi $ consist lower order acoustical terms, we define the principle acoustical terms (P.A.) of $\pre {Z} {} {\sigma_{2,1}}$ by keeping only spatial derivatives for the deformation tensors as follows.
     \bee\label{paz}
     \bes
     [\pre {Z} {} {\sigma_{2,1}}]_{P.A.}&=\biggl[\dfrac{1}{2}T(tr\szgpi)L\fai-\lie_{T}\zgpi_{L}\cdot\sd\fai-\dfrac{1}{2}
     \s{div}\zgpi_{\dl}L\fai\\
     &-\dfrac{1}{2}\s{div}\zgpi_{L}\dl\fai+\dfrac{1}{2}\sd\zgpi_{L\dl}\cdot\sd\fai+
     \s{div}(\mu\pre {(Z)} {} {\hat{\til{\s{\pi}}}})\cdot\sd\fai\biggl]_{P.A.}.
     \end{split}
     \ee
     Similarly, one can define the principle acoustical terms (P.A.) for the deformation tensors.
     Then we study the top order acoustical terms in $Z_{n-1}\cdots Z_2\pre {Z_1} {} {\sigma_{2,1}}$ and there are 3 possibilities corresponding to the different choices of $\fai_n$.
 \begin{itemize}
 \item[Case1:] 
      For $Z=Q$, it holds that
     \beeq
     &&[tr\sqgpi]_{P.A.}=2\Omega t\cdot\til{tr}\chi,\hs [\qgpi_{L}]_{P.A.}=0,\hs
     [\qgpi_{L\dl}]_{P.A.}=0,\\
     &&[\qgpi_{\dl}]_{P.A.}=[2\Omega t(2\zeta_A+\sd\mu)]_{P.A.}=2\Omega t\sd\mu,\hs [\mu\pre {(Q)} {} {\hat{\til{\s{\pi}}}}]_{P.A.}=2\mu\Omega t\hat{\chi}.
     \eeq
     Substituting these into \eqref{paz} and noticing \eqref{Ttrchi} and \eqref{Codazzi2} yield 
     \bee
     \bes
     [\pre {(Q)} {} {\sigma_{2,1}}]_{P.A.}&=[\Omega tL\fai(T\til{tr}\chi-\s{\de}
     \mu)+2\mu\Omega t\s{div}\hat{\chi}\cdot\sd\fai]_{P.A.}\\
     &=\mu \Omega t\sd tr\chi\cdot\sd\fai.
     \end{split}
     \ee
     Hence, for the variation $R_i^{\a'}T^lQ\fai$, i.e. $Z_1=Q,$ $Z_{n-1}\cdots Z_2=R_i^{\a'}T^{l}$, the corresponding top order acoustical terms are
     \bee\label{paq}
     \left\{\bes
     &\mu\Omega t\sd R_i^{\a}tr\chi\cdot\sd\fai,\hs if\ l=0,\\
     &\mu\Omega t R_i^{\a'+1}T^{l-1}\s{\de}\mu\cdot\sd\fai,\hs if\ l\ge 1.
     \end{split}\right.
     \ee
 \item[Case2:] For $Z=T$, it holds that
       \beeq
       &&[tr\stgpi]_{P.A.}=-2\Omega\eta^{-2}\mu\til{tr}\chi,\hs [\tgpi_L]_{P.A.}=-\Omega\sd\mu,\hs [\tgpi_{\dl}]_{P.A.}=-\Omega\eta^{-2}\mu\sd\mu,\\
       &&[\tgpi_{L\dl}]_{P.A.}=-2\Omega T\mu,\hs [\mu\pre {(T)} {} {\hat{\til{
       \s{\pi}}}}]_{P.A.}=-2\Omega\eta^{-2}\mu^2\hat{\chi}.
       \eeq
       Similarly, for the variation $R_i^{\a'}T^{l+1}\fai$, the corresponding top order acoustical terms are
       \bee\label{pqt}
       \left\{\bes
       &\Omega T\fai\cdot R_i^{\a}\s{\de}\mu-\Omega\eta^{-2}\mu^2\sd\fai\cdot\sd R_i^{\a}tr\chi,\hs if\ l=0,\\
       &\Omega T\fai\cdot R_i^{\a'}T^l\s{\de}\mu-\Omega\eta^{-2}\mu^2\sd\fai\cdot R_i^{\a'+1} T^
       {l-1}\s{\de}\mu,\hs if\ l\ge 1.
       \end{split}\right.
       \ee
  \item[Case3:] For $Z=R_i$, it holds that
       \beeq
       &&[tr\srgpi]_{P.A.}=2\eta^{-1}\Omega\lam_i \til{tr}\chi,\hs [\rgpi_L]_{P.A.}=-\Omega R_i\cdot\chi,\hs [\rgpi_{L\dl}]=-2\Omega R_i\mu,\\
       &&[\rgpi_{\dl}]_{P.A.}=\Omega\eta^{-2}\mu R_i\cdot\chi+2\Omega\lam_i\eta^{-1}
       \sd\mu,\hs [\mu\pre {(R_i)} {} {\hat{\til{\s{\pi}}}}]_{P.A.}=2\eta^{-1}\mu\Omega\lam_i\hat{\chi}.
       \eeq
       Hence, for the variation $R_i^{\a+1}\fai$, the corresponding top order acoustical term is
       \bee\label{par}
       \Omega T\fai\cdot\sd R_i^{\a}tr\chi+\eta^{-1}\mu\Omega\lam_i\sd\fai\cdot
       \sd R_i^{\a}tr\chi.
       \ee
 \end{itemize}
 \hs Note that $|\sd\fai|\les\da M$ and $|T\fai|\les M$. Then, it follows from \eqref{paq}, \eqref{pqt} and \eqref{par} that to estimate the top order acoustical terms, one needs only to estimates
 \[
 T\fai\cdot \s{d}R_i^{a}tr\chi,\hs T\fai\cdot R_i^{\a'}T^l\s{\de}\mu,
 \]
 where $|\a|+1=|\a'|+l+1=N_{top}$.
 \subsection{\textbf{
  Estimates for lower order terms}}
 \hs In this subsection, we will give the estimates for lower order terms (the terms of order $\leq N_{top}$), which will be necessary for the estimates of top order acoustical terms. To this end,  
 it will be shown later that the key role in this subsection is to estimate the difference between the acoustical coordinates and the rectangular coordinates.\\ 
\hs In order to estimate the derivatives of the solution up to $k-th$ order, we first define the following modified energies and fluxes as follows. Let $0\leq k\leq N_{top}$ and $(b_k)$ be a sequence of nonnegative integers to be chosen later. Recall $\mu_m^u$ defined by \eqref{defmumu} and define:
\bee
\bes
\wi{E}_{0,k+1}(t, u )&:=\sup_{\tau\in [-2,t]}\{\mu^{ u }_m(\tau)^{2b_{k+1}}E_{0,k+1}(\tau, u )\}:=\sup_{\tau\in [-2,t]}\{\mu^{ u }_m(\tau)^{2b_{k+1}}\sum_{\fai}\sum_{|\a|=k-1}\da^{2l}E_0(Z^{\a+1}\fai)\},\\
\wi{F}_{0,k+1}(t, u )&:=\sup_{\tau\in [-2,t]}\{\mu^{ u }_m(\tau)^{2b_{k+1}}F_{0,k+1}(\tau, u )\}:=\sup_{\tau\in [-2,t]}\{\mu^{ u }_m(\tau)^{2b_{k+1}}\sum_{\fai}\sum_{|\a|=k-1}\da^{2l}F_0(Z^{\a+1}\fai)\},\\
\wi{E}_{1,k+1}(t, u )&:=\sup_{\tau\in [-2,t]}\{\mu^{ u }_m(\tau)^{2b_{k+1}}E_{1,k+1}(\tau, u )\}:=\sup_{\tau\in [-2,t]}\{\mu^{ u }_m(\tau)^{2b_{k+1}}\sum_{\fai}\sum_{|\a|=k-1}\da^{2l}E_1(Z^{\a+1}\fai)\},\\
\wi{F}_{1,k+1}(t, u )&:=\sup_{\tau\in [-2,t]}\{\mu^{ u }_m(\tau)^{2b_{k+1}}F_{1,k+1}(\tau, u )\}:=\sup_{\tau\in [-2,t]}\{\mu^{ u }_m(\tau)^{2b_{k+1}}\sum_{\fai}\sum_{|\a|=k-1}\da^{2l}F_1(Z^{\a+1}\fai)\},
\end{split}
\ee
where $l$ is the number of $T's$ in $Z^{\a+1}\in\{Q,T,R_i\}$ and for $k=0$, we set $E_{\a,1}=E_{\a}$ and $ F_{\a,1}=F_{\a}$ ($\a=0,1$). Denote $E_{0,\leq k+1}:=\sum_{l\leq k}E_{0,l+1}$ and $\wi{E}_{0,\leq k+1}(t,u):=\sup_{\tau\in [-2,t]}\{\mu^{ u }_m(\tau)^{2b_{k+1}}E_{0,\leq k+1}(\tau, u )\}$. Similarly, one can define $E_{1,\leq k+1},\wi{E}_{1,\leq k+1},F_{0,\leq k+1},\wi{F}_{0,\leq k+1},F_{1,\leq k+1}$ and $\wi{F}_{1,\leq k+1}$. Here, we multiply the power of $\mu_m$ in order to control the blow up rate of the energies and fluxes as a shock forms.
\begin{lem}
Let $A(t,u)$ be the area of $S_{t,u}$. Then there is a universal constant $C$ such that
\[
C^{-1}\leq A(t,u)\leq C.
\]
\end{lem}
\begin{pf}
It follows from the definition of $\chi$ and Lemma\ref{2ndff} that
\[
|tr\chi|=|\frac{1}{2}\lie_L\sg|=|\dfrac{1}{\sqrt{\det\sg}}(L\sqrt{\det\sg})|\les1.
\]
Denote $\mu_{\sg}(-2,0)$ be the area form of the sphere $S_{-2,0}$.
Then,
\bee\label{dmusgdmu-2}
d\mu_{\sg}=\dfrac{\sqrt{\det\sg(t,u)}}{\sqrt{\det\sg(-2,0)}}d\mu_{\sg}(-2,0)=
\dfrac{\sqrt{\det\sg(t,u)}}{\sqrt{\det\sg(-2,u)}}\dfrac{\sqrt{\det\sg(-2,u)}}{\sqrt{\det\sg(-2,0)}}d\mu_{\sg}(-2,0).
\ee
Note that
\bee
|\ln\dfrac{\sqrt{\det\sg(t,u)}}{\sqrt{\det\sg(-2,u)}}|=|\int_{-2}^t\dfrac{1}{\sqrt{\det\sg}}(L\sqrt{\det\sg})dt'|\les 1,
\ee
which implies
\bee
|\dfrac{\sqrt{\det\sg(t,u)}}{\sqrt{\det\sg(-2,u)}}|\les 1.
\ee
Then, integrating \eqref{dmusgdmu-2} over $S^2$ yields
\bee
\int_{S^2}d\mu_{\sg(-2,0)}\les \int_{S^2}d\mu_{\sg(t,u)}\les \int_{S^2}d\mu_{\sg(-2,0)}.
\ee
That is, $C^{-1}\leq A(t,u)\leq C$.
\end{pf}
\begin{lem}\label{ellipticenergy}
For any smooth function $\fai$ which vanishes on $C_0$, it holds that
\bee
\int_{S_{t,u}}\fai^2\les\tilde{\da}\int_{\Si_t^{u}}\mu(L\fai)^2+(\dl\fai)^2.
\ee
\end{lem}
\begin{pf}
On the hypersurface $\Si_t^u$, one can choose a new coordinate system on $S_{t,u}$ still denoted to be $(\ta_1,\ta_2)$ such that $T=\dfrac{\pa}{\pa u}|_{t,u,\ta}$. Then,
\beeq
\int_{S_{t,u}}\fai^2&\les&\int_{S_{t,u}}\left(\int_0^{u}T\fai(t,u',\ta)du'\right)^2d\mu_{\sg(t,u)}\\
&\les&\tilde{\da}\int_{S_{t,u}}\int_0^{u}(T\fai(t,u',\ta))^2du'd\mu_{\sg(t,u)}\\
&\les&\tilde{\da}\int_0^{u}\int_{S_{t,u'}}(T\fai(t,u',\ta))^2d\mu_{\sg(t,u')}du',
\eeq
where in last inequality, the boundness of area of $S_{t,u'}$ is used. It follows from the definition of $T$ that $(T\fai)^2\les \mu^2(L\fai)^2+(\dl\fai)^2$, which completes the proof.
\end{pf}
As a corollary, it holds that for all $k\leq N_{top}-1$,
\[
\sum_{|\be|\leq k}\int_{S_{t,u}}(R_i^{\be}\fai)^2\les\tilde{\da} (E_{0,\leq k+1}(t,u)+E_{1,\leq k+1}(t,u)).
\]
\hs One has the following relation between two rotational vector fields, which can be verified directly.
 \begin{lem}
 For any two rotational vector fields $R_i$ and $R_j$, it holds that
 \bee\label{rirj}
 \bes
 [R_i,R_j]&=-\ep_{lij}R_l+\lam_j\ep_{ilm}y^l\s{d}_mx\cdot\sg^{-1}-\lam_i
 \ep_{jlm}y^l\s{d}_mx\cdot\sg^{-1}\\
 &+\lam_i\left(\s{k}-\eta^{-1}\chi'+\sg\dfrac{\eta^{-1}-1}{u-t}\right)\cdot R_j
 -\lam_j\left(\s{k}-\eta^{-1}\chi'+\sg\dfrac{\eta^{-1}-1}{u-t}\right)\cdot R_i.
 \end{split}
 \ee
 \end{lem}
\hs 
Given any $k$ and $j$, define $x^j_{i_1\cdots\i_k}$ and $\pre {(k)} {} {\da}^j_{i_i\cdots i_k}$ as
\[
 R_{i_k}\cdots R_{i_1}x^j=\mari{R}_{i_k}\cdots\mari{R}_{i_1}x^j-\pre {(k)} {} {\da}^j_{i_i\cdots i_k}=x^j_{i_1\cdots\i_k}-\pre {(k)} {} {\da}^j_{i_i\cdots i_k},
\]
where for $k=0$, we set $\pre {(0)} {} {\da}^j=0$.
Obviously, $x^j_{i_1\cdots\i_k}$ are linear functions of the rectangular coordinates and in particular, $\dfrac{\pa(x^j_{i_1\cdots\i_k})}{\pa x^l}:=\pre {(k)} {} {c}^j_{l,i_i\cdots i_k}$ are constants. $\pre {(k)} {} {\da}^j_{i_i\cdots i_k}$ can be computed as follows. Since $\pre {(0)} {} {\da}^j=0$ and $\pre {(1)} {} {\da}^j_i=\lam_i\hat{T}^j,$ it follows that
\bee
\bes
x^j_{i_1\cdots\i_k}-\pre {(k)} {} {\da}^j_{i_i\cdots i_k}=R_{i_k}\cdots R_{i_1}x^j&=R_{i_k}x^j_{i_1\cdots i_{k-1}}-R_{i_k}\pre {(k-1)} {} {\da}^j_{i_i\cdots i_{k-1}},\\
R_{i_k}x^j_{i_1\cdots i_{k-1}}&=x^j_{i_1\cdots i_k}-\pre {(k-1)} {} {c}^j_{l,i_i\cdots i_{k-1}}\lam_{i_k}\hat{T}^{l},
\end{split}
\ee
which implies
\bee\label{recursionda}
\pre {(k)} {} {\da}^j_{i_i\cdots i_k}=R_{i_k}\pre {(k-1)} {} {\da}^j_{i_i\cdots i_{k-1}}+\pre {(k-1)} {} {c}^j_{l,i_i\cdots i_{k-1}}\lam_{i_k}\hat{T}^{l}.
\ee
Applying Lemma\ref{recur} to \eqref{recursionda} yields
\bee\label{recursionda2}
\pre {(k)} {} {\da}^j_{i_i\cdots i_k}=\sum^{k}_{m=1}\pre {(m-1)} {} {c}^j_{l,i_1\cdots i_{m-1}}R_{i_k}\cdots R_{i_{k-m+2}}(\lam_{i_{m}}\hat{T}^l).
\ee
It follows from the definition that $x^j_{i_1\cdots\i_k}$ are linear functions of $x^j$ with uniform coefficients, so can be bounded by $r$, then by some universal constant $C$.
\begin{lem}\label{RikRi1yjLinfty}
For $k=0,1,\cdots N_{\infty}-1$, if $||R_{i_k}\cdots R_{i_i}y^j||_{\supnormda}\les\da$, then the following estimates hold for sufficiently small $\da$
\beeq
 &&||R_{i_k}\cdots R_{i_i}\lam_i||_{\supnormda}\les\da,\\
&&||\pre {(k+1)} {} {\da}^j_{ii_i\cdots i_k}||_{\supnormda}\les\da.
\eeq
\end{lem}
\begin{remark}
Since 
$|x^j_{i_1\cdots i_k}|\les 1$, one can bound $R_{i_k}\cdots R_{i_i}x^j$ in terms of $R_{i_k}\cdots R_{i_i}y^j$ provided that the assumption of Lemma\ref{RikRi1yjLinfty} holds.
\end{remark}
\begin{pf}
The case for $k=0$ holds trivially due to the definition, \eqref{estlami} and \eqref{estyi}.\\ 
\hs In the case $k\ge 1$, since $\lam_i=\bar{g}(\mari{R},\hat{T})=\ep_{imk}x^my^k$ and $\hat{T}^l=y^l
+\dfrac{x^l}{u-t}$, one can prove the desired estimates by using an induction argument and \eqref{recursionda2}.
\end{pf}
Now we turn to the estimates for $R_{i_k}\cdots R_{i_1}y^j$. It follows from Lemma\ref{2ndff} that
\bee
R_iy^j=\left(
\ta_{AB}-\dfrac{\da_{AB}}{u-t}\right)R_i^A\cdot\s{d}_Bx^j,\quad |R_iy^j|\les |\ta'|\les\da M,
\ee
and then, 
\bee
R_{i_k}\cdots R_{i_1}R_iy^j=\lie_{R_{i_k}}\cdots\lie_{R_{i_1}}
\left[\left(
\s{k}_{AB}-\eta^{-1}\chi'_{AB}
+\sg_{AB}\dfrac{\eta^{-1}-1}{u-t}\right)R_i^A\cdot\s{d}_Bx^j\right].
\ee
\hs Therefore, to estimate $||R_i^{\a}y^j||_{\supnormda} $ for $|\a|=k+1\leq N_{\infty}-1$, it suffices to estimate the following terms with $|\be|\leq|\a|-1$:
\been[(1)]
\item $|\lie_{R_j}^{\be}\s{k}_{AB}|=|R_j^{\be}X_A^kX_B^l\pa_k\fai_l|\les
      |R_j^{\be+1}\fai|+\text{l.o.ts}\les\da$, where the first inequality can be shown by an induction argument,
\item $\lie_{R_i}^{\be}\s{d}x^j\sim\lie_{R^i}^{\be}R_j$,
\item $\lie_{R_i}^{\be}\sg_{AB}
    =\lie_{R_i}^{\be-1}\left(-2\lam_j\eta^{-1}
    (\eta\s{k}_{AB}-\chi'_{AB}+\dfrac{\sg_{AB}}{u-t})\right)$, $\lie_{R_i}^{\be}\chi'_{AB}$, $\lie_{R_i}^{\be}R_j$.
\een
For the term $\lie_{R_i}^{\be}R_j$, it follows from \eqref{rirj} that
\bee
\lie_{R_i}^{\be}R_j=(\text{constants})\cdot R_k+\lie_{R_i}^{\leq\be-1}\{\lam_i,\chi',\s{d}x^l,y^i\}R_k.
\ee
Therefore, it follows from Lemma\ref{RikRi1yjLinfty} and an induction argument that
\bee\label{RiayjlieRibechi'}
R_i^{\a}y^j=O(1)\cdot \lie_{R_i}^{\be}\chi'+O_2^{\leq |\a|}.
\ee
Then, it suffices to estimate $\lie_{R_i}^{\be}\chi'$ with $|\be|\leq N_{\infty}-2$.
\begin{lem}\label{Riachi'Linfty}
For sufficiently small $\da$ and all $|\a|\leq  N_{\infty}-2=[\dfrac{N_{top}}{2}]+1$, it holds that
\[
||\lie_{R_i}^{\a}\chi'||_{\supnormda}\les_{M}\da.
\]
\end{lem}
\begin{pf}
The case for $|\a|=0$ follows from Lemma\ref{2ndff}. The general case $|\a|\geq1$ will be treated by an induction argument as follows. 
Commuting \eqref{Lchi'AB} with $\lie_{R_i}^{\a}$ yields
\bee\label{zax}
\begin{split}
\lie_{L}\lie_{R_i}^{\a}\chi'&=[\lie_{L},\lie_{R_i}^{\a}]\chi'+(e+a\eta^{-1}\hat{T}^i\fai_i)\cdot\lie_{R_i}^{\a}\chi'+2\chi'\cdot\lie_{R_i}^{\a}\chi'\\
&+\sum_{|\be_1|+|\be_2|=|\a|,|\be_1|>0}R_i^{\be_1}e\cdot\lie_{R_i}^{\be_2}\chi'+\sum_{|\be_1|+|\be_2|+|\be_3|=|\a|,|\be_2|>0}
\lie_{R_i}^{\be_1}\chi'\cdot\lie_{R_i}^{\be_2}\sg\cdot\lie_{R_i}^{\be_3}\chi'\\
&+\sum_{|\be_1|+|\be_2|=|\a|,|\be_1|>0}R_i^{\be_1}(a\eta^{-1}\hat{T}^i\fai_i)\cdot\lie_{R_i}^{\be_2}\chi'+\lie_{R_i}^{\a}\left(-\dfrac{(e+a\eta^{-1}\hat{T}\fe)\sg_{AB}}{ u -t}-\a'_{AB}\right).
\end{split}
\ee
The right hand side of \eqref{zax} will be estimated as follows.
\been[(1)]
\item For the terms $\lie_{R_i}^{\be}\chi'$ where $|\be|<|\a|$, they can be dealt by an induction process.
\item For $\lie_{R_i}^{\be}\sg_{AB}$ with $|\be|\leq|\a|$, one can bound it by an induction process since
      \begin{equation}
      \lie_{R_i}^{\be}\sg_{AB}
    =\lie_{R_i}^{\be-1}\left(-2\lam_j\eta^{-1}
    (\eta\s{k}_{AB}-\chi'_{AB}+\dfrac{\sg_{AB}}{u-t})\right).
      \end{equation}
\item For $R_i^{\be}e$ and $R_i^{\be}(a\eta^{-1}\hat{T}^j\fai_j)$ with $|\be|\leq|\a|$, it follows from an induction process and the bootstrap assumptions that
      \begin{align}
      &|R_i^{\be}e|=|R_i^{\be}\left(\dfrac{1}{2\eta^2}(\dfrac{\p'}{\p})'Lh+\eta^{-1}
      \hat{T}^i(L\fai_i) \right)|\les|R_i^{\be}Lh|+|R_i^{\be}L\fai_i|+\text{l.o.ts}\les\da M,\\
      &|R_i^{\be}(a\eta^{-1}\hat{T}^j\fai_j)|\les |R_i^{\be}\fai+R_i^{\be}(\hat{T}^j)\fai_j|+\text{l.o.ts}\les\da M.
      \end{align}
      Similarly,
       \bee
       \bes
      |\lie_{R_i}^{\a}\a'_{AB}|&=|\lie_{R_i}^{\a}\left(-\dfrac{1}{2}\dfrac{dH}{dh}
      \s{D}^2_{AB}h+\dfrac{1}{2}\dfrac{dH}{dh}\eta^{-1}\s{k}_{AB}Lh+\a_{AB}^{[N]}
      \right)|\les|R_i^{\a+2}\fai|+\text{l.o.ts}\les\da M.
      \end{split}
      \ee
\item For the commutator, it holds that
    \bee\label{lieLlieRia}
    [\lie_{L},\lie_{R_i}^{\a}]\chi'=\sum_{|\be_1|+|\be_2|=|\a|-1}
    \lie_{R_i}^{\be_1}\cdot\lie_{\pre {R_i} {} {\s{\pi}_{L}}}\cdot\lie_{R_i}^{\be_2}\chi',
    \ee
    where $\pre {R_i} {} {\s{\pi}}_{L A}=\lam_i\eta\mu^{-1}\zeta_A+z^{\lam}\ep_{i\lam l}X_A^l-R_i^{B}\chi'_{AB}$.
     Thus, \eqref{lieLlieRia} can be bounded by
     \bee
     O(\da)\cdot\lie_{R_i}^{\a}\chi'+O_2^{\leq|\a|+1},
     \ee
     due to an induction process and \eqref{estrgpi}.
\een
Therefore, it follows that
\bee\label{lieLlieRiachi'}
\lie_{L}\lie_{R_i}^{\a}\chi'=(e+a\eta^{-1}\hat{T}^i\fai_i+O(\da))\cdot\lie_{R_i}^{\a}\chi'+2\chi'\cdot\lie_{R_i}^{\a}\chi'+O^{\leq|\a|+2}_2,\\
\ee
which implies
\bee
||\lie_{R_i}^{\a}\chi'||_{\supnormda}\les_M\da,
\ee
by integrating \eqref{lieLlieRiachi'} along integral curves of $L$ due to the estimates $|e|+|\fai_i|+|\chi'|\les\da M$.
\end{pf}
\begin{remark}
The similar argument and estimates hold for $R_i$ replaced by $Z_i\in\{Q,T,R_i\}$.
\end{remark}
Thus, we conclude that:
\begin{prop}\label{yxLinfty}
For sufficiently small $\da$ and all $|\a|\leq N_{\infty}-2$, it holds that
\[Z^{\a+1}_iy^j,\da Z_i^{\a+1}x^j,\lie_{Z_i}^{\a}\chi'\in O^{|\a|+1}_{2-2l},
\]
where $l$ is the number of $T$'s in $Z_i^{\a}$.
\end{prop}
We now turn the $L^{\infty}$ estimate for $\mu.$
\begin{prop}\label{ZiamuLinfty}
For sufficiently small $\da$ and all $|\a|\leq N_{\infty}-2$, 
it holds that
\[||Z_i^{\a+1}\mu||_{\supnormda}\les_M\da^{-l},
\]
where $l$ is the number of $T$'s in $Z_i^{\a+1}$.
\end{prop}
\begin{pf}
The case, $|\a|=0$, has been given in Lemma\ref{estmu},\ref{ltfe}. Commuting $Z^{\a+1}_i$ with \eqref{transportmu} yields
\[ L (\da^lZ^{\a+1}_i\mu)=\da^l[ L ,Z^{\a+1}_i]\mu+\da^lZ^{\a+1}_im+e\da^l Z^{\a+1}_i\mu+\sum_{|\be_1|+|\be_2|=|\a|+1,|\be_1|>0}\da^{l_1}Z^{\be_1}e\cdot\da^{l_2}Z^{\be_2}\mu,
\]
where $l_i$ is the number of $T$'s in $Z^{\be_i}$ for $i=1,2$. Since $\da^l[ L ,Z^{\a+1}_i]\mu=\da^{l}\sum_{|\be_1|+|\be_2|=|\a|}Z_i^{\be_1}\cdot\zpi_{ L }\cdot Z_i^{\be_2}\mu$, it can be bounded as
\bee
O(\da)\cdot Z^{\a+1}\mu+O^{\leq|\a|+1}_0,
\ee
due to the same estimates in \eqref{lieLlieRia} and an induction process. Similarly to the argument in the proof of Lemma\ref{Riachi'Linfty}, one can estimate remaining terms except $\da^lZ_i^{\a+1}m$, which can also be bounded by $O^{\leq|\a|+2}_0$ (see the remark below). Therefore, it follows that
\bee\label{LdalZa+1}
 L (\da^lZ^{\a+1}_i\mu)=(e+O(\da))\da^l Z^{\a+1}_i\mu+O^{\leq|\a|+2}_0.
 \ee
Integrating \eqref{LdalZa+1} along integral curves of $L$ yields $||Z^{\a+1}_i\mu||_{\supnormda}\les_M\da^{-l}.$
\end{pf}
\begin{remark}\label{estTmfe}
We explain the detials for estimating $Z_i^{\a+1}m$ here. For the case when $Z_i^{\a+1}=R_i^{\a+1}$ or $Z_i^{\a+1}=R_i^{\a}Q$, then the proof will be same as in Lemma\ref{Riachi'Linfty}. However, when $Z_i^{\a+1}=R_i^{\a'}T_i^{m+1}$ with $|\a'|+m=|\a|$, things will be more complicated since the estimates for $T^m\fe$ for $m>1$ are required. In fact, \textbf{it can be shown that by a bootstrap argument and an induction argument that $||R_i^{\a'}T^{m+1}\fe||_{\supnormda}\les\da^{-m}M$ for all $|\a'|+m\leq N_{\infty}-1$ and as a corollary, $||R_i^{\a'}T^{m}\mu||_{\supnormda}\les\da^{-m}M$}, which will be sketched as follows (see also Lemma\ref{ltfe}).\\
\hs Let $|\a'|+m=N_{\infty}-1$ and we will proceed with induction on $m$. For $m=0$, it follows from Lemma\ref{LTidTi}, the bootstrap assumptions and the argument above that $|R_i^{\a'}T\fe|\les\da M$ and $|R_i^{\a'}\mu|\les M$.\\
\hs Assume that $|R_i^{\be}T^{k+1}\fe|\les\da^{-k}M$ and $|R_i^{\be}T^{k}\mu|\les\da^{-k}M$ for all $|\be|+k=N_{\infty}-1$ and $0\leq k\leq m-1$. Suppose that $|R_i^{\a'}T^{m+1}\fe|\les\da^{-m}M$. Then, it follows that
\bee
\bes
R_i^{\a'}T^{m+1}\fe&=R_i^{\a'}T^{m}(\eta^{-1}\mu\hat{T}^i\fai_i)=\eta^{-1}R_i^{\a'}T^{m}(\mu)\hat{T}^i\fai_i+\eta^{-1}\mu\fai_iR_i^{\a'}T^{m}(\hat{T}^i)+O^{\leq N_{\infty}-1}_{2-2m}.\\
\end{split}
\ee
\been
\item For $R_i^{\a'}T^m\mu$, commuting $R_i^{\a'}T^m$ with \eqref{transportmu} yields
\bee\label{tkmu}
\bes
LR_i^{\a'}T^m\mu&=R_i^{\a'}T^mm+eR_i^{\a'}T^m\mu+\sum_{|\a_1|+|\a_2|=|\a'|,\be_1+\be_2=m,|\a_2|+\be_1>0}R_i^{\a_1}T^{\be_1}e\cdot R_i^{\a_2}T^{\be_2}\mu\\
&+\sum_{\be_1+\be_2=|\a'|-1}R_i^{\be_1}\srpi_LR_i^{\be_2}T^m\mu+\sum_{\be_1+\be_2=m-1}R_i^{\a'}T^{\be_1}\stpi_{L}T^{\be_2}\mu.
\end{split}
\ee
It holds that $|R_i^{\a'}T^mm|\les|R_i^{\a'}T^{m+1}\fai|+|R_i^{\a'}T^{m+1}\fe|\les\da^{-m}M$. By induction, the right hand side of \eqref{tkmu} can be bounded by
$(e+O(\da))R_i^{\a'}T^m\mu+\da^{-m}M$. Then, integrating along integral curves of $L$ yields $|R_i^{\a'}T^m\mu|\les\da^{-m}M$.
\item For $R_i^{\a'}T^m(\hat{T}^i)$, it follows from Lemma\ref{LTidTi} that
\bee
R_i^{\a'}T^m(\hat{T}^i)=R_i^{\a'}T^{m-1}\left(-X_B(\eta^{-1}\mu)\cdot\s{d}_Bx^i\right).
\ee
By the same argument in deriving the estimates for $|R_{i_k}\cdots R_{i_1}y^j|$, it can be shown that $|R_i^{\a'}T^{m}(\hat{T}^i)|=|O(1)R_i^{\a'+1}T^{m-1}\mu|+\text{l.o.ts}\les\da^{-m+1}M$.
\een
Hence, $|R_i^{\a'}T^{m+1}\fe|\les\da^{-m+1}M^2<\da^{-m}M$ for sufficiently small $\da$, which recovers the assumption $|R_i^{\a'}T^{m+1}\fe|\les\da^{-m}M$.
\end{remark}
\begin{remark}
It follows from \eqref{computationqpi}, \eqref{computationtpi}, \eqref{computationrpi1} and \eqref{rpiLATAAB} that for all $|\a|\leq N_{\infty}-2$
\[\lie_{Z_i}^{\a}\pre {Z_j} {} {\s{\pi}_{L}},\lie_{Z_i}^{\a}\pre {Z_j} {} {\s{\pi}}, \da\lie_{Z_i}^{\a}\pre {Z_j} {} {\s{\pi}_{\dl}}\in O^{|\a|+1}_{2-2l},
\]
where $l$ is the number of $T$'s in $Z_i^{\a}$.
\end{remark}
\hs Next, we turn to the $L^2$ estimates of above entities, which follows from the same framework in deriving $L^{\infty}$ estimates.
\begin{lem}
For $k=0,1,2...[\frac{l}{2}]:=[\frac{N_{top}}{2}]$, if $\max_{j,i_1...i_k}||R_{i_k}\cdots R_{i_1}y^j||_{\supnormu}\les\da,$ then for sufficiently small $\da$, it holds that
\begin{align*}
&\max_{j,i_1...i_k}||R_{i_k}\cdots R_{i_1}\lam_j||_{\normu2}\les\sum_{k=0}^l\max_{j,i_1...i_k}||R_{i_k}\cdots R_{i_1}y^j||_{\normu2},\\
&\max_{j,i_1...i_k}||\pre {(k+1)} {} {\da}^j_{ii_i\cdots i_k}||_{\normu2}\les\sum_{k=0}^l\max_{j,i_1...i_k}||R_{i_k}\cdots R_{i_1}y^j||_{\normu2}.
\end{align*}
\end{lem}
\begin{pf}
The lemma holds trivially for $l=0$. The general case $l\geq 1$ can be proved by using an induction argument.
\end{pf}
Next, we estimate $||\lie_{R_i}^{\a}y^j||_{\normu2}$ with $|\a|\leq N_{top}$. It suffices to estimate $||\lie_{R_i}^{\be}\chi'||_{\normu2}$ with $|\be|\leq N_{top}-1$ due to \eqref{RiayjlieRibechi'} and an induction argument. In fact, one has:
\begin{prop}\label{Riachi'L2}
For sufficiently small $\da$ and all $|\a|\leq N_{top}-1$, 
 it holds that
\[||\lie_{R_i}^{\a}\chi'||_{\normu2}\les
\sum_{|\a'|\leq|\a|}||\lie_{R_i}^{\a'}\chi'||_{L^{2}(\Si_{-2}^{\tilde{\da}})}+\int_{-2}^t\mu_m^{-\frac{1}{2}}(t')
\sqrt{E_{1,\leq |\a|+2}(t', u )}dt'.
\]
\end{prop}
\begin{pf}
Since for any $S_{t,u}$ $(0,2)$ tensor $\vartheta$, it holds that
\bee
\bes
L|\vartheta|^2&=L(\sg^{AC}\sg^{BD}\vta_{AB}\vta_{CD})=2(\vta,\lie_L\vta)-4\vta\cdot\chi\cdot\vta\\
&=2(\vta,\lie_L\vta)-4\vta\cdot\chi'\cdot\vta+4\dfrac{|\vta|^2}{u-t},
\end{split}
\ee
which implies that
\bee\label{Lu-t2vta}
\bes
L((u-t)^2|\vta|)&=(u-t)^2L|\vta|-2(u-t)|\vta|\\
&\les (u-t)^2\left(|\lie_L\vta|+|\chi'|\cdot|\vta|+2\dfrac{|\vta|}{u-t}-2\dfrac
{|\vta|}{u-t}\right)\\
&=(u-t)^2(|\lie_L\vta|+|\chi'|\cdot|\vta|).
\end{split}
\ee
Taking $\vta=\lie_{R_i}^{\a}\chi'$ and integrating \eqref{Lu-t2vta} along integral curves of $L$ yield
\[||\lie_{R_i}^{\a}\chi'||_{\normu2}\les||\lie_{R_i}^{\a}\chi'||_{L^{2}(\Si_{-2}^{ u })}+
\int_{-2}^t|||\chi'|\cdot|\lie_{R_i}^{\a}\chi'|+|\lie_{ L }\lie_{R_i}^{\a}\chi'|||_{L^{2}(\Si_{\tau}^{ u })}d\tau.
\]
By \eqref{zax}, Lemma\ref{Riachi'Linfty}, Proposition\ref{yxLinfty} and \ref{ZiamuLinfty}, one can estimate $||\lie_{ L }\lie_{R_i}^{\a}\chi'||_{L^2(\Si_t^u)}$ through $L^{\infty}-L^2$ argument outlined as follows.\\
\hs The basic idea of estimating terms on the right hang side of \eqref{zax} is based on the following.
\been[i)]
\item For the linear terms, one estimates them in $L^2$ norm. For the nonlinear terms factorized at least two terms, one estimates the highest order term in $L^2$ norm. Since $N_{\infty}-2=[\dfrac{N_{top}}{2}]+1$, there is at most one term of order $>N_{\infty}-2$ in each factor.
\item For the terms of order $\leq N_{\infty}-2$ in each factor, one can estimate them in $L^{\infty}$ norm.
\een
Therefore, it follows that
\bee\label{lieRiachi'before}
||\lie_{R_i}^{\a}\chi'||_{\normu2}\les\sum_{|\a'|\leq|\a|}||\lie_{R_i}^{\a'}\chi'||_
 {L^2(\Si_{-2}^{\tilde{\da}})}+\int_{-2}^t\da||\lie_{R_i}^{\a}\chi'||_{L^2(\Si_{t'}^u)}+
 \mu_m^{-\frac{1}{2}}(t')
\sqrt{E_{1,\leq |\a|+2}(t', u )}dt'.
\ee
Then applying Gronwall inequality to \eqref{lieRiachi'before} yields
 \[||\lie_{R_i}^{\a}\chi'||_{\normu2}\les\sum_{|\a'|\leq|\a|}||\lie_{R_i}^{\a'}\chi'||_
 {L^2(\Si_{-2}^{\tilde{\da}})}+\int_{-2}^t\mu_m^{-\frac{1}{2}}(t')
\sqrt{E_{1,\leq |\a|+2}(t', u )}dt'.\]
\end{pf}
Using a similar argument as above, one can also obtain the $L^2$ estimates for $\mu$ as follows.
\begin{prop}\label{Zia+1muL2}
For sufficiently small $\da$ and all $|\a|\leq N_{top}-1$, 
it holds that
\[\da^l||Z^{\a+1}_i\mu||_{\normu2}\les\sum_{|\a'|\leq|\a|}\da^{l'}||Z^{\a'+1}_i\mu||_{L^2(\Si_{-2}^{\tilde{\da}})}+
\int_{-2}^t\sqrt{E_{0,\leq|\a|+2}(\tau, u )}+\mu_m^{-\frac{1}{2}}(\tau)\sqrt{E_{1,\leq|\a|+2}(\tau, u )}d\tau,
\]
where $l'$ is the number of $T$'s in $Z_i^{\a'+1}$.
\end{prop}
Thus, we conclude that:
\begin{prop}\label{lotL2est}
For sufficiently small $\da$ and all $|\a|\leq N_{top}-1$, the $\normu2$ norms of following
\[Z^{\a+1}_iy^j,Z_i^{\a+1}x^j,\lie_{Z_i}^{\a}\chi',\lie_{Z_i}^{\a}\pre {Z_j} {} {\s{\pi}_{ L }},\lie_{Z_i}^{\a}\pre {Z_j} {} {\s{\pi}}, \lie_{Z_i}^{\a}\pre {Z_j} {} {\s{\pi}_{L}},
\]
are bounded by $\sum_{|\a'|\leq|\a|}||\lie_{Z_i}^{\a'}\chi'||_{L^{2}(\Si_{-2}^{\tilde{\da}})}+\sum_{|\a'|\leq|\a|}||Z^{\a'+1}_i\mu||_{L^2(\Si_{-2}^{\tilde{\da}})}+\da^{-l}\int_{-2}^t\sqrt{E_{0,\leq|\a|+2}(t', u )}+\mu_m^{-\frac{1}{2}}(t')
\sqrt{E_{1,\leq |\a|+2}(t', u )}dt',$
where $l$ is the number of $T$'s in $Z_i^{\a}$.
\end{prop}
\section{\textbf{Regularization of the transport equations for $\mu$ and tr$\chi$ and estimates for the top order acoustical terms}}\label{section7}
\hs As discussed in the end of Section 6.2, the key to derive the top order energy estimates is to give precise estimates for $\s{d}R^{\a}_{i}tr\chi$ and $R^{\a'}_{i}T^{l}\s{\de}\mu,$ where $|\a|+1=|\a'|+l+1=N_{top}.$ Note that the right hand side of \eqref{transporttrchi} is of order $2$ while $Ltr\chi$ is also of order 2, which impies that one has to regularize the equation for $tr\chi$. Similarly, the equation for $\mu$ needs to be regularized due to the fact that $m$ and $e$ are both of order 1 while $L\mu$ is of order 1. In this section, we carry out the details for regularizating the transport equation for $tr\chi$ while the computation for $\mu$ will be shown in sketch.
\subsection{\textbf{
Estimates for the top order angular derivatives of tr$\chi$}}
\hs The idea of regularization is moving the highest order term in \eqref{transporttrchi} \textbf{in frames.} Note that \eqref{transporttrchi} can be written as $L tr\chi=\dfrac{1}{2}\dfrac{dH}{dh}\s{\de}h+\text{l.o.ts}$ and the equation $\square_gh:=\p$ (see Lemma\ref{squareh}) allows one to write $\s{\de}h$ as
$\s{\de}h=\mu^{-1}L(\dl h)+\text{l.o.ts}$. Thus, $Ltr\chi=\frac{1}{2}\frac{dH}{dh}\mu^{-1}L(\dl h)+\text{l.o.ts}$ and one can move the term $L(\dl h)$ to the left hand side, which regularizes the transport equation for $tr\chi$.\\
\hs Direct computation leads to the following lemma.
\begin{lem}\label{squareh}
The enthalpy $h$ satisfies the following equation:
\bee\label{boxgh}
\bes
\square_gh&=\p_1+\mu^{-1}L\fai_i\dl\fai_i-\sdfai^2+\Omega^{-1}\dfrac{d\Omega}{dh}
\dfrac{1}{2\mu}\left[Lh\cdot\dl(h-a\fe)+\dl h\cdot L(h-a\fe)\right]\\
&-\Omega^{-1}\dfrac{d\Omega}{dh}\sd(h)\cdot\sd(h-a\fe)+a\square_g\fe:=\p,
\end{split}
\ee
where
\bee
\p_1=-\dfrac{2a\eta'}{\p\eta}\dfrac{\pa\fai_j}{\pa x^j}[\fai_0-(\fai_i)^2]+
\dfrac{a}{\p^2}\left(\dfrac{\pa}{\pa t}-\fai_i\dfrac{\pa}{\pa x^i}\right)[\fai_0-\dfrac{1}{2}(\fai_i)^2].
\ee
\hs Note that the right hand side of \eqref{boxgh} is of order $1$. Moreover, it holds that
\bee
|\mu\p_1|\les(|T\fai_j|+|\sd\fai_j|)|[\fai_0-(\fai_i)^2]|+
|\dl[\fai_0-\dfrac{1}{2}(\fai_i)^2]|.
\ee
\end{lem}
\hs It follows from \eqref{decomposewavenonconformal} and Lemma\ref{squareh} that
\bee\label{musdeh}
\mu\s{\de}h=L(\dl h)+\dfrac{1}{2}(tr\dc\cdot Lh+tr\chi\cdot\dl h)+2\zeta\cdot\sd h+\mu\p.
\ee
and then
\bee\label{mutra}
\bes
\mu tr\a&=-\dfrac{1}{2}\dfrac{dH}{dh}(L(\dl h))
-\dfrac{dH}{dh}\zeta\cdot\s{d}h-\dfrac{1}{2}\dfrac{dH}{dh}\mu\p+\mu tr\a^{[N]}.
\end{split}
\ee
Thus,
\bee\label{mutrxf}
L\left(\mu tr\chi+\breve{f}\right)=2L\mu tr\chi-\mu|\chi|^2+\breve{g},
\ee
where
\bee
\bes
\breve{f}&=-\dfrac{1}{2}\dfrac{dH}{dh}\dl h,\\
\breve{g}&=-\dfrac{1}{2}\dfrac{d^2H}{dh^2}Lh\cdot\dl h+\dfrac{1}{2}\dfrac{dH}
{dh}\zeta\cdot\sd h+\dfrac{1}{2}\dfrac{dH}{dh}\mu\p-\mu tr\a^{[N]}.
\end{split}
\ee
\hs Note that the right hand side of \eqref{mutrxf} does not invole $\mu^{-1}$ and $\breve{g}$ is of order $1$, which regularizes the equation for $Ltr\chi$. Next we will derive the transport equations for high order angular derivatives of $tr\chi$. Set $F_{\a}=\mu\sd R_i^{\a}tr\chi+\sd R_i^{\a}\breve{f}.$ Then, $F_{\a}$ satisfies the following equation:
\bee\label{Fa}
LF_{\a}+(tr\chi-2\mu^{-1}L\mu)F_{\a}=(\dfrac{1}{2}tr\chi-2\mu^{-1}L\mu)\sd R_i^{\a}\breve{f}-\mu\sd R_i^{\a}|\hat{\chi}|^2+g_{\a},
\ee
where $g_{\a}$ is the commutation term given by
\bee\label{expressionga}
\bes
g_{\a}&=\lie_{R_i}^{\a}g_0+\sum_{|\be_1|+|\be_2|=|\a|,|\be_1|>0}
R_i^{\be_1}(\dfrac{1}{2}tr\chi-2\mu^{-1}L\mu)\sd R_i^{\be_2}\breve{f}
-\sum_{|\be_1|+|\be_2|=|\a|,|\be_1|>0}R_i^{\mu}\cdot\sd R_i^{\be_2}|\hat{\chi}|^2\\
&+\sum_{|\be_1|+|\be_2|=|\a|-1}\lie_{R_i}^{\be_1}\cdot
\lie_{\srpi_{L}}\cdot\lie_{R_i}^{\be_2}F_0-\sum_{|\be_1|+|\be_2|=|\a|,|\be_1|>0}
R_i^{\be_1}(tr\chi-2\mu^{-1}L\mu)\lie_{R_i}^{\be_{2}}F_0\\
&-\sum_{|\be_1|+|\be_2|=|\a|,|\be_1|>0}\lie_L(R_i^{\be_1}\mu\cdot\sd R_i^{\be_2}tr\chi)-(tr\chi-2\mu^{-1}L\mu)\sum_{|\be_1|+|\be_2|=|\a|,|\be_1|>0}R_i^{\be_1}\mu
\cdot \sd R_i^{\be_2}tr\chi.
\end{split}
\ee
\hs Since $|\a|=N_{top}-1$, then $F_{\a}$ is of order $N_{top}+1$. The top order terms (order of $N_{top}+1$) on the right hand side of \eqref{Fa} come from
\been
\item $\sd R_i^{\a}\breve{f}$, $\mu\sd R_i^{\a}|\hat{\chi}|^2$,
\item $g_{\a}$, which includes $\lie_{R_i}^{\be_1}\cdot
      \lie_{\srpi_{L}}\cdot\lie_{R_i}^{\be_2}F_0$, $R_i\mu\cdot\sd R_i^{\a-1}Ltr\chi$ and $\lie_{R_i}^{\a}g_0$. In $\lie_{R_i}^{\a}g_0$, the top order terms are:
      \bee\label{topordertermsRiag0}
      \left\{\bes
      &(\sd\mu)R_i^{\a}(Ltr\chi+|\chi|^2),\\
      &tr\chi\cdot \sd R_i^{\a}(2L\mu+\breve{f}),\\
      &\sd R_i^{\a}\breve{g}.
      \end{split}\right.
      \ee
\een
\hs Now we turn to the estimates for $F_{\a}$. For any $S_{t,u}$ $1-$form $\xi$, it holds that
\bee
|\xi|L|\xi|=(\xi,\lie_L\xi)-\xi\cdot\hat{\chi}\cdot\xi-\dfrac{1}{2}tr\chi|\xi|^2.
\ee
\hs Taking $\xi=F_{\a}$ yields
\bee\label{LFa}
\bes
L|F_{\a}|&\leq(2\mu^{-1}L\mu-\dfrac{3}{2}tr\chi+|\hat{\chi}|)|F_{\a}|
+(-\dfrac{1}{2}tr\chi+2\mu^{-1}|L\mu|)|\sd R_i^{\a}\breve{f}|\\
&+\mu|\sd R_i^{\a}|\hat{\chi}|^2|+|g_{\a}|.\\
\end{split}
\ee
Indeed, the term $(2\mu^{-1}L\mu-\dfrac{3}{2}tr\chi+|\hat{\chi}|)|F_{\a}|$ can be absorbed by Gronwall inequality as follows.
\been
\item $(-\dfrac{3}{2}tr\chi+|\hat{\chi}|)|F_{\a}|$ can be bounded by
$C|F_{\a}|$;
\item for the term $\mu^{-1}L\mu$: if $\mu\ge\dfrac{1}{10}$, $|\mu^{-1}L\mu|\les 1$;
      if $\mu<\dfrac{1}{10}$, then by Proposition\ref{keymu1}, $L\mu<0$, which means that this term can be ignored.
\een
Then, applying Gronwall inequality to \eqref{LFa} yields
\bee\label{Fal2}
\bes
||F_{\a}||_{\normu2}&\les ||F_{\a}||_{L^2(\Si_{-2}^{\tilde{\da}})}+\int_{-2}^t||
(\mu^{-1}|L\mu|+\dfrac{1}{2}tr\chi)\sd R_i^{\a}\breve{f}||_{L^2(\Si_{\tau}^{u})}d\tau\\
&+\int_{-2}^t||\mu\sd R_i^{\a}|\hat{\chi}|^2||_{L^2(\Si_{\tau}^{u})}d\tau
+\int_{-2}^t||g_{\a}||_{L^2(\Si_{\tau}^{u})}d\tau\\
&:=||F_{\a}||_{L^2(\Si_{-2}^{\tilde{\da}})}+I_1+I_2+I_3,
\end{split}
\ee
where $I_1,I_2,I_3$ are the integrals given in order respectively, which will be estimated below.
\subsubsection{\textbf{Elliptic estimates on $S_{t,u}$ and estimates for $I_2$ and $I_3$}}
\hs We start with estimating for $I_2$. To this end, rewrite \eqref{Codazzi2} as
\bee
\s{div}\hat{\chi}-\dfrac{1}{2}\sd tr\chi=\be^{\ast}-\mu^{-1}(\zeta\cdot\chi-
\zeta\cdot tr\chi):=i,
\ee
where $i$ is of order $1$. This indicates that one is able to estimate $\sd R_i^{\a}\hat{\chi}$ in terms of $\sd R_i^{\a}tr\chi$ together with some lower order terms. To this end, we state the following two lemmas, whose proofs are given in\cite{christodoulou2014compressible}.
\begin{lem}
Set $\pre {(i_k\cdots i_1)} {} {\hat{\chi}}:=\hat{\lie}_{R_{i_k}}\cdots\hat{\lie}_{R_{i_1}}\hat{\chi}$. Then
$\pre {(i_k\cdots i_1)} {} {\hat{\chi}}$ satisfies the following elliptic system:
\bee\label{sdivhatchi}
\s{div}\left(\pre {(i_k\cdots i_1)} {} {\hat{\chi}}\right)=
\dfrac{1}{2}\s{d}R_{i_k}\cdots R_{i_1}tr\chi+\pre {(i_k\cdots i_1)} {} {i},
\ee
where
\bee
\bes
\pre {(i_k\cdots i_1)} {} {i}&=\left(R_{i_k}+\dfrac{1}{2}tr\pre {R_{i_k}} {} {\pi}\right)\cdots\left(R_{i_1}+\dfrac{1}{2}tr\pre {R_{i_1}} {} {\pi}\right)i\\
&+\sum_{m=0}^{k-1}\left(R_{i_{k}}+\dfrac{1}{2}tr\pre {R_{i_{k}}} {} {\pi}\right)\cdots\left(R_{i_{k-m+1}}+\dfrac{1}{2}tr\pre {R_{i_{k-m+1}}} {}
{\pi}\right)\pre {(i_{k-m}\cdots i_1)} {} {q},\\
\pre {(i_{k+1}\cdots i_1)} {} {q_a}&=\dfrac{1}{4}tr\pre {R_{i_{k+1}}} {} {\pi}\cdot
\sd R_{i_k}\cdots R_{i_1}tr\chi\\
&+\dfrac{1}{2}\pre {R_{i_{k+1}}} {} {\hat{\pi}}^{BC}\left(D_B\pre {i_k\cdots i_1} {} {\hat{\chi}}_{AC}+D_C\pre {i_k\cdots i_1} {} {\hat{\chi}}_{AB}-D_A\pre {i_k\cdots i_1} {} {\hat{\chi}}_{BC}\right)\\
&+\left(\s{div}\pre {R_{i_{k+1}}} {} {\hat{\pi}}^B\cdot\pre {i_k\cdots i_1} {} {\hat{\chi}}_{AB}\right).
\end{split}
\ee
Note that $\pre {(i_k\cdots i_1)} {} {i}$ is of order $k+1$ and when $k=0$, $\pre {(i_k\cdots i_1)} {} {i}=i$.
\end{lem}
\begin{lem}
For sufficiently small $\da$ and any symmetric $S_{t,u}$ trace free $(0,2)$ tensor $\ta$, it holds that
\bee
\int_{S_{t,u}}\mu^2(|\nabla\ta|^2+|\ta|^2)d\mu_{\sg}\les\int_{S_{t,u}}\mu^2|\s{div}
\ta|^2+|\sd\mu|^2|\ta|^2d\mu_{\sg}.
\ee
\end{lem}
Applying this lemma to \eqref{sdivhatchi} yields
\bee
\bes
||\mu\s{D}\pre {(i_k\cdots i_1)} {} {\hat{\chi}}||_{L^2(\Si_t^{u})}&\les
||\sd\mu||_{\supnormda}||\pre {(i_k\cdots i_1)} {} {\hat{\chi}}||_{\normu2}
+||\mu\s{div}\pre {(i_k\cdots i_1)} {} {\hat{\chi}}||_{\normu2}\\
&\les ||\mu \sd R_{i_k}\cdots R_{i_1}tr\chi||_{\normu2}+ (||\sd\mu||_{\supnormda}+||tr\srpi||_{\supnormda})||\pre {(i_k\cdots i_1)} {} {\hat{\chi}}||_{\normu2}\\
&+||\mu R_{i_k}\cdots R_{i_1}\fai_{\a}||_{\normu2}.
\end{split}
\ee
For $I_2$, since
\bee\label{sdriahatchi}
\bes
\mu\sd R_i^{\a}|\hat{\chi}|^2&=\sum\mu\sd\left(\hat{\lie_{R_i}^{\be_1}}\hat{\chi}\cdot\hat{\lie_{R_i}
^{\be_2}}\hat{\chi}\cdot\hat{\lie_{R_i}^{\be_3}}\sg^{-1}\cdot\hat{\lie_{R_i}
^{\be_4}}\sg^{-1}\right)=2\mu\sd\hat{\lie_{R_i}^{\a}}\hat{\chi}\cdot\hat{\chi}+\text{l.o.ts},
\end{split}
\ee
then
\bee
\bes
\int_{-2}^t||\mu\sd R_i^{\a}|\hat{\chi}|^2||_{L^2(\Si_{\tau}^{u})}d\tau&\les
\int_{-2}^t|\hat{\chi}|\cdot||\sd\hat{\lie_{R_i}^{\a}}\hat{\chi}||_{L^2(\Si_{\tau}
^{u})}+\da\sqrt{E_{1,\leq|\a|+2}(\tau,u)}\
d\tau+||\lie_{R_i}^{\a}\chi'||_{L^2(\Si_{-2}^{\tilde{\da}})}\\
&\les\int_{-2}^t\da||F_{\a}||_{L^2(\Si_{\tau}^{u})}+\da
||\sd R_i^{\a}\breve{f}||_{L^2(\Si_{\tau}^{u})}+\da\sqrt{E_{1,\leq|\a|+2}(\tau,u)}
\ d\tau+||\lie_{R_i}^{\a}\chi'||_{L^2(\Si_{-2}^{\tilde{\da}})}.
\end{split}
\ee
It follows from the definition of $\breve{f}$ that
\bee
\bes
\sd R_i^{\a}\breve{f}&=\sd R_i^{\a}\left(-\dfrac{1}{2}\dfrac{dH}{dh}\dl h\right)=-\dfrac{1}{2}\dfrac{dH}{dh}\left(\sd R_i^{\a}\dl \fai_0+\fai_i\sd R_i^{\a}
\dl\fai_i\right)+\text{l.o.ts},
\end{split}
\ee
and then
\bee\label{sdRiaf}
\int_{-2}^t||\sd R_i^{\a}\breve{f}||_{L^2(\Si_{\tau}^{u})}\les
\int_{-2}^t\sqrt{E_{0,\leq |\a|+2}(\tau,u)}\ d\tau.
\ee
Hence,
\bee
I_2\les \da\int_{-2}^t||F_{\a}||_{L^2(\Si_{\tau}^{u})}+
\sqrt{E_{0,\leq|\a|+2}(\tau,u)}+\sqrt{E_{1,\leq|\a|+2}(\tau,u)}\ d\tau.
\ee
\hs For $I_3$, it follows from the discussion \eqref{expressionga}-\eqref{topordertermsRiag0} that one needs only to consider the following cases:
\been
\item $||\lie_{R_i}^{\be_1}\lie_{\srpi_{L}}F_{\be_2}||_{\normu2}\les \da^{\frac{1}{2}}||F_{\a}||_{\normu2}$.
\item Since $Ltr\chi+|\chi|^2=e tr\chi-tr\a'+a\eta^{-1}\hat{T}^i\fai_i tr\chi$, it holds that
      \bee
      \bes
      &|||\sd\mu|\cdot R_i^{\a}(Ltr\chi+|\chi|^2)||_{\normu2}\les \da||\lie_{R_i}^{\a}\chi'||_{L^2(\Si_{-2}^{\tilde{\da}})}+
      \int_{-2}^t\mu_m^{-\frac{1}{2}}\sqrt{E_{1,\leq|\a|+2}(\tau,u)}\ d\tau+
      \mu_m^{-\frac{1}{2}}\sqrt{E_{1,\leq|\a|+2}(t,u)}.
      \end{split}
      \ee
\item Note that $\sd R_i^{\a}\breve{g}$ involves both angular and spatial derivatives of $\fai$. For the angular derivatives, one can simply bound them by the corresponding energies and fluxes. It suffices to bound terms involving spatial derivatives such as: $Lh\cdot\dl h$.
      \bee
      \bes
      \sd R_i^{\a}(Lh\cdot\dl h)&=\sd R_i^{\a}L\fai_0\cdot\dl h+\sd R_i^{\a}\dl\fai_0\cdot Lh+\text{l.o.ts},\\
       \int_{-2}^t||\sd R_i^{\a}L\fai_0\cdot\dl h||_{L^2(\Si_{\tau}^{u})}d\tau&\les \left(\int_0^uF_{1,\leq|\a|+2}(t,u')du'\right)^{\frac{1}{2}},\\
      \int_{-2}^t||\sd R_i^{\a}\dl\fai_0\cdot Lh||_{L^2(\Si_{\tau}^{u})}d \tau
      &\les\da\int_{-2}^t\sqrt{E_{0,\leq|\a|+2}(\tau,u)}\ d\tau.
      \end{split}
      \ee
\een
\hs Collecting the results above yields
     \bee
     \bes
     I_3&\les \da^{\frac{1}{2}}\int_{-2}^t||F_{\a}||_{L^2(\Si_{\tau}^{u})}d\tau+
     \da||\lie_{R_i}^{\a}\chi'||_{L^2(\Si_{-2}^{\tilde{\da}})}+\int_{-2}^t\mu_m^{-\frac{1}{2}}
     \sqrt{E_{1,\leq|\a|+2}(\tau,u)}\ d\tau\\
     &+\da\int_{-2}^t\sqrt{E_{0,\leq|\a|+2}(\tau,u)}\ d\tau+\left(\int_0^uF_{1,\leq|\a|
     +2}(t,u')du'\right)^{\frac{1}{2}}.
     \end{split}
     \ee
\subsubsection{\textbf{Estimate for $I_1$}}
    \hs Finally, we turn to the estimate for $I_1$, which is the most difficult term and we need the crucial lemma (Lemma\ref{crucial}) in this paper. Note that
    \bee\label{I1bound}
    \bes
    I_1 &\les \int_{-2}^t||(\mu^{-1}|L\mu|+1)\sd R_i^{\a}\breve{f}||_{L^2(\Si_{\tau}^{u}\cap W_{shock})}+||(\mu^{-1}|L\mu|+1)\sd R_i^{\a}\breve{f}||_{L^2(\Si_{\tau}^{u}\cap W_{n-s})}d \tau.
    \end{split}
    \ee
\hs Define:
    \beeq
    t_0&=&\inf\{t\in[-2,t^{\ast})|\mu_m(t)<\dfrac{1}{10}\},\\
    M(t)&=&\max_{(u,\ta),(t,u,\ta)\in W_{shock}}|(\mu^{-1} L\mu)_-(t,u,\ta)|,\\
    I_{b}(t)&=&\int_{-2}^t\mu^{-b}(\tau)M(\tau)d\tau.
    \eeq
    Then, the following crucial lemma holds, which will be proved at the end of this subsection.
    \begin{lem}\label{crucial}
   For sufficiently large $b>2$ and all $t\in[t_0,t^{\ast})$, it holds that
    \bee
    I_b(t)\les b ^{-1}\mu_m^{-b}(t).
    \ee
    \end{lem}
    Assume that this lemma holds. We then continue the estimate for $I_1$.
    Applying Lemma\ref{crucial} to \eqref{I1bound} and taking $b=b_{|\a|+2}$ yield
    \bee
    \bes
    I_1&\les I_{\ba2}\sqrt{\wi{E}_{0,|\a|+2}(t,u)}+\int_{-2}^{t_0}||(\mu^{-1}|L\mu|+1)\sd R_i^{\a}\breve{f}||_{L^2(\Si_{\tau}^{u})}d \tau\\
    &\les\dfrac{1}{b_{|\a|+2}}\mu_m^{-\ba2}
    \sqrt{\wi{E}_{0,|\a|+2}(t,u)}+\int_{-2}^t\mu_m^{-\frac{1}{2}}\sqrt{E_{1,\leq|\a|+2}(\tau,u)}d\tau.
    \end{split}
    \ee
    Collectig the estimates for $I_1$, $I_2$ and $I_3$ yields
    \bee\label{Fabound1}
    \bes
    ||F_{\a}||_{\normu2}&\les||F_{\a}||_{L^2(\Si_{-2}^{\tilde{\da}})}+\da||\lie_{R_i}^{\a}\chi'||_{L^2(\Si
    _{-2}^{\tilde{\da}})}+(\da+\da^{\frac{1}{2}})
    \int_{-2}^t||F_{\a}||_{L^2(\Si_{\tau}^{u})}
    d\tau\\
    &+\dfrac{1}{b_{|\a|+2}}\mu_m^{-b_{|\a|+2}}
    \sqrt{\wi{E}_{0,|\a|+2}(t,u)}+\int_{-2}^t\mu_m^{-\ba2}\left(
    \da\sqrt{\wi{E}_{0,\leq|\a|+2}}+\mu_m^{-\frac{1}{2}}
    \sqrt{\wi{E}_{1,\leq|\a|+2}}\right)d\tau\\
    &+\mu_m^{-\ba2}\left(\int_{0}^u\wi{F}_{1,|\a|+2}(t,u')du'\right)^{\frac{1}{2}}.
    \end{split}
    \ee
    Applying Gronwall inequality to \eqref{Fabound1} yields
    \bee
    \bes
    ||F_{\a}||_{\normu2}&\les||F_{\a}||_{L^2(\Si_{-2}^{\tilde{\da}})}+\da||\lie
    _{R_i}^{\a}\chi'||_{L^2(\Si_{-2}^{\tilde{\da}})}+\dfrac{1}{b_{|\a|+2}}\mu_m^{-b_{|\a|+2}}
    \sqrt{\wi{E}_{0,|\a|+2}(t,u)}\\
    &+\int_{-2}^t\mu_m^{-\ba2}\left(\mu_m^{-\frac{1}{2}}(\tau)
    \sqrt{\wi{E}_{1,\leq|\a|+2}(\tau,u)}+\da\sqrt{\wi{E}_{0,\leq|\a|+2}(\tau,u)}\right)\ d\tau\\
    &+\mu_m^{-\ba2}\left(\int_0^u\wi{F}_{1,\leq|\a|+2}(t,u')du'\right)^{\frac{1}{2}}.
    \end{split}
    \ee
    This together with \eqref{sdRiaf} yields
    \bee\label{boundssdRiatrchi}
    \bes
    ||\mu\sd R_i^{\a}tr\chi||_{\normu2}&\les||F_{\a}||_{L^2(\Si_{-2}^{\tilde{\da}})}+\da||\lie
    _{R_i}^{\a}\chi'||_{L^2(\Si_{-2}^{\tilde{\da}})}+\dfrac{1}{b_{|\a|+2}}\mu_m^{-b_{|\a|+2}}
    \sqrt{\wi{E}_{0,|\a|+2}(t,u)}\\
    &+\int_{-2}^t\mu_m^{-\ba2}\left(\mu_m^{-\frac{1}{2}}(\tau)
    \sqrt{\wi{E}_{1,\leq|\a|+2}(\tau,u)}+\da\sqrt{\wi{E}_{0,\leq|\a|+2}(\tau,u)}\right)\ d\tau\\
    &+\mu_m^{-\ba2}\left(\int_0^u\wi{F}_{1,\leq|\a|+2}(t,u')du'\right)^{\frac{1}{2}}.
    \end{split}
    \ee
    It remains to prove Lemma\ref{crucial}. The key step is to bound $\mu$ accurately from both below and above.
    \begin{pf}
     Note that in the shock region, $L\mu<0$. Thus, $-\eta_m(\tau):=\min_{(u,\ta)}L\mu(\tau,u,\ta)$ for $\tau\in[t_0,t^{\ast})$ where $\eta_m(\tau)>0$, which can be achieved at some $(u_m,\ta_m)$. Moreover, for any $t$, one can choose $(u_t,\ta_t)$ such that $\mu_m(t)=\mu(t,u_t,\ta_t)$. \textbf{Fix} $s\in [t_0,t]$. It follows from
    \bee
    \bes
    e^{at}(-t)L\mu(t,u_s,\ta_s)&=e^{as}(-s) L\mu(s,\mu_s,\ta_s)+I\\
    &=e^{as}s\eta_m(s)+\left[e^{as}sL\mu(s,u_m,\ta_m)-e^{as}sL\mu(s,u_s,\ta_s)\right]+I\\
    &:=e^{as}s\eta_m(s)+II\\
    \end{split}
    \ee
    and Proposition\ref{accruatemu1} that $I$ is bounded by $O(\da M^2)$. Furthermore, it follows from the argument of Proposition\ref{accruatemu1} that
    \bee
    \bes
    1+O(\da)&=\mu(-2,u_s,\ta_s)=\mu(s,u_s,\ta_s)+e^{as}sL\mu(s,u_s,\ta_s)\int_s^{-2}\dfrac{1}{e^{a\tau}(-\tau)}+O(\da M^2),\\
    1+O(\da)&=\mu(-2,u_m,\ta_m)=\mu(s,u_m,\ta_m)+e^{as}sL\mu(s,u_m,\ta_m)\int_s^{-2}\dfrac{1}{e^{a\tau}(-\tau)}+O(\da M^2).
    \end{split}
    \ee
    Thus, this shows that $\left|e^{as}s(L\mu(s,u_s,\ta_s)-L\mu(s,u_m,\ta_m))\right|\les\da M^2$, which implies $|II|\les\da M^2$.\\
    \hs Therefore,
    \bee
    \bes
    \mu_m(t)&=\mu_m(s,u_t,\ta_t)+\int_s^tL\mu(\tau,u_t,\ta_t)d\tau\\
    &\geq \mu_m(s,u_s,\ta_s)+\int_s^t\dfrac{e^{as}s}{e^{a\tau}(-\tau)}\eta_m(s)d\tau+\int_s^t
    \dfrac{II}{e^{a\tau}(-\tau)}d\tau\\
    &\geq\mu_m(s)+\left(e^{as}(-s)\eta_m(s)-\dfrac{1}{b}\right)\int_s^t\dfrac{1}{e^{a\tau}\tau}d\tau,
    \end{split}
    \ee
    for $\da M^2<\dfrac{1}{b}$. On the other hand, it holds that
    \bee\label{upper}
    \bes
    \mu_m(t)&\leq \mu_m(s,u_s,\ta_s)+\int_s^tL\mu(\tau,u_s,\ta_s)d\tau\\
    &=\mu_m(s)+\int_s^t\dfrac{e^{as}s}{e^{a\tau}(-\tau)}\eta_m(s)d\tau+\int_s^t\dfrac{II}{
    e^{a\tau}(-\tau)}d\tau\\
    &\leq \mu_m(s)+(e^{as}(-s)\eta_m(s)+\dfrac{1}{b})\int_s^t\dfrac{1}{e^{a\tau}\tau}d\tau.
    \end{split}
    \ee
    Hence, it follows from Proposition\ref{keymu1} that
    \bee
    \begin{split}
    I_b(t)&\leq \int_{-2}^t\left[\mu_m(s)+(e^{as}(-s)\eta_m(s)-\dfrac{1}{b})\int_s^{\tau}\dfrac{1}{e^{ax}x}dx\right]^{-(b+1)}\dfrac{1}{e^{a\tau}(-\tau)}d\tau\\
    &=-\int_{t_1'}^{t_2'}\left[\mu_m(s)+(e^{as}(-s)\eta_m(s)-\dfrac{1}
    {b})t'\right]^{-(b+1)}dt'\\
    &\leq \dfrac{1}{b(e^{as}(-s)\eta_m(s)-\frac{1}{b})}\left[\mu_m(s)+
    (e^{as}(-s)\eta_m(s)-\dfrac{1}{b})\int_s^t\frac{1}{e^{a\tau}\tau}d\tau\right]^{-b}.
    \end{split}
    \ee
    This together with \eqref{upper} yields
    \bee
    \bes
    I_b(t)&\leq \dfrac{1}{a(e^{as}(-s)\eta_m(s)-\frac{1}{b})}\left[\dfrac{\mu_(s)+
    (e^{as}(-s)\eta_m(s)-\dfrac{1}{b})\int_s^t\frac{1}{e^{a\tau}\tau}}{\mu_(s)+
    (e^{as}(-s)\eta_m(s)+\dfrac{1}{b})\int_s^t\frac{1}{e^{a\tau}\tau}}\right]
    ^{-b}\mu_m^{-b}(t)\\
    &\les\dfrac{1}{b}\left(\dfrac{e^{as}(-s)\eta_m(s)-\frac{1}{b}}{e^{as}(-s)\eta_m(s)+\frac{1}
    {b}}\right)^{-b}\mu_m^{-b}(t)\ra \dfrac{1}{b}e^{2e^{as}s\eta_m(s)}\mu_m^{-b}(t)\ (as\ b\to \infty)\\
    &\les \dfrac{1}{b}\mu_m^{-b}(t).
    \end{split}
    \ee
    \end{pf}
    There are some corollaries of Lemma\ref{crucial}.
\begin{cor}
For sufficiently large $b$ and all $t\in[t_0,t^{\ast})$, it holds that
\[\int_{-2}^t\mu_m^{-b-1}dt'\les\dfrac{1}{b}\mu_m^{-b}(t).
\]
\end{cor}
\begin{pf}
The proof is exactly the same as for Lemma\ref{crucial}.
\end{pf}
\begin{cor}
For sufficiently small $\da$ and fixed $b>2$, there exists a constant $C_0$ independent of $b$ and $\da$ such that for all $\tau\in[-2,t]$, the following bound holds:
\[\mu^b_m(t)\leq C_0\mu_m^b(\tau).\]
\end{cor}
\begin{pf}
In Proposition\ref{keymu1}, we have proved that $\mu$ is decreasing along integral curves of $L$ in the shock region.\\
\hs Define: $t_1=\max\{t|\mu_m(\tau)\geq1-\frac{1}{b},\ \text{for}\ \tau\in[-2,t]\}$. Then one can use the same argument to show that $\mu$ is decreasing at very beginning (for $t\geq t_1$) along integral curves of $L$ for sufficiently small $\da$.
\been
\item For $\tau\leq t_1$, it follows from definition of $t_1$ that
      \[
      \mu_m^b(t)\leq 1\leq C_0\left(1-\dfrac{1}{b}\right)^b\leq C_0\mu_m^b(\tau);
      \]
\item while for $\tau\geq t_1$, since $\mu$ is decreasing, then
      \[
      \mu_m^b(t)=\mu_m^b(t,u_t,\ta_t)\leq\mu_m^b(t,u_{\tau},\ta_{\tau})\leq
      C_0\mu_m^b(\tau).
      \]
\een
\end{pf}
    \subsection{\textbf{
    Estimates for the top order spatial derivatives of $\mu$}}
     \hs In this subsection, we derive the estimates for $R_i^{\a'}T^l\s{\de}\mu$ by making use of the basic transport equation \eqref{transportmu}. 
     We start with the following commutation lemma.
    \begin{lem}\label{commutationlde}
    For any function $\fe$, it holds that
    \bee\label{lsdetsde}
    \bes
    [L,\s{\de}]\fe+tr\chi\cdot\s{\de}\fe&=-2\hat{\chi}\cdot\hat{\s{D}}^2\fe-2
    \s{div}\hat{\chi}\cdot\sd\fe,\\
    [T,\s{\de}]\fe+\eta^{-1}\mu tr\ta\cdot\s{\de}\fe&=-2\eta^{-1}\mu\hat{\ta}\cdot\hat{\s{D}}^2\fe
    -2\s{div}(\eta^{-1}\mu\hat{\ta})\cdot\sd\fe.
    \end{split}
    \ee
    \end{lem}
    \begin{pf}
    Since in $(t,u,\ta)$, $L=\dfrac{\pa}{\pa t}$, then
    \bee\label{ldefe}
    \bes
    L\s{\de}\fe&=\sg^{AB}\left(\dfrac{\pa^2}{\pa\ta^A\pa\ta^B}\dfrac{\pa\fe}{\pa t}-\Ga_{AB}^C\dfrac{\pa}{\pa\ta^C}\dfrac{\pa\fe}{\pa t}\right)\\
    &+-2\chi^{AB}\s{D}_A(\sd_B\fe)-\dfrac{\pa}{\pa t}\Gamma_{AB}^C\sg^{AB}\dfrac{\pa\fe}{\pa\ta^C}.
    \end{split}
    \ee
    Note that
    \bee
    \bes
    \dfrac{\pa}{\pa t}\Gamma_{AB}^C&=\s{D}_A\chi_B^C+\s{D}_B\chi_A^C-\s{D}^C\chi_{AB},\\
    \sg^{AB}\dfrac{\pa}{\pa t}\Gamma_{AB}^C&=2div\chi^C-\s{D}^Ctr\chi=2div
    \hat{\chi}^C.
    \end{split}
    \ee
    Substituting this into \eqref{ldefe} yields
    \bee
    L\s{\de}\fe=\s{\de}L\fe-2\chi^{AB}\s{D}_A(\sd_B\fe)-2div\hat{\chi}\cdot
    \sd\fe.
    \ee
    The second identity in\eqref{lsdetsde} can be proved similarly.
    \end{pf}
    \hs Now we regularize the equaiton for $\mu$. Commuting $\s{\de}$ with \eqref{transportmu} yields:
    \bee\label{lmusdemu}
    \bes
    L(\mu\s{\de}\mu)&=(m+2\mu e)\s{\de}\mu+\mu[L,\s{\de}]\mu+\mu\s{\de}m+\mu^2\s{\de}e+2\sd\mu\cdot\sd e.
    \end{split}
    \ee
    Note that
    \bee
    \bes
    \mu\s{\de}m&=L\breve{f}_0+\dfrac{1}{2}tr\chi\cdot \breve{f}_0+\dfrac{1}{2}\dfrac{dH}{dh}\eta^{-2}\mu^2\sd tr\chi\cdot\sd h+m\s{\de}\mu+a\frac{\mu^2}{\eta^2}(\sd x^i)\fai_i\cdot\sd tr\chi+m_0,
    \end{split}
    \ee
    where $\breve{f}_0=\dfrac{1}{2}\dfrac{dH}{dh}T\dl h$. Note that $m_0$ does not contain any acoustical terms of order 2, and
    \bee
    \bes
    \mu^2\s{\de}e&=L(\mu\breve{f}_1)+\dfrac{1}{2}tr\chi\breve{f}_1+\mu\left[\dfrac{1}{\eta}\dfrac{d\eta}{dh}\sd h+
    \dfrac{1}{\eta}\hat{T}^i\sd\fai_i-\dfrac{1}{\eta^2}(\sd h-aT\fe-\eta\hat{T}^i\sd\fai_i)\right]\mu\sd tr\chi+e_0,
    \end{split}
    \ee
    where $\breve{f}_1=\dfrac{1}{\eta}\dfrac{d\eta}{dh}L\dl h+\dfrac{1}{\eta}\hat{T}^iL\dl\fai_i$ and $e_0$ does not contain any acoustical terms of order 2.\\
    \hs Set $\breve{f}=\breve{f}_0+\mu\breve{f}_1$. Then, it follow from Lemma\ref{commutationlde} that
    \bee\label{lf0mu}
    L(\mu\s{\de}\mu)=(2\mu^{-1}L\mu-tr\chi)\mu\s{\de}\mu+L\breve{f}+\dfrac{1}{2}
    tr\chi\cdot\breve{f}-2\mu\hat{\chi}\cdot\hat{\s{D}}^2\mu+\breve{g},
    \ee
    where the principle acoustical part of $\breve{g}$ is given by
    \bee
    \left[\breve{g}\right]_{P.A.}=
    \left(-\sd\mu+\dfrac{2\mu}{\eta}\hat{T}^i\sd\fai_i-\dfrac{2\mu}{\eta^2}\sd h+\dfrac{a\mu}{\eta^2}(T\fe+(\sd x^i)\fai_i)\right)\mu\sd tr\chi.
    \ee
    \hs Define $F_0=\mu\s{\de}\mu-\breve{f}$. Then $F_{\a',l}$ satisfies the following transport equation:
    \bee\label{lfal}
    \bes
    LF_{\a',l}&+(tr\chi-2\mu^{-1}L\mu)F_{\a',l}=(-\dfrac{1}{2}tr\chi+2\mu^{-1}L\mu)
    R_i^{\a'}T^l\breve{f}-2\mu\hat{\chi}\cdot R_i^{\a'}T^l(\hat{\s{D}}^2\mu)+g_{\a',l},\\
    g_{\a',l}&=R_i^{\a'}g_{0,l}+\sum_{k=0}^{\a'-1}R_i^{\a'}\srpi_LF_{\a'-k-1,l}
    +\sum_{k=0}^{\a'-1}R_i^{\a'}y_{\a'-k-1,l},
    \end{split}
    \ee
    where $y_{\a',l}$ is given by
    \bee
    \bes
    y_{\a',l}&=\left[2(R_iL\mu)-\mu R_itr\chi\right]R_i^{\a'}T^l\s{\de}\mu+
    \dfrac{1}{2}R_i(tr\chi)R_i^{\a'}T^l\breve{f}\\
    &-2R_i(\mu\hat{\chi})\cdot R_i^{\a'}T^l(\hat{\s{D}}^2\mu)-4tr\rpi\cdot
    \mu\hat{\chi}\cdot R_i^{\a'}T^l(\hat{\s{D}}^2\mu).
    \end{split}
    \ee
    \hs  
    By the expression of $g_{\a',l}$, one can see that $y_{\a',l}$ is of order $|\a'|+l+2$, and hence $R_i^ky_{\a'-k-1,l}$ is of order $|\a'|+l+1=N_{top}$, which is a lower order term. Thus, the top order acoustical terms in $g_{\a',l}$ come from: $R_i^{\a'}g_{0,l}$ and $R_i^k\srpi_L F_{\a'-k-1,l}$, where the latter one is equalient to $R_i^{\a'+1}T^{l-1}\s{\de}\mu$. In $R_i^{\a'}g_{0,l}$, the top order acoustical terms come from:
    \been
    \item  $R_i^{\a'}T^k\Lambda F_{0,l-k-1}$, which is equivalent to $\mu R_i^{\a'+1}T^{l-1}\s{\de}\mu$.
    \item $R_i^{\a'}T^l\breve{g}$, which is equivalent to $\mu R_i^{\a'}T^l\sd
    tr\chi$, and the following two cases hold:
          \bee\label{discussionga'l}
          \left\{\bes
          \ &\mu\sd R_i^{\a}tr\chi,\hs if\ l=0,\\
          \ &\mu R_i^{\a'+1}T^{l-1}\s{\de}\mu,\hs if\ l\geq 1.
          \end{split}\right.
          \ee
    \een
    It follows from \eqref{lfal} that
    \bee
    \bes
    L|F_{\a',l}|&\leq (2\mu^{-1}L\mu-\dfrac{3}{2}tr\chi+|\hat{\chi}|)|F_{\a',l}|^2+|\dfrac{1}{2}
    tr\chi-2\mu^{-1}L\mu||R_i^{\a'}T^l\breve{f}|\\
    &+2\mu|\hat{\chi}|\cdot|R_i^{\a'}T^l(\hat{\s{D}}^2\mu)|+|g_{\a',l}|.
    \end{split}
    \ee
    Applying Gronwall inequality as in deriving \eqref{Fal2} yields
    \bee
    \bes
    ||F_{\a',l}||_{\normu2}&\leq ||F_{\a',l}||_{L^2(\Si_{-2}^{\tilde{\da}})}+\int_{-2}^t||(\dfrac{1}{2}tr\chi+
    2\mu^{-1}|L\mu|)R_i^{\a'}T^l\breve{f}||_{L^2(\Si_{\tau}^{u})}\\
    &+2\mu|\hat{\chi}|\cdot||R_i^{\a'}T^l\hat{\s{D}}^2\mu||_{L^2(\Si_{\tau}^{u})
    }+||g_{\a',l}||_{L^2(\Si_{\tau}^{u})}d\tau\\
    &=||F_{\a',l}||_{L^2(\Si_{-2}^{\tilde{\da}})}+I_1+I_2+I_3,
    \end{split}
    \ee
    where $I_1,I_2,I_3$ are the integrals given in order respectively, which will be estimated below.
    \subsubsection{\textbf{Elliptic estimates on $S_{t,u}$ and estimates for $I_1$, $I_2$, $I_3$}}
    \hs Note that $I_2$ involves $\hat{\s{D}}^2\mu=\s{D}^2\mu-\dfrac{1}{2}(\s{\de}\mu) g$. The following lemma implies that one could bound $\s{D}^2\mu$ by $\s{\de}\mu$, whose proof is given in\cite{christodoulou2014compressible}. 
    \begin{lem}\label{ellipticmu}
    For any $S_{t,u}$ function $\fe$, it holds that
    \bee
    \int_{S_{t,u}}\mu^2\left(\dfrac{1}{2}|\nabla^2\fe|^2+K|\sd\fe|^2\right)d
    \mu_{\sg}\leq 2\int_{S_{t,u}}\mu^2|\s{\de}\fe|^2d\mu_{\sg}+3\int_{S_{t,u}}
    |\sd\mu|^2\cdot|\sd\fe|^2d\mu_{\sg},
    \ee
    where $K$ is the Gauss curvature of $S_{t,u}$.
    \end{lem}
    \hs Applying this lemma to $I_2$ yields
    \bee\label{riatld2mu}
    \bes
    &||\mu R_i^{\a'}T^l(\hat{\s{D}}^2\mu)||_{\normu2}\leq
    ||\mu R_i^{\a'}T^l\s{\de}\mu||_{\normu2}+||\mu\sd\mu\sd R_i^{\a'}T^l\mu||_{\normu2}+\|\mu[\hat{\s{D}}^2,R_i^{\a'}T^l]\mu\|_{\normu2}\\
    &\leq ||F_{\a',l}||_{\normu2}+|| R_i^{\a'}T^l\breve{f}||_{\normu2}
    +||R_i^{\a'+1}T^l\mu||_{L^2(\Si_{-2}^{\tilde{\da}})}\\
    &+\da^{-l}\int_{-2}^t\sqrt{E_{0,\leq|\a|+2}(\tau,u)}+\mu_m^{-\frac{1}{2}}(\tau)
    \sqrt{E_{1,\leq|\a|+2}(\tau,u)}d\tau.
    \end{split}
    \ee
    \hs It follows from the definition of $\breve{f}$ that
    \bee
    || R_i^{\a'}T^l\breve{f}||_{\normu2}\les\da^{-(l+1)}\left(
    \sqrt{E_{0,\leq|\a|+2}(t,u)}+\da\sqrt{E_{1,\leq|\a|+2}(t,u)}\right).
    \ee
    Therefore,
    \bee
    \bes
    \da^{(l+1)}I_2&\leq\da\int_{-2}^t\da^{l+1}||F_{\a',l}||_{L^2(\Si_{\tau}^{u})}
    +\sqrt{E_{0,\leq|\a|+2}(\tau,u)}+\mu_m^{-\frac{1}{2}}(\tau)\sqrt{E_{1,\leq
    |\a|+2}(\tau,u)}d\tau\\
    &+\da^{l+2}||R_i^{\a'+1}T^l\mu||_{L^2(\Si_{-2}^{\tilde{\da}})}.
    \end{split}
    \ee
    \hs For $I_3$, by the discussion around \eqref{discussionga'l}, the top order acoustical terms in $g_{\a',l}$ are $\mu\sd R_i^{\a}tr\chi$ and $\mu R_i^{\a'+1}T^{l-1}\s{\de}\mu$. The remaining terms in $I_3$ can be bounded directly by the energies and fluxes. Hence, it follows from \eqref{boundssdRiatrchi} that
    \bee
    \bes
    \da^{(l+1)}I_3&\les\da||F_{\a}||_{L^2(\Si_{-2}^{\tilde{\da}})}+\da^2||\lie
    _{R_i}^{\a}\chi'||_{L^2(\Si_{-2}^{\tilde{\da}})}+\da\int_{-2}^t\da^{l+1}||F_{\a',l}||_{L^2(\Si_{\tau}^{u})}\\
    &+\da\int_{-2}^t\mu_m^{-\ba2}\left(\mu_m^{-\frac{1}{2}}(\tau)
    \sqrt{\wi{E}_{1,\leq|\a|+2}}+\da\sqrt{\wi{E}_{0,\leq|\a|+2}}\right)d\tau\\
    &+\da\mu_m^{-\ba2}\left(\int_0^u\wi{F}_{1,\leq|\a|+2}(t,u')du'\right)^{\frac{1}{2}}+\da\int_{-2}^t\mu_m^{-\ba2}\left(
    \sqrt{\wi{E}_{0,\leq|\a|+2}}+\da\mu_m^{-\frac{1}{2}}\sqrt{\wi{E}_{1,\leq|\a|+2}}
    \right)d\tau.
    \end{split}
    \ee
    \hs For $I_1$, the treatment is similar as before. Applying Lemma\ref{crucial} yields
    \bee
    \bes
    I_1&\leq\left(\int_{t_0}^t+\int_{-2}^{t_0}\right)||(\frac{1}{2}tr\chi+2\mu^{-1}|L\mu|)R_i^{\a'}T^l\breve{f}||_{L^2(\Si_{
    \tau}^{u})}d\tau\\
    &\les(\dfrac{1}{\ba2}+1)\mu_m^{-\ba2}\da^{-(l+1)}\left(
    \sqrt{\wi{E}_{0,\leq|\a|+2}(t,u)}+\da\sqrt{\wi{E}_{1,\leq|\a|+2}(t,u)}\right).
    \end{split}
    \ee
    Collecting the estimates for $I_1,I_2,I_3$ yields
    \bee
    \bes
    \da^{l+1}||F_{\a',l}||_{\normda2}&\les\da^{l+1}||F_{\a',l}||_{L^2(\Si_{-2}^{\tilde{\da}})}
    +\da^{l+2}||R_i^{\a'+1}T^l\mu||_{L^2(\Si_{-2}^{\tilde{\da}})}+\da^2||\lie_{R_i}^{\a}\chi'
    ||_{L^2(\Si_{-2}^{\tilde{\da}})}\\
    &+\da\int_{-2}^t\da^{l+1}||F_{\a',l}||_{L^2(\Si_{\tau}^{u})}d\tau+\dfrac{1}{\ba2}\mu_m^{-\ba2}
    \left(\sqrt{\wi{E}_{0,\leq|\a|+2}(t,u)}+\da\sqrt{\wi{E}_{1,\leq|\a|+2}(t,u)}\right)\\
    &+\da\int_{-2}^t\mu_m^{-\ba2}\left(\sqrt{\wi{E}_{0,\leq|\a|+2}}+\mu_m^{
    -\frac{1}{2}}\sqrt{\wi{E}_{1,\leq|\a|+2}}\right)d\tau\\
    &+\da\mu_m^{-\ba2}\left(\int_0^u\wi{F}_{0,\leq|\a|+2}(t,u')du'\right)^{\frac{1}{2}}.
    \end{split}
    \ee
    Applying Gronwall inequality yields
    \bee\label{boundsFa'l}
    \bes
    \da^{l+1}||F_{\a',l}||_{\normu2}&\les\da^{l+1}||F_{\a',l}||_{L^2(\Si_{-2}^{\tilde{\da}})}
    +\da^{l+2}||R_i^{\a'+1}T^l\mu||_{L^2(\Si_{-2}^{\tilde{\da}})}+\da^2||\lie_{R_i}^{\a}\chi'
    ||_{L^2(\Si_{-2}^{\tilde{\da}})}\\
    &+\dfrac{1}{\ba2}\mu_m^{-\ba2}
    \left(\sqrt{\wi{E}_{0,\leq|\a|+2}(t,u)}+\da\sqrt{\wi{E}_{1,\leq|\a|+2}(t,u)}\right)\\
    &+\da\int_{-2}^t\mu_m^{-\ba2}\left(\sqrt{\wi{E}_{0,\leq|\a|+2}}+\mu_m^{
    -\frac{1}{2}}\sqrt{\wi{E}_{1,\leq|\a|+2}}\right)d\tau\\
    &+\da\mu_m^{-\ba2}\left(\int_0^u\wi{F}_{0,\leq|\a|+2}(t,u')du'\right)^{\frac{1}{2}}.
    \end{split}
    \ee
    \hs As a corollary, the same estimate holds for $\mu R_i^{\a'}T^l\s{\de}\mu$ due to
    \bee\label{boundsRia'Tlsdemu}
    \da^{l+1}||\mu R_i^{\a'}T^l\s{\de}\mu||_{\normu2}\leq
    \da^{l+1}||F_{\a',l}||_{\normu2}+\da^{l+1}||R_i^{\a'}T^l\breve{f}||_{\normu2}.
    \ee
    \section{\textbf{Top order energy estimates}}\label{section8}
\subsection{\textbf{Contributions of the top order acoustical terms to the error integrals}}
\hs As we discussed in Section6.2, to complete top order energy estimates, we need to estimate the following error integrals:
\bee\label{topordererrorintegral}
\bes
&-\int_{W_t^u}\til{\p}_{|\a|+2}\cdot K_0R_i^{\a+1}\fai/K_0R_i^{\a'}T^{l+1}\fai
dt'du'd\mu_g,\\
&-\int_{W_t^u}\til{\p}_{|\a|+2}\cdot K_1R_i^{\a+1}\fai/K_1R_i^{\a'}T^{l+1}\fai
dt'du'd\mu_g.
\end{split}
\ee
It has been shown in Section6.2 that the contributions of top order acoustical terms to the error integrals are
\[
T\fai\cdot \sd R_i^{\a}tr\chi,\hs T\fai\cdot R_i^{\a'}T^l\s{\de}\mu.
\]
\subsubsection{\textbf{Contributions associated with $K_0$}}
We first consider the following space-time integral:
      \bee\label{sdRiatrchi}
      \int_{W_t^u}T\fai\cdot\sd R_i^{\a}tr\chi\cdot\dl R_i^{\a+1}\fai\ dt'du'd\mu_{\sg}.
      \ee
      \hs Then,
      \bee\label{K0main}
      \bes
      &\int_{W_t^u}T\fai\cdot\sd R_i^{\a}tr\chi\cdot\dl R_i^{\a+1}\fai\ dt'du'd\mu_g\\
      \leq&\int_{-2}^t\sup(\mu^{-1}|T\fai|)\|\mu\sd R_i^{\a}tr\chi\|_{L^2(\Si_{t'}^{u})}\|\dl R_i^{\a+1}\fai\|_{L^2(\Si_{t'}^{u})}dt'.
      \end{split}
      \ee
      \hs For $\dl R_i^{\a+1}\fai$, it follows from the definition of $E_1$ that
      \bee\label{dlRia+1fai}
      \|\dl R_i^{\a+1}\fai\|_{L^2(\Si_{t}^{u})}\leq\mu_m^{-\ba2}\sqrt{\wi{E}_{0,\leq|\a|+2}}.
      \ee
      \hs For the term $\mu\sd R_i^{\a}tr\chi$, it follows from \eqref{boundssdRiatrchi} that
      \bee\label{K0trx}
      \bes
      \|\mu\sd R_i^{\a}tr\chi\|_{\normu2}&\les\da\|\lie_{R_i}^{\a+1}\chi'\|_{L^2(\Si_{-2}^{\tilde{\da}})}+\sqrt{\wi{E}_{1,\leq|\a|+2}}(t,u)\int_{-2}^t\mu_m^{-\ba2-\frac{1}{2}}(t')\ dt'
      \\
      &+\da\sqrt{\wi{E}_{0,\leq|\a|+2}}(t,u)\int_{-2}^t\mu_m^{-\ba2}(t')\ dt'+\mu_m^{-\ba2}\left(\int_0^u\wi{F}_{1,\leq|\a|+2}(t,u')du'\right)^{\frac{1}{2}}.
      \end{split}
      \ee
      \hs To deal with the integral involving $\mu^{-1}$ above, we have to consider it in the shock region and in the non-shock region. It follows from crucial Lemma\ref{crucial} that for $t\geq t_0$ 
      \bee
      \int_{-2}^{t_0}\mu_m^{-\ba2-1}\leq\mu_m^{-\ba2},\quad \int_{-2}^t\mu_m^{-\ba2-1}\leq\dfrac{C}{\ba2}\mu_m^{-\ba2}.
      \ee
      \hs Hence,
      \bee\label{K0trx2}
      \bes
      \|\mu\sd R_i^{\a}tr\chi\|_{\normu2}&\les(\dfrac{C}{\ba2-\frac{1}{2}}+1)\mu_m^{-\ba2+\frac{1}{2}}
      \sqrt{\wi{E}_{1,\leq|\a|+2}(t,u)}+(\dfrac{C\da}{\ba2-1}+1)\mu_m^{-\ba2+1}\sqrt
      {\wi{E}_{0,\leq|\a|+2}(t,u)}\\
      &+\mu_m^{-\ba2}\left(\int_0^u\wi{F}_{1,\leq|\a|+2}(t,u')du'\right)^{\frac{1}{2}}.
      \end{split}
      \ee
      Here and in the following, we omit the initial data involving $\mu$ and $\chi$. It follows from \eqref{K0main}, \eqref{dlRia+1fai} and \eqref{K0trx2} that for $t\geq t_0$
      \bee\label{K0trchiL2final1}
      \bes
      &\int_{-2}^t\sup(\mu^{-1}|T\fai|)\|\mu\sd R_i^{\a}tr\chi\|_{L^2(\Si_{t'}^{u})}\|\dl R_i^{\a+1}\fai\|_{L^2(\Si_{t'}^{u})}dt'\les\dfrac{CM\da}{\ba2-\frac{1}{2}}\mu_m^{-2\ba2}\wi{E}_{0,\leq
      |\a|+2}(t,u)\\
      &+\dfrac{CM}{\ba2-\frac{1}{2}}\mu^{-2\ba2}\wi{E}_{1,\leq|\a|+2}(t,u)+\dfrac{CM}{2\ba2}\mu_m^{-2\ba2}\int_0^u\wi{F}_{1,\leq|\a|+2}(t,u')du'.
      \end{split}
      \ee
      \hs Later we will choose $\ba2$ suitably large. While for $t\leq t_0$, it holds that
      \bee\label{K0trchiL2final2}
      \bes
      &\int_{-2}^t\sup(\mu^{-1}|T\fai|)\|\mu\sd R_i^{\a}tr\chi\|_{L^2(\Si_{t'}^{u})}\|\dl R_i^{\a+1}\fai\|_{L^2(\Si_{t'}^{u})}dt'\les M\int_{-2}^t\mu_m^{-2\ba2}\left(\da\wi{E}_{0,\leq|\a|+2}+\wi{E}_{1,\leq|\a|+2}
      \right)dt'\\
      &+M\mu_m^{-2\ba2}\int_0^u\wi{F}_{1,\leq|\a|+2}(t,u')du'.
      \end{split}
      \ee
Next we consider the space-time integral:
      \bee\label{Ria'TlsdemuK0}
      \da^{2l+2}\int_{W_t^u}T\fai\cdot R_i^{\a'}T^l\s{\de}\mu\cdot\dl R_i^{\a'} T^{l+1}\fai\ dt'du'd\mu_{\sg},
      \ee
      which can be bounded as
      \bee\label{K0main2}
      \da^{2l+2}\int_{-2}^t\sup(\mu^{-1}|T\fai|)\cdot\|\mu R_i^{\a'}T^l\s{\de}\mu\|_{L^2(\Si_{t'}^{u})}\cdot\|\dl R_i^{\a'}T^{l+1}\fai\|_
      {L^2(\Si_{t'}^{u})}dt'.
      \ee
      \hs It follows from \eqref{boundsFa'l} and \eqref{boundsRia'Tlsdemu} that
      \bee\label{boundsRia'Tlsdemu2}
      \bes
      \da^{l+1}\|\mu R_i^{\a'}T^l\s{\de}\mu\|_{\normu2}&\les(\dfrac{C}{\ba2}\mu_m^{-\ba2}+1)
      \left(\sqrt{\wi{E}_{0,\leq|\a|+2}(t,u)}+\da\sqrt{\wi{E}_{1,\leq|\a|+2}(t,u)}\right)\\
      &+\da\mu_m^{-\ba2}\left(\int_0^u\wi{F}_{1,\leq|\a|+2}(t,u')du'\right)^{
      \frac{1}{2}}.
      \end{split}
      \ee
      Note that $\da^{l+1}\|\dl R_i^{\a'}T^{l+1}\fai\|_{\normu2}\leq
      \mu_m^{-\ba2}\sqrt{\wi{E}_{0,\leq|\a|+2}}$. This together with \eqref{K0main2} and \eqref{boundsRia'Tlsdemu2} implies that for $t\geq t_0$
      \bee\label{K0muL2final1}
      \bes
      &\da^{2l+2}\int_{-2}^t\sup(\mu^{-1}|T\fai|)\cdot\|\mu R_i^{\a'}T^l\s{\de}\mu\|_{L^2(\Si_{t'}^{u})}\cdot\|\dl R_i^{\a'}T^{l+1}\fai\|_
      {L^2(\Si_{t'}^{u})}dt'\\
      &\les\dfrac{CM}{2\ba2-1}\mu_m^{-2\ba2}\wi{E}_{0,\leq|\a|+2}
      (t,u)+\dfrac{CM\da}{2\ba2-1}\mu_m^{-2\ba2}\wi{E}_{1,\leq|\a|+2}(t,u)\\
      &+\dfrac{CM}{2\ba2-1}\mu_m^{-2\ba2}\int_0^u\wi{F}_{1,\leq|\a|+2}(t,u')du',
      \end{split}
      \ee
      and for $t\leq t_0$
      \bee\label{K0muL2final2}
      \bes
      &\da^{2l+2}\int_{-2}^t\sup(\mu^{-1}|T\fai|)\cdot\|\mu R_i^{\a'}T^l\s{\de}\mu\|_{L^2(\Si_{t'}^{u})}\cdot\|\dl R_i^{\a'}T^{l+1}\fai\|_
      {L^2(\Si_{t'}^{u})}dt'\\
      &\les CM\int_{-2}^t\mu_m^{-2\ba2}\left(\wi{E}_{0,\leq|\a|+2}+\da\wi{E}_{1,\leq|\a|+2}
      \right)+C\da M\mu_m^{-2\ba2}\int_0^u\wi{F}_{1,\leq|\a|+2}(t,u')du'.
      \end{split}
      \ee

\subsubsection{\textbf{Contributions associated with $K_1$}}
We need two elementary lemmas of calculus on the manifold.
\begin{lem}\label{change}
Let $f$ and $g$ be arbitrary functions defined on $S_{t,u}$ and $X$ be an $S_{t,u}$ tangential vector field. Then, it holds that
\[\int_{S_{t,u}}f(Xg)d\mu_{\sg}=-\int_{S_{t,u}}\left\{g(Xf)+\dfrac{1}{2}(tr\pre X {} {\s{\pi}})fg\right\}d\mu_{\sg}.
\]
\end{lem}
\begin{pf}
\bee
\int_{S_{t,u}}X(fg)d\mu_{\sg}=\int_{S_{t,u}}div(fgX)d\mu_{\sg}-\int_{S_{t,u}}(div X)fgd\mu_{\sg}.
\ee
\hs By the divergence theorem, the first term on the right hand side vanishes. For the second term, note that
$div X=g^{AB}\s{D}_{A}X_{B}=\dfrac{1}{2}tr\pre {X} {} {\s{\pi}}$, which completes the proof.
\end{pf}
\begin{lem}\label{integralbypartsl}
For any $S_{t,u}$ function $f$, it holds that
    \bee\label{patintegralbyparts}
    \dfrac{\pa}{\pa t}\int_{S_{t,u}}fd\mu_{\sg}=\int_{S_{t,u}}(L+tr\chi)fd\mu_{\sg}.
    \ee
\end{lem}
\begin{pf}
Noting that $L=\dfrac{\pa}{\pa t}|_{t,u,\ta}$ and $\dfrac{\pa}{\pa t}d\mu_{\sg}
=tr\chi d\mu_{\sg}$, one can obtain \eqref{patintegralbyparts} by direct computations.
\end{pf}
 We start with the space-time integral:
      \bee\label{K1main}
      \int_{W_t^u}T\fai\cdot R_i^{\a+1}tr\chi\cdot LR_i^{\a+1}\fai\ dt'du'd\mu_{\sg}.
      \ee
      Note that one could treat \eqref{K1main} with similar argument in estimating the contributions to $K_0$ by using flux to bound $LR_i^{\a+1}\fai$ due to the fact that there is no need to deal with the growth of time. However, we will adapt Christodoulou's another approach which seems to apply to other cases involving long time behavior of the solutions. The idea is a classical idea which changes the above integral into some hypersurface integrals and a controllable space-time integral.\\
      \hs Due to Lemma\ref{integralbypartsl}, the space-time integral \eqref{K1main} can be rewritten as
      \bee
      \bes
      &\int_{\Si_t^u}T\fai\cdot R_i^{\a+1}tr\chi\cdot R_i^{\a+1}\fai-
      \int_{\Si_{-2}^u}T\fai\cdot R_i^{\a+1}tr\chi\cdot R_i^{\a+1}\fai\\
      &-\int_{W_t^u}(L+tr\chi)\left(T\fai\cdot R_i^{\a+1}tr\chi\right)\cdot
      R_i^{\a+1}\fai\ dt'du'd\mu_{\sg}:=I_0-I_1-I_2,
      \end{split}
      \ee
      where $I_1$ can be bounded by the "initial energy" to be checked later.\\
      \hs It follows from Lemma\ref{change} that 
      \bee
      \bes
      I_0
      &=-\int_{\Si_t^u}R_i\left(T\fai\cdot R_i^{\a+1}\fai\right)\cdot R_i^{\a}tr\chi+
      R_i^{\a}tr\chi\cdot T\fai\cdot R_i^{\a+1}\fai\cdot\dfrac{1}{2}tr\srpi\\
      &=-\int_{\Si_t^u}T\fai\cdot R_i^{\a}tr\chi\cdot R_i^{\a+2}\fai-
      \int_{\Si_t^u}R_i^{\a+1}\fai\cdot R_i^{\a}tr\chi\left(
      R_i T\fai+T\fai\cdot\dfrac{1}{2}tr\srpi\right)\\
      &:=-A_1-A_2.
      \end{split}
      \ee
      \hs Since $|tr\srpi|\les\da^{\frac{1}{2}}$, then it suffices to estimate $A_1$. Since the term $R_i^{\a}tr\chi$ in $A_1$ is not a top order term, then it follows from Proposition\ref{Riachi'L2} that
      \bee
      \bes
      |A_1|&\les M\|R_i^{\a}tr\chi\|_{\normu2}\cdot\|R_i^{\a+2}\fai\|_{\normu2}\\
      &\les\dfrac{CM}{\ba2-\frac{1}{2}}\mu_m^{-2\ba2}\wi{E}_{1,\leq|\a|+2}(t,u)
      +\int_{-2}^t\mu_m^{-2\ba2}\wi{E}_{1,\leq|\a|+2}(t',u)dt'.
      \end{split}
      \ee
      \hs For $I_2$, it holds that
      \bee
      \bes
      I_2&=\int_{W_t^u}(L+tr\chi) R_i^{\a+1}tr\chi\cdot T\fai\cdot R_i^{\a+1}\fai
      \ dt'du'd\mu_{\sg}\\
      &+\int_{W_t^u}LT\fai\cdot R_i^{\a+1}tr\chi\cdot R_i^{\a+1}\fai\ dt'du'd\mu_{\sg}:=B_1+B_2.
      \end{split}
      \ee
      \hs For $B_1$, note that
      \bee
      (L+tr\chi) R_i^{\a+1}tr\chi=R_i(L+tr\chi) R_i^{\a}tr\chi-(R_i tr\chi+\srpi_L)(R_i^{\a}tr\chi).
      \ee
      \hs Since $|R_i tr\chi|+|\srpi_L|\les\da$ and $R_i^{\a}tr\chi$ is a lower order term, it suffices to estimate the contribution from $R_i(L+tr\chi) R_i^{\a}tr\chi$. It follows from Lemma\ref{change} that
      \bee
      \bes
      &\int_{W_t^u}R_i(L+tr\chi)R_i^{\a}tr\chi\cdot T\fai\cdot R_i^{\a+1}\fai\ dt'du'd\mu_{\sg}\\
      &=-\int_{W_t^u}(L+tr\chi)R_i^{\a}tr\chi\left
      [R_i(T\fai\cdot R_i^{\a+1}\fai)+T\fai\cdot R_i^{\a+1}\fai\cdot\dfrac{1}{2}tr\srpi\right]\ dt'du'd\mu_{\sg}\\
      &=-\int_{W_t^u}(L+tr\chi)R_i^{\a}tr\chi\cdot R_i^{\a+2}\fai\cdot T\fai+
      (L+tr\chi)R_i^{\a}tr\chi\cdot R_i^{\a+1}\fai\cdot\left(T\fai\cdot\dfrac{1}{2}tr\srpi+R_iT\fai\right)\ dt'du'd\mu_{\sg}\\
      &:=-B_{11}-B_{12}.
      \end{split}
      \ee
      It suffices to estimate $B_{11}$ due to $|tr\srpi|\les\da M$.\\
      \hs Note that
      \bee
      \bes
      Ltr\chi&=etr\chi-|\chi|^2-tr\a'-a\eta^{-1}\hat{T}\fe tr\chi:=\p_0,\\
      LR_i^{\a}tr\chi&=R_i^{\a}\p_0+\sum_{|\be_1|+|\be_2|=|\a|-1}
      R_i^{\be_1}\cdot\srpi_L\cdot R_i^{\be_2}tr\chi,
      \end{split}
      \ee
      where the terms in the summation are lower order terms and thus can be bound by $C\da^{\frac{1}{2}}\mu_m^{-\ba2}(\sqrt{\wi{E}_{0,\leq|\a|+2}}+\sqrt{\wi{E}_{1,\leq|\a|+2}})$ due to Proposition\ref{lotL2est}, so are all the terms $
R_i^{\a}(|\chi|^2),\ R_i^{\be_1}e\cdot R_i^{\be_2}tr\chi,R_i^{\a}tr\chi $, $R_i^{\a}tr\a'$ in $R_i^{\a}\p_0$ \textbf{except the top order term}
     \bee
     R_i^{\a}\s{\de}h.
     \ee
     \hs Corresponding to this term, we estimate the following space-time integral as:
     \bee
     \bes
     &\int_{W_t^u}R_i^{\a}\s{\de}h\cdot T\fai\cdot R_i^{\a+2}\fai\les
     \int_{-2}^tM\|\sd R_i^{\a+1}\fai_{\a}\|^2_{L^2(\Si_{t'}^{u})}dt'\\
     &\les\dfrac{CM}{2\ba2-1}\mu_m^{-2\ba2}\wi{E}_{1,\leq|\a|+2}(t,u)+
     \int_{-2}^tM\mu_m^{-2\ba2}\wi{E}_{1,\leq|\a|+2}(t',u).
     \end{split}
     \ee
     \hs Therefore, 
     \bee
     |B_{11}|\les\dfrac{CM}{2\ba2-1}\mu_m^{-2\ba2}\wi{E}_{1,\leq|\a|+2}(t,u)+
     \int_{-2}^tM\mu_m^{-2\ba2}\wi{E}_{1,\leq|\a|+2}(t,u).
     \ee
     \hs $B_2$ can be bounded directly as follows.
     \bee
     \bes
     |B_2|&\les\dfrac{CM}{2\ba2-\frac{1}{2}}\mu_m^{-\ba2}\wi{E}_{1,\leq|\a|+2}+\dfrac{CM\da}{2\ba2-1}
     \mu_m^{-\ba2}\wi{E}_{0,\leq|\a|+2}+\dfrac{CM}{2\ba2}\mu_m^{-\ba2}\int_0^u\wi{F}_{1,\leq|\a|+2}(t,u')du'.
     \end{split}
     \ee
     \hs Collecting the results above yields
     \bee\label{Ria+1trchiK1final}
     \bes
     \eqref{K1main}&\les \dfrac{CM}{2\ba2-\frac{1}{2}}\mu_m^{-2\ba2}\wi{E}_{1,\leq|\a|+2}+\dfrac{CM\da}{2\ba2-\frac{1}{2}}
     \mu_m^{-\ba2}\wi{E}_{0,\leq|\a|+2}\\
     &+\dfrac{CM}{2\ba2}\mu_m^{-\ba2}\int_0^u\wi{F}_{1,\leq|\a|+2}(t,u')\ du'+\int_{-2}^tM\mu_m^{-2\ba2}\left(\wi{E}_{1,\leq|\a|+2}+\wi{E}_{0,\leq|\a|+2}\right).
     \end{split}
     \ee
Next we consider the space-time integral:
      \bee\label{Ria'TlsdemuK1}
      \da^{2l+2}\int_{W_t^u}T\fai\cdot R_i^{\a'}T^l\s{\de}\mu\cdot LR_i^{\a'}T^{l+1}\fai\ dt'du'd\mu_{\sg},
      \ee
      which can be estimated similarly. Indeed, as for \eqref{K1main}, \eqref{Ria'TlsdemuK1} can be written as
      \bee\label{K1main2}
      \bes
      &\da^{2l+2}\int_{\Si_t^u}T\fai\cdot R_i^{\a'}T^l\s{\de}\mu\cdot R_i^{\a'}T^{l+1}\fai-\da^{2l+2}\int_{\Si_{-2}^u}T\fai\cdot R_i^{\a'}T^l\s{\de}\mu
      \cdot R_i^{\a'}T^{l+1}\fai\\
      -&\da^{2l+2}\int_{W_t^u}(L+tr\chi)
      \left[T\fai\cdot R_i^{\a'}T^l\s{\de}\mu\right]R_i^{\a'}T^{l+1}\fai\ dt'du'd\mu_{\sg}\\
      &:=J_0-J_1-J_2,
      \end{split}
      \ee
      where $J_1$ can be bounded by the "initial energy".\\
      \hs Similarly, it follows from Lemma\ref{change} that
      \bee
      \bes
      J_0&=-\da^{2l+2}\int_{\Si_t^u}T\fai\cdot R_i^{\a'-1}T^l\s{\de}\mu\cdot R_i^{\a'+1}
      T^{l+1}\fai\\
      &-\da^{2l+2}\int_{\Si_t^u}R_i^{\a'-1}T^l\s{\de}\mu\cdot R_i^{\a'}T^{l+1}\fai
      \left(R_iT\fai+T\fai\dfrac{1}{2}tr\srpi\right):=-A_1-A_2,
      \end{split}
      \ee
      where it suffices to estimate $A_1$. Since the term $R_i^{\a'-1}T^l\s{\de}\mu$ is not a top order term, then due to Proposition\ref{lotL2est}, $A_1$ can be bounded as
      \bee
      \bes
      |A_1|&\les M\mu_m^{-\ba2-\frac{1}{2}}\sqrt{\wi{E}_{1,\leq|\a|+2}(t,u)}\cdot
      \int_{-2}^t\da\mu_m^{-\ba2}\left(\mu_m^{-\frac{1}{2}}\sqrt{\wi{E}_{1,\leq|\a|+2}(t',u)}+\sqrt{\wi{E}_{0,\leq|\a|+2}(t',u)}\right) dt'\\
      &\les(\dfrac{CM\da}{2\ba2-\frac{1}{2}}\mu_m^{-2\ba2}+1)\left(\wi{E}_{1,\leq|\a|+2}(t,u)+\wi{E}_{0,\leq|\a|+2}(t,u)\right).
      \end{split}
      \ee
      \hs For $J_2$, it holds that
      \bee
      \bes
      J_2&=\da^{2l+2}\int_{W_t^u}(L+tr\chi)(R_i^{\a'}T^l\s{\de}\mu)\cdot T\fai
      \cdot R_i^{\a'}T^{l+1}\fai\ dt'du'd\mu_{\sg}\\
      &+\da^{2l+2}\int_{W_t^u}LT\fai\cdot R_i^{\a'}T^l\s{\de}\mu\cdot R_i^{\a'}T^{l+1}\fai\ dt'du'd\mu_{\sg}:=B_1+B_2.
      \end{split}
      \ee
      \hs For $B_2$, it can be divided into the integral on $[t_0,t]$, which can be bounded by
      \bee
      \bes
      &\leq\dfrac{CM}{2\ba2}\mu_m^{-2\ba2}\wi{E}_{0,\leq|\a|+2}+\dfrac{CM\da^2}
      {2\ba2}\mu_m^{-2\ba2}\wi{E}_{1,\leq|\a|+2}+CM\da^2 \mu_m^{-2\ba2}\int_0^u\wi{F}_{1,\leq|\a|+2}(t,u')du',
      \end{split}
      \ee
      and on $[-2,t_0]$, which is bounded by
      \bee
      \bes
      &\int_{-2}^tM\mu_m^{-2\ba2}\left(\wi{E}_{0,\leq|\a|+2}(t',u)+\da^2\wi{E}_
      {1,\leq|\a|+2}(t',u)\right)\ dt'+M\da^2\mu_m^{-2\ba2}\int_0^u\wi{F}_{1,\leq|\a|+2}(t,u')du'.
      \end{split}
      \ee
      \hs Note that
      \bee
      (L+tr\chi)R_i^{\a'}T^l\s{\de}\mu=R_i(L+tr\chi)R_i^{\a'-1}T^l\s{\de}\mu-(\srpi_L+R_itr\chi)R_i^{\a'-1}T^l\s{\de}\mu.
      \ee
      \hs Since $|\srpi|+|R_itr\chi|\les\da$, it suffices to consider the contribution from $R_i(L+tr\chi)R_i^{\a'-1}T^l\s{\de}\mu$. It holds that
      \bee
      \bes
      &\da^{2l+2}\int_{W_t^u}R_i(L+tr\chi)R_i^{\a'-1}T^l\s{\de}\mu\cdot T\fai\cdot
      R_i^{\a'}T^{l+1}\fai\ dt'du'd\mu_{\sg}\\
      &=-\da^{2l+2}\int_{W_t^u}(L+tr\chi)R_i^{\a'-1}T^l\s{\de}\mu\cdot
      T\fai\cdot R_i^{\a'+1}T^{l+1}\cdot T\fai\ dt'du'd\mu_{\sg}\\
      &-\da^{2l+2}\int_{W_t^u}(L+tr\chi)R_i^{\a'-1}T^l\s{\de}\mu\cdot R_i^{\a'}T^{l+1}\fai\left(R_iT\fai+T\fai\dfrac{1}{2}tr\srpi\right)\
      dt'du'd\mu_{\sg}\\
      &:=B_{11}+B_{12},
      \end{split}
      \ee
      where it suffices to estimate $B_{11}$.\\
      \hs We decompose $B_{11}$ as:
      \bee
      \bes
      B_{11}&=\da^{2l+2}\int_{W_t^u}LR_i^{\a'-1}T^l\s{\de}\mu\cdot R_i^{\a'+1}
      T^{l+1}\fai\cdot T\fai\ dt'du'd\mu_{\sg}\\
      &+\int_{W_t^u}(tr\chi)R_i^{\a'-1}T^l\s{\de}\mu\cdot R_i^{\a'+1}T^{l+1}\fai\cdot T\fai\ dt'du'd\mu_{\sg}:=B_{111}+B_{112},
      \end{split}
      \ee
      where $B_{112}$ has been estimated before since $R_i^{\a'-1}T^l\s{\de}\mu$ is not a top order term. It remains to estimate $B_{111}$.\\
      \hs Commuting $R_i^{\a'-1}T^l\s{\de}$ with \eqref{transportmu} yields
      \bee
      \bes
      LR_i^{\a'-1}T^l\s{\de}\mu&=\dfrac{1}{2}\dfrac{dH}{dh}R_i^{\a'-1}T^{l+1}\s{\de}h
      +\dfrac{1}{\eta}\dfrac{d\eta}{dh}LR_i^{\a'-1}T^l\s{\de}h+\dfrac{1}{\eta}\hat{T}^iLR_i^{\a'-1}T^l\s{\de}\fai_i+\text{l.o.ts},
      \end{split}
      \ee
      where
      \bee
      \bes
      &\da^{l+1}\|R_i^{\a'-1}T^{l+1}\s{\de}h\|_{\normu2}+\da^{l+1}\|LR_i^{\a'-1}T^l\s{\de}h\|_{\normu2}\les\mu_m^{-\ba2-\frac{1}{2}}
      \sqrt{\wi{E}_{1,\leq|\a|+2}(t,u)}.
      \end{split}
      \ee
      \hs Hence, $B_{111}$ can be bounded as
      \bee
      \bes
      |B_{111}|&\les\int_{-2}^tM
      \left[\mu_m^{-\ba2-\frac{1}{2}}\sqrt{\wi{E}_{1,\leq|\a|+2}}
      +\int_{-2}^{t'}\mu_m^{-\ba2}(\sqrt{\wi{E}_{0,\leq|\a|+2}}+
      \mu_m^{-\frac{1}{2}}\sqrt{\wi{E}_{1,\leq|\a|+2}})\right]\\
      &\cdot \mu_m^{-\ba2-\frac{1}{2}}\sqrt{\wi{E}_{1,\leq|\a|+2}}dt'\\
      &\les\dfrac{CM}{2\ba2}\mu_m^{-2\ba2}(\wi{E}_{1,\leq|\a|+2}+\wi{E}_{0,\leq
      |\a|+2})+\int_{-2}^tM\mu_m^{-2\ba2}(\wi{E}_{1,\leq|\a|+2}+\wi{E}_{0,\leq|\a|+2})
      dt'.
      \end{split}
      \ee
      \hs Consequently,
      \bee\label{Ria'TlsdemuK1final}
      \bes
      \eqref{Ria'TlsdemuK1}&\les\dfrac{CM(1+\da)}{2\ba2}\mu_m^{-2\ba2}\wi{E}_{1,\leq|\a|+2}
      +\dfrac{CM}{2\ba2}\wi{E}_{0,\leq
      |\a|+2}\\
      &+CM\mu_m^{-2\ba2}\int_0^u\wi{F}_{1,\leq|\a|+2}(t,u') du'+\int_{-2}^tM\mu_m^{-2\ba2}(\wi{E}_{1,\leq|\a|+2}+\wi{E}_{0,\leq|\a|+2})
      dt',
      \end{split}
      \ee
where the contributions from the initial data are omitted.
\subsection{\textbf{The top order energy estimates}}
\hs With the estimates for contributions from the top order acoustical terms to \eqref{topordererrorintegral}, we can complete the top order energy estimates associated to $K_0$ and $K_1$.
\subsubsection{\textbf{Estimates associated with $K_1$}}
\hs Let $Z_i\in\{Q,T,R_i\}$. It follows from \eqref{energy1},\eqref{energy2} and Remark\ref{energyremark} that
\bee\label{toporderK11}
\bes
&\sum_{|\a'|\leq|\a|}\da^{2l'}\left(E_1(Z^{\a'+1}_i\fai)(t,u)
+F_1(Z^{\a'+1}_i\fai)(t,u)+K(Z^{\a'+1}_i\fai)(t,u)\right)\\
&\leq C\sum_{|\a'|\leq|\a|}\da^{2l'}(E_1(Z^{\a'+1}_i\fai)+E_0(Z^{\a'+1}_i\fai))(-2,u)+
C\sum_{|\a'|\leq|\a|}\da^{2l'}\int_{W_t^{u}}Q_{1}^{|\a'|+2},
\end{split}
\ee
where $l'$ is the number of $T$'s in $Z_i^{\a'+1}$ and $Q_1^{|\a|+2}:=(-\p_{|\a|+2}K_1(Z_i^{\a+1}\fai)+\frac{1}{2}T^{\a\be}(Z_i^{\a+1}\fai)\pi_{1,\a\be})$. The error integrals $\int_{W_t^{u}}Q_{1}^{|\a|+2}$ contain
\been[i)]
\item contributions from $T^{\a\be}(Z_i^{\a+1}\fai) $ which can be bounded by high order energies and fluxes directly and from the deformation tensors associated to $K_1$, which have been treated in previous sections;
\item major terms of the form
    \bee\label{integralQ1a+2}
    \int_{W^{u}_t}\wi{\p}_{|\a|+2}\cdot K_1 Z^{\a+1}_i\fai\ dt'du' d\mu_{\sg},
    \ee
     where the contributions from the top order acoustical terms have been estimated in last subsection, see\eqref{Ria+1trchiK1final} and \eqref{Ria'TlsdemuK1final}.\\
    \hs It remains to deal with the contributions from lower order terms. For example, for the terms $Z_i^{\be}\chi'$ and $Z_i^{\be+1}\mu$ with $|\be|\leq N_{top}-1$ in $\til{\p}_{|\a|+2}$ (see expression\eqref{recursion2}), it follows from Proposition\ref{lotL2est} that the corresponding integral \eqref{integralQ1a+2} can be bounded as
    \bee
    \bes
    &\dfrac{CM}{2\ba2}\mu_m^{-2\ba2}(t)\left(\wi{E}_{1,\leq|\a|+2}
    +\wi{E}_{0,\leq|\a|+2}+\int_0^u\wi{F}_{1,\leq|\a|+2}\right).
    \end{split}
    \ee
    For the term $\zpi_L\cdot\sd Z_i^{\a}\fai$ with $|\a|=N_{top}$, it follows from the definition of $K(t,u)$ that the corresponding integral \eqref{integralQ1a+2} can be bounded as
    \bee
    \left(\int_0^uF_{1,\leq|\a|+2}\right)
    ^{\frac{1}{2}}\cdot\da^{\frac{1}{2}}\mu_m^{-\ba2}\sqrt{K_{\leq|\a|+2}}\les \mu_m^{-2\ba2}(t)\left(\int_0^u\wi{F}_{1,\leq|\a|+2}+
    C\da^{\frac{1}{2}}K_{\leq|\a|+2}\right),
    \ee
    where $K_{\leq|\a|+2}$ is defined as
    \[
    K_{\leq|\a|+2}(t,u)=\sup_{\tau\in [-2,t]}\{\mu^{ u }_m(\tau)^{2b_{k+1}}\sum_{\fai}\sum_{|\a|=k-1}\da^{2l}K(Z_i^{\a+1}\fai)\},
    \]
     with $l$ being number of $T$'s in $Z_i^{\a+1}$. For the contribution from the term $(Z_{n-1}+\pre {Z_{n-1}} {} {\da})\cdots(Z_{1}+\pre {Z_{1}} {} {\da})\til{\p}_{1}$ in $\til{\p}_{|\a|+2}$, it follows from \eqref{recursion1'} that the corresponding integral \eqref{integralQ1a+2} can be bounded as
    \bee
    \bes
    &\int_{-2}^t\left(\da||Z_i^{\a}T\fai||_{L^2(\Si_{t'}^u)}+\da\mu||Z_i^{\a}\sd\fai||_{L^2(\Si_{t'}^u)}
    +||Z_i^{\a}\dl\fai||_{L^2(\Si_{t'}^u)}\right)\cdot\left(\int_0^uF_{1,\leq|\a|+2}\right)
    ^{\frac{1}{2}} dt'\\
    \leq &\dfrac{CM}{2\ba2}\mu_m^{-2\ba2}(t)\left(\wi{E}_{0,\leq|\a|+2}+\int_0^u\wi{F}_{1,\leq|\a|+2}\right).
    \end{split}
    \ee
    These actually hold for all lower order terms. In conclusion, for the contributions from the lower order terms, the corresponding integral \eqref{integralQ1a+2} is bounded as
    \bee\label{contributionlotK1}
    \dfrac{CM}{2\ba2}\mu_m^{-2\ba2}(t)\left(\wi{E}_{1,\leq|\a|+2}
    +\wi{E}_{0,\leq|\a|+2}+\int_0^u\wi{F}_{1,\leq|\a|+2}\right)+
    C\da^{\frac{1}{2}}\mu_m^{-2\ba2}(t)K_{\leq|\a|+2}.
    \ee
\een
\hs It follows that 
    \bee\label{toporderK12}
    \bes
    &\mu_m^{2\ba2}\sum_{|\a'|\leq||\a|}\da^{2l'}\left(E_1(Z^{\a'+1}_i\fai)(t,u)
     +F_1(Z^{\a'+1}_i\fai)(t,u)+K(Z^{\a'+1}_i\fai)(t,u)\right)\\
     \leq&\dfrac{CM}{2\ba2}\wi{E}_{0,\leq|\a|+2}+\dfrac{CM(1+\da)}{2\ba2}
     \wi{E}_{1,\leq|\a|+2}+CM\int_0^u\wi{F}_{1,\leq|\a|+2}du'\\
     +&CM\int_{-2}^t\wi{E}_{1,\leq|\a|+2}+\wi{E}_{0,\leq|\a|+2}\ dt'+C\da^{\frac{1}{2}}K_{\leq|\a|+2}.
    \end{split}
    \ee
    \hs Note that the right hand side above is increasing in $t$. Thus, \eqref{toporderK12} implies 
    \bee\label{toporderK13}
    \bes
    &\wi{E}_{1,\leq|\a|+2}(t,u)+\wi{F}_{1,\leq|\a|+2}(t,u)+K_{\leq|\a|+2}(t,u)\leq\text{contributions from the initial data}\\
    &+\dfrac{CM}{2\ba2}\wi{E}_{0,\leq|\a|+2}(t,u)+\dfrac{CM(1+\da)}{2\ba2}\wi{E}_
    {1,\leq|\a|+2}(t,u)+C\da^{\frac{1}{2}}K_{\leq|\a|+2}(t,u)\\
    &+CM\left(\dfrac{1}{2\ba2}+1\right)\int_0^u\wi{F}_{1,\leq|\a|+2}\ du'+M\int_{-2}^t
    \wi{E}_{1,\leq|\a|+2}+\wi{E}_{0,\leq|\a|+2}\ dt',
    \end{split}
    \ee
    where the contributions from the initial data are $C\sum_{|\a'|\leq|\a|}\da^{2l'}(E_1(Z_i^{\a'+1}\fai)+E_0(Z_i^{\a'+1}\fai))(-2,u)$+other initial data, which can be omitted, see Remark\ref{otherinitialdata}. Choosing $\ba2$ such that $\dfrac{CM}{\ba2}\leq\dfrac{1}{4}$ and keeping only $\wi{E}_1$ on the left hand side of \eqref{toporderK13} yield for sufficiently small $\da$
    \bee
    \bes
    \wi{E}_{1,\leq|\a|+2}(t,u)&\leq C\sum_{|\a'|\leq|\a|}\da^{2l'}(E_1(Z_i^{\a'+1}\fai)+E_0(Z_i^{\a'+1}\fai))
    (-2,u)+CM\int_{-2}^t\wi{E}_{1,\leq|\a|+2}(t',u)\ dt'\\
    &+\dfrac{CM}{\ba2}\wi{E}_{0,\leq|\a|+2}(t,u)+CM\int_{-2}^t\wi{E}_{0,\leq|\a|+2}(t',u)\ dt'+CM\left(
    \dfrac{1}{2\ba2}+1\right)\int_0^u\wi{F}_{1,\leq|\a|+2}(t,u')\ du'.
    \end{split}
    \ee
    \hs Applying Gronwall inequality yields
    \bee\label{toporderK14}
    \bes
    \wi{E}_{1,\leq|\a|+2}(t,u)&\leq C\sum_{|\a'|\leq|\a|}\da^{2l'}(E_1(Z_i^{\a'+1}\fai)+E_0(Z_i^{\a'+1}\fai))
    (-2,u)+\dfrac{CM}{\ba2}\wi{E}_{0,\leq|\a|+2}(t,u)\\
    &+M\int_{-2}^t\wi{E}_{0,\leq|\a|+2}(t',u)\ dt'+CM\left(
    \dfrac{1}{2\ba2}+1\right)\int_0^u\wi{F}_{1,\leq|\a|+2}(t,u')\ du'.
    \end{split}
    \ee
    \hs Applying the same argument to $\wi{F}_1$ yields
    \bee\label{toporderK1final}
    \bes
    &\wi{E}_{1,\leq|\a|+2}(t,u)+\wi{F}_{1,\leq|\a|+2}(t,u)+K_{\leq|\a|+2}(t,u)\\
    \leq&C\sum_{|\a'|\leq|\a|}\da^{2l'}(E_1(Z_i^{\a'+1}\fai)+E_0(Z_i^{\a'+1}\fai))(-2,u)+C\sum_{
    |\a'|\leq|\a|+1}\da^{2l'}\|\lie_{Z_i}^{\a'}\chi'\|_{L^2(\Si_{-2}^{\tilde{\da}})}
    \\
    &+C\sum_{|\a'|+l'\leq|\a|+2}\da^{2l'}\|Z_i^{\a'}T^l\mu\|_{
    L^2(\Si_{-2}^{\tilde{\da}})}+\dfrac{CM}{2\ba2}\wi{E}_{0,\leq|\a|+2}(t,u)+M\int_{-2}^t\wi{E}_{0,\leq|\a|+2}
    (t',u)\ dt'.
    \end{split}
    \ee
    This completes the top order energy estimates associated to $K_1$.
    \subsubsection{\textbf{Estimates associated with $K_0$}}
    \hs Following similar arguments for \eqref{toporderK11}, one can obtain the following energy inequality for $E_0$ and $F_0$:
\bee\label{toporderK01}
\bes
&\sum_{|\a'|\leq||\a|}\da^{2l'}\left(E_0(Z^{\a'+1}_i\fai)(t, u)+F_0(Z^{\a'+1}_i\fai)(t, u)\right)\\
\leq&C\sum_{|\a'|\leq|\a|}\da^{2l'}(E_0(Z^{\a'+1}_i\fai)+E_1(Z^{\a'+1}_i\fai))(-2, u)+C\sum_{|\a'|\leq|\a|}\da^{2l'}\int_{W_t^{ u}}Q_{0}^{|\a|+2},
\end{split}
\ee
where $l'$ is the number of $T$'s in $Z_i^{\a'+1}$ and $Q_0^{|\a|+2}:=(-\p_{|\a|+2}K_0(Z_i^{\a+1}\fai)+\frac{1}{2}T^{\a\be}(Z_i^{\a+1}\fai)\pi_{0,\a\be})$.
As in the discussion below \eqref{toporderK11}, it can be checked that the error integrals $\int_{W_t^{u}}Q_{0}^{|\a|+2}$ contain
\been[i)]
\item contributions from $T^{\a\be}(Z_i^{\a+1}\fai) $ which can be bounded by high order energies and fluxes directly and from the deformation tensors associated to $K_0$, which have been treated in previous sections;
\item major terms of the form
    \bee\label{integralQ0a+2}
    \int_{W^{u}_t}\wi{\p}_{|\a|+2}\cdot K_0 Z^{\a+1}_i\fai\ dt'du' d\mu_{\sg},
    \ee
     where the contributions from the top order acoustical terms have been estimated in last subsection, see \eqref{K0trchiL2final1},\eqref{K0trchiL2final2},\eqref{K0muL2final1} and \eqref{K0muL2final1}.\\
    \hs For the contributions from lower order terms, the corresponding integral \eqref{integralQ0a+2} can be bounded as due to \eqref{toporderK1final}:
    \bee\label{contributionslotK0}
    \bes
    &CM\int_{-2}^t\left(\sqrt{E_{0,\leq|\a|+1}}+\mu_m^{-\frac{1}{2}}\sqrt
    {E_{1,\leq|\a|+2}}\right)\cdot(\sqrt{E_{0,\leq|\a|+2}})\ dt'+C\da^{\frac{1}{2}}\mu_m^{-2\ba2}K_{\leq|\a|+2}\\
    \leq&\dfrac{CM}{2\ba2}\mu_m^{-2\ba2}\wi{E}_{0,\leq|\a|+2}(t,u)+CM\mu_m^{-2\ba2}
    \int_{-2}^t\wi{E}_{0,\leq|\a|+2}(t',u)\ dt'.
    \end{split}
    \ee
\een
    \hs Substituting \eqref{contributionslotK0} and the estimates estlibalished in the previous subsections into \eqref{toporderK01} yields
    \bee
    \bes
    &\mu_m^{2\ba2}\sum_{|\a'|\leq|\a|}\da^{2l'}\left(E_0(Z^{\a'+1}_i\fai)(t, u)+F_0(Z^{\a'+1}_i\fai)(t, u)\right)\\
\leq&C\sum_{|\a'|\leq|\a|}\da^{2l'}\left(E_0(Z^{\a'+1}_i\fai)(-2, u)+E_1(Z_i^{\a'+1}\fai)(-2,u)\right)\\
    +&\dfrac{CM(1+\da)}{2\ba2}\wi{E}_{0,\leq|\a|+2}(t,u)
    +CM(1+\da)\int_{-2}^t\wi{E}_{0,\leq|\a|+2}(t',u)\ dt'.
    \end{split}
    \ee
    \hs Applying the same argument as before yields
    \bee\label{toporderK03}
    \wi{E}_{0,\leq|\a|+2}(t,u)+\wi{F}_{0,\leq|\a|+2}(t,u)\leq C\sum_{|\a'|\leq|\a|}\da^{2l'}\left(E_0(Z^{\a'+1}_i\fai)(-2, u)+E_1(Z_i^{\a'+1}\fai)(-2,u)\right).
    \ee
    Define a quantity depending on intial data as:
    \bee
    \bes
    \mathcal{D}_{|\a|+2}&=\sum_{|\a'|\leq|\a|}\da^{2l'}\left(E_0(Z^{\a'+1}_i\fai)(-2, u)+E_1(Z_i^{\a'+1}\fai)(-2,u)\right)\\
    &+\sum_{
    |\a'|\leq|\a|+1}\da^{2l'}\|\lie_{Z_i}^{\a'}\chi'\|_{L^2(\Si_{-2}^{\tilde{\da}})}
    +\sum_{|\a'|+l'\leq|\a|+1}\da^{2l'}\|Z_i^{\a'}T^l\mu\|_{
    L^2(\Si_{-2}^{\tilde{\da}})}.
    \end{split}
    \ee
    Then, it follows from \eqref{toporderK03} and \eqref{toporderK1final} that
    \bee\label{toporderen}
    \bes
    \wi{E}_{1,\leq|\a|+2}(t,u)+\wi{F}_{1,\leq|\a|+2}(t,u)+K_{\leq|\a|+2}(t,u)
    &\leq C\mathcal{D}_{|\a|+2},\\
    \wi{E}_{0,\leq|\a|+2}(t,u)+\wi{F}_{0,\leq|\a|+2}(t,u)&\leq
    C\mathcal{D}_{|\a|+2}.
    \end{split}
    \ee
    \begin{remark}\label{otherinitialdata}
    In fact, the derivatives of $\chi'$ and $\mu$ involved in $\mathcal{D}_{|\a|+2}$ can be absorbed by corresponding $E_0$ and $E_1$ due to the fact that on $\Si_{-2}$, $\mu=\eta$ and $\chi'=(1-\eta)\frac{\sg_{AB}}{r}+\eta\sd\fai$.
    \end{remark}
    \section{\textbf{Decent scheme}}\label{section9}
    \hs In Section\ref{section8}, we have established the energy estimates for $N_{top}$ derivatives of $\fai$ \eqref{toporderen} for sufficiently small $\da$. Unfortunately, this is a weighted energy estimate and the weights may go to $0$ when a shock forms due to the definition $\wi{E}_{0,k+1}(t,u)=\sup_{\tau\in [-2,t]}\{\mu_m(\tau)^{2b_{k+1}}E_{0,k+1}(\tau,u)\}$. Thus, it's not enough to recover the bootstrap assumptions. The key point is that to recover the bootstrap assumptions, we don't need the top order energy estimates. \textbf{By choosing the weights suitably and lowering one order of derivatives of $\fai$ in the energies and fluxes, the power of $\mu$ will be decreased one. After several steps, the power of $\mu$ is eliminated while the order of derivatives of $\fai$ are large enough to recover the bootstrap assumptions.} In this section, it will be shown first that the corresponding energy inequality \eqref{toporderen} holds for $N_{top}-1$ order derivatives of $\fai$ and the power of $\mu$ in the weights will be decreased one. Then, after several steps which are far less than $N_{top}$, the power of $\mu$ in the weights will be eliminated.\\
\hs We start with a lemma which refines Corollary 7.2.
\begin{lem}\label{refine}
Let $k$ and $c$ be positive constants with $k>1$. Then, if $\da$ is sufficiently small depending on an upper bound for $c$ and a lower bound for $k$, it holds that, for any $\tau\in[-2,t]$,
\[\mu_m^c(t)\leq k\mu_m^c(\tau).
\]
\end{lem}
\begin{pf}
We use the notations in the proof of Corollary 7.2 and choose $b$ suitably large such that $b>\frac{1}{C\da}$.
\been[(1)]
\item If $\tau\leq t_1$, then $\mu_m(\tau)\geq1-\dfrac{1}{b}$ and it follows that
\[
 \mu_m(\tau)^{-c}\leq(1-\dfrac{1}{b})^{-c}\leq(1-C\da)^{-c}\leq k\leq k\cdot\mu_m^{-c}(t),
\]
for $\da$ sufficiently small (depending on an upper bound for $c$ and a lower bound for $k$).
\item If $\tau\geq t_1$, then
\[\mu_m(t)\leq\mu(t,u_{\tau},\ta_{\tau})\leq\mu(\tau,u_{\tau},\ta_{\tau})=\mu_{m}(\tau)\leq k\mu_m(\tau),
\]
\een
which completes the proof.
\end{pf}
\subsection{\textbf{The next-to-top order energy estimates}}
\hs Recall that the top order is $|\a|+2=N_{top}+1$ and set $\ba2=b_{|\a|+1}+1$. The following computations are similar to the top order energy estimates with $|\a|+2=N_{top}+1$ replaced by $|\a|+1=N_{top}$.\\
\hs Corresponding to \eqref{K1main}, we need to consider the space-time integral:
      \bee\label{nttK11}
      \bes
      &\int_{W_t^u}R_i^{\a}tr\chi\cdot T\fai\cdot LR_i^{\a}\fai\ dt'du'd\mu_{\sg}\\
      \les M&\left(\int_{-2}^t\|R_i^{\a}tr\chi\|^2_{L^2(\Si_{t'}^u)}\ dt'\right)^
      {\frac{1}{2}}\cdot\left(\int_0^u\|LR_i^{\a}\fai\|^2_{L^2(C_{u'}^t)}du'\right)
      ^{\frac{1}{2}}.
      \end{split}
      \ee
      \hs Then, it follows from Proposition\ref{Riachi'L2} that
      \bee
      \bes
      &\int_{W_t^u}R_i^{\a}tr\chi\cdot T\fai\cdot LR_i^{\a}\fai\ dt'du'd\mu_{\sg}\\
      &\les\da M\mu_m^{-2\nba}\wi{E}_{1,\leq|\a|+2}+\da^{-1}M\mu_m^{-2\nba}\int_0^u\wi{F}_{
      1,\leq|\a|+1}\ du'\\
      &\les\da M\mu_m^{-2\nba}\mathcal{D}_{|\a|+2}+\da^{-1}M\mu_m^{-2\nba}\int_0^u\wi{F}_
      {1,\leq|\a|+1}\ du'.
      \end{split}
      \ee
      Here and in the following, we omit the initial data.\\
      \hs Next, Corresponding to \eqref{Ria'TlsdemuK1}, we consider the contribution from $R_i^{\a'-1}T^l\s{\de}\mu$. It follows from Proposition\ref{Zia+1muL2} that
      \bee\label{nttK12}
      \bes
      &\da^{2l+2}\int_{W_t^u}R_i^{\a'-1}T^l\s{\de}\mu\cdot T\fai\cdot LR_i^{\a'-1}
      T^{l+1}\fai\ dt'du'd\mu_{\sg}\\
      &\les M\da^{l+1}\left(\int_{-2}^t\|R_i^{\a'-1}T^l\s{\de}\mu\|^2_{L^2(\Si_{t'}^{u}
      )}\right)^{\frac{1}{2}}\cdot\da^{l+1}\left(\int_0^u\|LR_i^{\a'-1}T^{l+1}\fai\|
      ^2_{L^2(C_{u'}^t)}\right)^{\frac{1}{2}}\\
      &\leq C\da M\mu_m^{-2\nba}\wi{E}_{0,\leq|\a|+2}+C\da M\mu_m^{-2\nba}\wi{E}_{1,\leq|\a|+2}+CM\mu_m^{-2\nba}\int_0^u\wi{F}_{1,\leq|
      \a|+1}\ du'\\
      &\leq C\da M\mathcal{D}_{|\a|+2}+CM\mu_m^{-2\nba}\int_0^u\wi{F}_{1,\leq
      |\a|+1}du'.
      \end{split}
      \ee
\hs Corresponding to \eqref{sdRiatrchi}, we consider the space-time integral:
      \bee\label{nttK01}
      \bes
      &\int_{W_t^u}R_i^{\a}tr\chi\cdot T\fai\cdot\dl R_i^{\a}\fai\ dt'du'd\mu_{\sg},
      \end{split}
      \ee
      which can be bounded by
      \bee
      \bes
      M&\left(\int_{-2}^t\|R_i^{\a}tr\chi\|^2_{L^2(\Si_{t'}^u)}\ dt'\right)
      ^{\frac{1}{2}}\cdot\left(\int_{-2}^t\|\dl R_i^{\a}\fai\|^2_{L^2(\Si_{t'}^{u})}dt'
      \right)^{\frac{1}{2}}\\
      &\les CM\da\mu_m^{-\ba2}\wi{E}_{1,\leq|\a|+2}+CM\da\mu_m^{-\nba}\wi{E}_{0,\leq|\a|+1}\\
      &\les CM\da\mu_m^{-\nba}\mathcal{D}_{|\a|+2}+CM\da\mu_m^{-\nba}\wi{E}_{0,\leq|\a|+1}.
      \end{split}
      \ee
      \hs Next, corresponding to \eqref{Ria'TlsdemuK0}, we consider the space-time integral:
      \bee\label{nttK02}
      \bes
      &\da^{2l+2}\int_{W_t^u}R_i^{\a'-1}T^l\s{\de}\mu\cdot T\fai\cdot\dl R_i^{\a'-1}T^{l+1}\fai\ dt'du'd\mu_{\sg}\\
      \les M&\da^{l+1}\left(\int_{-2}^t\|R_i^{\a'-1}T^l\s{\de}\mu\|^2_{L^2
      (\Si_{t'}^u)}\ dt'\right)^{\frac{1}{2}}\cdot\da^{l+1}\left(
      \int_{-2}^t\|\dl R_i^{\a'-1}T^{l+1}\fai\|_{L^2(\Si_{t'}^u)}\ dt'\right)
      ^{\frac{1}{2}}.
      \end{split}
      \ee
      \hs It follows from the definition of $E_0$ and Proposition\ref{Zia+1muL2} that
      \bee
      \bes
      &\da^{2l+2}\int_{W_t^u}R_i^{\a'-1}T^l\s{\de}\mu\cdot T\fai\cdot\dl R_i^{\a'-1}T^{l+1}\fai\ dt'du'd\mu_{\sg}\\
      &\leq C\da M\mu_m^{-2\nba}\left(\wi{E}_{0,\leq|\a|+2}
      +\wi{E}_{1,\leq|\a|+2}\right)+C\da M\mu_m^{-2\nba}\wi{E}_{0,\leq|\a|+1}\\
      &\leq C\da M\mu_m^{-2\nba}\mathcal{D}_{|\a|+2}+C\da M\mu_m^{-\nba}
      \wi{E}_{0,\leq|\a|+1}.
      \end{split}
      \ee

We then turn to the next-to-top order energy estimates.\\

\hs Similar to \eqref{toporderK11}, it holds that
      \bee\label{energyinequalitynexttotop}
      \bes
      &\mu_m^{2\nba}\sum_{|\a'|\leq|\a|-1}\da^{2l'}\left(E_1(Z_i^{\a'+1}\fai)+F_1(Z_i^{\a'+1}\fai)
      +K(Z_i^{\a'+1}\fai)\right)(t,u)\\
      \leq&C\sum_{|\a'|\leq|\a|-1}\da^{2l'}(E_{1}(Z_i^{\a'+1}\fai)+E_{0}(Z_i^{\a'+1}\fai))(-2,u)
      +C\mathcal{D}_{|\a|+2}+CM\int_0^u\wi{F}_{1,\leq|\a|+1}\ du'\\
      +&C\da M\wi{E}_{0,\leq|\a|+1}+C\int_{-2}^t\wi{E}_{0,\leq|\a|+1}\ dt'+C\da M\wi{E}_{1,\leq|\a|+1}+C\da^{\frac{1}{2}}K_{\leq|\a|+1}.
      \end{split}
      \ee
      \hs Applying the argument as showed in the top order energy estimates part and using Lemma\ref{refine} yield
      \bee
      \bes
      &(\wi{E}_{1,\leq|\a|+1}+\wi{F}_{1,\leq|\a|+1}+K_{\leq|\a|+1})(t,u)\\
      \leq&C\mathcal{D}_{|\a|+2}
      +C\da M\wi{E}_{0,\leq|\a|+1}+C\int_{-2}^t\wi{E}_{0,\leq|\a|+1}\ dt'.
      \end{split}
      \ee
 \hs Similar to \eqref{toporderK01}, the following estimate holds for $E_0$ and $F_0$
      \bee
      \bes
      &\mu_m^{2\nba}\sum_{|\a'|\leq|\a|-1}\da^{2l'}\left(
      E_0(Z_i^{\a'+1}\fai)+F_0(Z_i^{\a'+1}\fai)\right)(t,u)\\
      \leq&C\sum_{|\a'|\leq|\a|-1}\da^{2l'}({E}_1(Z_i^{\a'+1}\fai)+{E}_0(Z_i^{\a'+1}\fai))(-2,u)
      +C\mathcal{D}_{|\a|+2}+C\da M\wi{E}_{0,\leq|\a|+1}\\
      +&C\int_{-2}^t
      \wi{E}_{0,\leq|\a|+1}\ dt'+C\da^{\frac{1}{2}}K_{\leq|\a|+1}(t,u)+C\wi{E}_{1,\leq|\a|+1}+
      C\int_0^u\wi{F}_{1,\leq|\a|+1}\ du'\\
      \leq&C\mathcal{D}_{|\a|+2}+C\da M\wi{E}_{0,\leq|\a|+1}+C\int_{-2}^t
      \wi{E}_{0,\leq|\a|+1}\ dt'.
      \end{split}
      \ee
      \hs Applying the same argument as before yields
      \bee
      \bes
      \left(\wi{E}_{0,\leq|\a|+1}+\wi{F}_{0,\leq|\a|+1}\right)(t,u)\leq
      C\mathcal{D}_{|\a|+2}+CM\int_{-2}^t\wi{E}_{0,\leq|\a|+1}\ dt'.
      \end{split}
      \ee
      \hs Applying Gronwall inequality yields
      \bee\label{nttefinal}
      \bes
      \wi{E}_{0,\leq|\a|+1}(t,u)+\wi{F}_{0,\leq|\a|+1}(t,u)&\leq C\mathcal{D}_{|\a|+2},\\
      \wi{E}_{1,\leq|\a|+1}(t,u)+\wi{F}_{1,\leq|\a|+1}(t,u)+K_{\leq|\a|+1}(t,u)&\leq
      C\mathcal{D}_{|\a|+2}.
      \end{split}
      \ee
\subsection{\textbf{The decent scheme}}
\hs Set
    \bee
    b_{|\a|+2-n}=\ba2-n,\quad b_{|\a|+1-n}=\ba2-1-n,\quad \ba2=[\ba2]+\dfrac{3}{4}.
    \ee
    Let the argument for the next-to-top order estimates be the $1-$st step and $n-$th step be the corresponding energy estimates with $\ba2$ and $\nba$ replaced by $b_{|\a|+2-n}$ and $b_{|\a|+1-n}$, respectively. Then, \textbf{as long as $b_{|\a|+1-n}>0$}, i.e. $n\leq[\nba]$, the $n-$th step can be proceeded exactly in the same way as the $1-$st step.\\
    \hs In the $n-$th step, one has to consider the following integral:
    \bee\label{nthmu}
    \int_{-2}^t\mu_m^{-2b_{|\a|+2-n}+1}\ dt'.
    \ee
    \hs Since $n\leq[\nba]$, then $-2b_{|\a|+2-n}+1\leq-\dfrac{5}{2}$. Thus, the integral \eqref{nthmu} can be bounded as $C\mu_m^{-2b_{|\a|+2-n}+2}=C\mu_m^{-2
    b_{|\a|+1-n}}$ due to Lemma\ref{refine} and then the same argument in the $1-$st step can be applied to the $n-$th step.\\
    \hs Therefore, it holds that
    \bee
    \bes
    \wi{E}_{0,\leq|\a|+1-n}+\wi{F}_{0,\leq|\a|+1-n}&\leq C\mathcal{D}_{|\a|+2},\\
    \wi{E}_{1,\leq|\a|+1-n}+\wi{F}_{1,\leq|\a|+1-n}+K_{\leq|\a|+1-n}&\leq C\mathcal
    {D}_{|\a|+2},
    \end{split}
    \ee
    for $n=0,1,\cdots,[\nba]$.\\
    \hs Up to now, $\wi{E}_{1,\leq|\a|+1-n}=\sup\{\mu_m^{-2b_{|\a|+1-n}}E_{1,\leq|\a|+1-n}\}$ and the power of $\mu$ is still not eliminated. Hence, we consider the next step. Set $n=[\ba2]$ to be the final step. Then $b_{|\a|+2-n}=\dfrac{3}{4}$, $b_{|\a|+1-n}=-\dfrac{1}{4}$, and the following integral needs to be considered:
    \bee
    \int_{-2}^t\mu_m^{-\frac{1}{2}}\ dt'.
    \ee
    It follows from the proof of Lemma\ref{crucial} that
\[
\int_{-2}^t\mu_m^{-\frac{1}{2}}\ dt'\leq\left(\mu_m(s)+(e^{as}(-s)\eta_m(s)-\frac{1}{2})\int_s^t\frac{1}{e^{a\tau}\tau}d\tau\right)^{-\frac{1}{2}}\leq C,
\]
for any \textbf{fixed} $s\in[-2,t]$.\\
\hs Set $b_{|\a|+1-n}=0$. Then, the final step can be proceeded in the same way as the previous steps and the following estimates hold for $n=[\ba2]$:
    \bee\label{desiredenergy}
    \bes
    \wi{E}_{0,\leq|\a|+1-n}+\wi{F}_{0,\leq|\a|+1-n}&\leq C\mathcal{D}_{|\a|+2},\\
    \wi{E}_{1,\leq|\a|+1-n}+\wi{F}_{1,\leq|\a|+1-n}+K_{\leq|\a|+1-n}&\leq C\mathcal
    {D}_{|\a|+2}.
    \end{split}
    \ee
    These are the desired energy estimates since it follows from the definition that
    \[
    \wi{E}_{0,\leq|\a|+1-[\ba2]}=\sup_{t'\in[-2,t]}\{\mu_m^{2b_{|\a|+1-[\ba2]}}(t')
    E_{0,|\a|+1-[\ba2]}(t')\}=\sup_{t'\in[-2,t]}\{E_{0,|\a|+1-[\ba2]}(t')\},
    \]
    so the power of $\mu$ has been eliminated.
    \begin{remark}
    The sequence $\{b_k\}$ is given by
    \bee
    \bes
    &b_1=\cdots=b_{|\a|+1-n}=0,\hs b_{|\a|+2-n}=\dfrac{3}{4},\\
    &b_{|\a|+3-n}=b_{|\a|+2-n}+1,\cdots\hs \ba2=\nba+1=[\ba2]+\dfrac{3}{4},
    \end{split}
    \ee
    for $n=[b_{|\a|+2}]$.
    \end{remark}
\section{\textbf{Recovery of the bootstrap assumptions and completion of the proof}}\label{section10}
\hs We need the following Sobolev type inequality to recover the bootstrap assumptions.
\begin{lem}\label{JA}(c.f. Lemma 17.1 in\cite{christodoulou2014compressible})
There exists a numerical constant $C$ 
such that
\bee
\|\fai\|_{L^{\infty}(S_{t,u})}\leq C\cdot S^{\frac{1}{2}}_{[2]}(\fai),
\ee
where
\bee
S_{[2]}(\fai)=\int_{S_{t,u}}\sum_{i,j}(|\fai|^2+|R_i\fai|^2+|R_iR_j\fai|^2)\ d\mu_{\sg}.
\ee
\end{lem}
\hs Define $S_n(t,u)$ to be
\bee
S_n(t,u)=\sum_{|\a'|\leq n}\int_{S_{t,u}}|\da^{l'}Z_i^{\a'}\fai_{\gamma}|^2\ d\mu_{\sg},
\ee
which is the sum of the integrals on $S_{t,u}$ of square of all the variations of $\fai$ up to order $n$ where $l'$ is the number of $T'$s in $Z_i^{\a'}$.\\ 
\hs It follows from Lemma\ref{ellipticenergy} that
\bee
S_{|\a|-[\ba2]+1}\leq C\tilde{\da}\left(E_{0,\leq|\a|+1-[\ba2]}+E_{1,\leq|\a|+1-[\ba2]}\right),
\ee
for all $(t, u)\in[-2,t^{\ast})\times[0,\tilde{\da}]$. The constant $C$ here depends on $\dfrac{\sqrt{\det \sg(t,u)}}{\sqrt{\det\sg(t,0)}}$ and is a numerical constant. 
It follows from \eqref{desiredenergy} that
\[E_{0,\leq|\a|+1-[\ba2]}+E_{1,\leq|\a|+1-[\ba2]}\leq C\cdot \mathcal{D}_{|\a|+2}.
\]
\hs Now for any variations of $\fai$ up to order $|\a|+1-[\ba2]-2$, it holds that
\[S_2(\fai)\leq S_{|\a|-[\ba2]+1}(t, u).
\]
\hs Thus, by Lemma\ref{JA}, the following estimate holds for all $|\a'|\leq|\a|-1-[\ba2]$:
\bee
\da^{l'}\sup|Z_i^{\a'}\fai|\leq C_0\da\sqrt{\mathcal{D}_{|\a|+2}},
\ee
for all $(t,u)\in[-2,t^{\ast})\times[0,\tilde{\da}]$. The constant $C_0$ here depends on some numerical constants.\\
\hs Hence, by choosing $M$ suitably large such that $C_0\sqrt{\mathcal{D}_{|\a|+2}}<M$, one can recover the bootstrap assumptions.\\
\hs Recall the definitions,
\beeq
&&s_{\ast}=\sup\{t|t\geq-2\ and\ \mu_m^{\tilde{\da}}(t)>0\},\\
&&t_{\ast}=\sup\{\tau|\tau\geq-2\ \text{such that smooth solution exists for all} (t,u)\in[-2,\tau)\times[0,\tilde{\da}]\ and\ \theta\in S^2\},\\
&&s^{\ast}=\min\{s_{\ast},\sigma\}, t^{\ast}=\min\{t_{\ast},s^{\ast}\},
\eeq
where $\sigma$ is a fixed small constant.\\
\hs To finish the proof, it suffices to show $t^{\ast}=s^{\ast}$, i.e. either $t^{\ast}=\sigma$, then the smooth solution exists on $[-2,\sigma]\times[0,\tilde{\da}]\times S^2$, or $t^{\ast}<\sigma$, then at least at one point on $\Si_{t^{\ast}}^{\tilde{\da}}$ such that $\mu_m(t^{\ast})=0$, and a shock forms in finite time.\\
\hs If $t^{\ast}< s^{\ast}$, then \textbf{$\mu_m$ is positive on $\Si_{t^{\ast}}^{\tilde{\da}}$}. Hence, the \textbf{Jacobian $\de$ of the transformation between the acoustical coordinates and the rectangular coordinates never vanishes on $\Si_{t^{\ast}}^{\tilde{\da}}$}, i.e. the transformation between two coordinates is regular on $\Si_{t^{\ast}}^{\tilde{\da}}$. Moreover, in the acoustical coordinates, $\fai$ and its derivatives $\fai_{\gamma}$ are regular on $\Si_{t^{\ast}}^{\tilde{\da}}$ due to the bootstrap assumptions. 
Therefore, in the rectangular coordinates, $\fai$ and its derivatives $\fai_{\gamma}$ are regular on $\Si_{t^{\ast}}^{\tilde{\da}}$ which belongs to the Sobolev space $H^3$. By the standard local well-posedness theory, one can obtain an extension of the solution to some $t_1>t^{\ast}$, which is a contradiction!\\
\bibliographystyle{plain}
\bibliography{ref} 
\end{document}